\documentclass[reqno]{amsart}

\usepackage{amsthm}
\usepackage{amssymb}
\usepackage{amsmath}
\usepackage{color}
\usepackage[T1]{fontenc}
\usepackage{mathrsfs}
\usepackage[colorlinks, linkcolor=blue, citecolor=blue, urlcolor=blue]{hyperref}

\newtheorem {Theorem} {Theorem}
\numberwithin{Theorem}{section}
\newtheorem {Lemma}[Theorem] {Lemma}
\newtheorem {Proposition}[Theorem] {Proposition}

\newcommand{\Hh}{{\mathcal H}}
\newcommand{\Pp}{{\mathcal P}}
\def\C{{\mathbb C}}
\def\R{{\mathbb R}}
\def\pv{\operatorname{p.v.}}

\begin{document}

\title[Finite-time blowup in a 1D MHD model]{Finite-time blowup in a one-dimensional magnetohydrodynamics model}


\author [Mimi Dai]{Mimi Dai}
\address{Department of Mathematics, Statistics, and Computer Science, University of Illinois Chicago, Chicago, IL 60607, USA}
\email{mdai@uic.edu}

\author [Xiaotong (Dawson) Yang]{Xiaotong (Dawson) Yang}
\address{Department of Mathematics, Statistics, and Computer Science, University of Illinois Chicago, Chicago, IL 60607, USA}
\email{xyang212@uic.edu} 

\thanks{The authors acknowledge support from NSF grant DMS--2308208. M. Dai also acknowledges support from the Simons Foundation.}





\begin{abstract}
We prove finite-time blowup of both the velocity and magnetic field in a one-dimensional magnetohydrodynamics (MHD) model starting from smooth odd initial data. The construction combines the dynamical rescaling approach with weighted stability estimates. Both the velocity and magnetic field remain in $L^2$ before the blowup time.

\bigskip


\hspace{0.02cm}Mathematics Subject Classification: 35Q35, 76B03, 76W05. 
\end{abstract}

\maketitle

\tableofcontents

\section{Introduction}

Magnetohydrodynamics (MHD) describes the interaction of an electrically conducting fluid with a magnetic field. In the incompressible ideal regime, with viscosity and magnetic resistivity neglected, the velocity $\mathbf u$ and magnetic field $\mathbf B$ satisfy
\begin{equation}\label{eq-mhd}
\begin{split}
\partial_t\mathbf u+(\mathbf u\cdot\nabla)\mathbf u
 -(\mathbf B\cdot\nabla)\mathbf B+\nabla\Pi&=0,\\
\partial_t\mathbf B+(\mathbf u\cdot\nabla)\mathbf B
 -(\mathbf B\cdot\nabla)\mathbf u&=0,\\
\nabla\cdot\mathbf u=\nabla\cdot\mathbf B&=0
\end{split}
\end{equation}
where $\Pi$ is the total pressure, including the magnetic pressure. The Lorentz force of the magnetic field acts on the fluid motion, while the induction equation transports and stretches the magnetic field. These coupled mechanisms distinguish MHD from the Euler equation, which is recovered from \eqref{eq-mhd} when $\mathbf B=0$.

The criterion of Caflisch, Klapper, and Steele~\cite{CKS1997} characterizes finite-time blowup of a solution $(\mathbf u, \mathbf B)$ of \eqref{eq-mhd} by the condition
\[
\int_0^T\left(\|\nabla\times\mathbf u(t)\|_{L^\infty}
 +\|\nabla\times\mathbf B(t)\|_{L^\infty}\right)\,dt=\infty,
\]
at time $T$. The interplay of transport, stretching, and velocity-magnetic coupling motivates the study of simpler systems that retain these features of \eqref{eq-mhd}.

An especially useful formulation employs the Els\"asser variables
$\mathbf p=\mathbf u+\mathbf B$ and $\mathbf m=\mathbf u-\mathbf B$, which satisfy
\[
\begin{split}
\partial_t\mathbf p+(\mathbf m\cdot\nabla)\mathbf p+\nabla\Pi&=0,\\
\partial_t\mathbf m+(\mathbf p\cdot\nabla)\mathbf m+\nabla\Pi&=0,\\
\nabla\cdot\mathbf p=\nabla\cdot\mathbf m&=0.
\end{split}
\]
We note each Els\"asser variable is advected by the other. Denoting $\boldsymbol\Omega=\nabla\times\mathbf p$ and $\boldsymbol{\tilde\Omega}=\nabla\times\mathbf m$, the Biot--Savart relations recover $\mathbf p$ and $\mathbf m$ from $\boldsymbol\Omega$ and $\boldsymbol{\tilde\Omega}$. This cross-advection structure, together with nonlocal recovery of the velocities, underlies the one-dimensional (1D) model introduced by the first author, Vyas, and Zhang~\cite{DVZ2023}. 

Nonlocal 1D models have a long history in the study of singularity formation for fluid equations. Constantin, Lax, and Majda~\cite{CLM1985} introduced an explicitly solvable model for the Euler equation that isolates vortex stretching and exhibits finite-time singularities. De Gregorio~\cite{DG1990} added a transport term to the Constantin-Lax-Majda model, and Okamoto, Sakajo, and Wunsch~\cite{OSW2008} studied a parameter-dependent family that varies the balance between transport and stretching. This family takes the form
\[
\partial_t\vartheta+a v\partial_x\vartheta
 =\vartheta\mathcal H\vartheta,\qquad v_x=\mathcal H\vartheta,
\]
where $\mathcal H$ is the Hilbert transform. The cases $a=0$ and $a=1$ are the Constantin--Lax--Majda and De Gregorio models, respectively.  Chen, Hou, and Huang~\cite{CHH2021} developed a dynamic-rescaling and nonlinear-stability approach to asymptotically self-similar blowup for the De Gregorio model. Chen~\cite{Chen2020} established self-similar blowup near $a=\frac12$, including the inviscid setting, using an explicit profile and weighted stability estimates. 

Reduced MHD models also reveal mechanisms leading to finite-time blowup scenarios. Dai and Friedlander~\cite{DF2022} studied two dyadic models for ideal MHD and proved finite-time blowup for positive data in the model with unidirectional energy transfer. In \cite{Dai2023} the author obtained blowup results for positive solutions with sufficiently large data in MHD and Hall-MHD shell models for suitable intermittency parameters. In a related but distinct setting, the author of \cite{Dai2025} introduced one-dimensional nonlocal reductions of electron MHD and established well-posedness and singularity-formation results. Electron MHD retains a Hall-type magnetic nonlinearity and is different from the ideal-MHD system considered here. These works illustrate the value of reduced models for understanding nonlinear mechanisms.

For the coupled 1D model studied below, the authors in \cite{DVZ2023} proved local well-posedness and a Beale--Kato--Majda-type continuation criterion in the periodic setting, and global regularity for a transport-only reduction. More recently, Sun~\cite{Sun2025} proved nonlinear stability near a sinusoidal equilibrium for a related periodic model with mixed stretching terms, in a different parameter regime. Here we work on the real line and construct solutions for which both the velocity and the magnetic field become unbounded in finite time.

\subsection{Main results}

Physical and normalized coordinates are denoted by $(x, t)$ and $(y, \tau)$, respectively. We also denote
\[
f(t) = f( \cdot , t), \qquad f(\tau) = f( \cdot , \tau), \qquad
\mathcal H s(0) = (\mathcal H s)(0, \tau),
\]
with $f( \cdot , 0)$ representing the initial profile.

The 1D MHD model on $\R$ introduced in \cite{DVZ2023} with transport parameter $a \in \R$ is given by
\begin{equation} \label{1.1}
\begin{split}
\Omega_t + a m \Omega_x = ( \Hh \Omega ) \tilde \Omega, \\
\tilde \Omega_t + a p \tilde \Omega_x = ( \Hh \tilde \Omega ) \Omega, \\
p_x = \Hh \Omega, \, \, m_x = \Hh \tilde \Omega,
\end{split}
\end{equation}
where the Hilbert transform is defined by
\[
\Hh f(x) = \frac{1}{\pi}  \pv \int_{\R} \frac{f(y)}{x - y} d y.
\]
We impose the parity and normalization conditions
\[
\Omega (- x, t) = - \Omega (x, t), \, \, \tilde \Omega (- x, t) = - \tilde \Omega (x, t), \, \, p(0, t) = m(0, t) = 0.
\]
Define
\[
\omega = \frac{\Omega + \tilde \Omega}{2} \, \, J = \frac{\Omega - \tilde \Omega}{2}, \, \, u = \frac{p + m}{2}, \, \, B = \frac{p - m}{2}.
\]
Then $u_x = \Hh \omega$ and $B_x = \Hh J$, and system~\eqref{1.1} becomes 
\begin{equation} \label{1.2}
\begin{split}
\omega_t + a (u \omega_x - B J_x) = \omega \Hh \omega - J \Hh J, \\
J_t + a (u J_x  - B \omega_x) = \omega \Hh J - J \Hh \omega. 
\end{split}
\end{equation}
For $J = B = 0$, system~\eqref{1.2} reduces to
\[
\omega_t + a u \omega_x = \omega \mathcal H \omega, \qquad u_x = \mathcal H \omega.
\]
We construct smooth odd data and a time $T \in (0, \infty)$ such that
\[
\lim_{t \to T} \| u(t) \|_\infty
= \lim_{t \to T} \| B(t) \|_\infty = \infty, \qquad
u(t), B(t) \in L^2(\R), \qquad 0 \leq t < T.
\]
The construction perturbs a scalar self-similar profile and controls the magnetic component by an adjoint eigenfunctional.

Theorems~\ref{1.4} and~\ref{1.6} treat $a = \frac12$ and $|a - \frac12| < \delta$ for some $\delta > 0$, respectively. Section~2 contains the magnetic linear estimates. Sections~3 and~4 prove the two theorems.

Denote $\Pp f(x) = \int_0^x \Hh f(y) d y$ and define the reference profile and its potential by 
\begin{equation} \label{1.3}
F(y) = - \frac{16 y}{3(1 + y^2)^2}, \, \, U(y) = \frac{8 y}{3(1 + y^2)}, \, \, Q = \Hh F = U'. 
\end{equation}

\begin{Theorem}\label{1.4}
Let $a=\tfrac12$. There exist $\lambda_*\in(-\tfrac{23}{36},-\tfrac12)$,
a real odd function $f\in C_c^\infty(\R)$, and $\varepsilon_0>0$ such that for each
$0<\varepsilon<\varepsilon_0$, the data
\[
(\omega,J)(\cdot,0)=(F,\varepsilon f),\qquad
(u,B)(\cdot,0)=(U,\varepsilon\mathcal P f)
\]
generate a smooth solution $(u(t),B(t))$ of \eqref{1.2} on $[0,T)$ with $T<\infty$. 
Moreover,  $u(t),B(t)\in L^2(\R)$ for $0\le t<T$; and 
for $t$ sufficiently close to $T$,
\begin{equation} \label{1.5}
c(T-t)^{-2/3}\le\|u(t)\|_\infty\le C(T-t)^{-2/3},
\qquad \|B(t)\|_\infty\ge c\varepsilon(T-t)^{-\gamma}
\end{equation}
with $c,C>0$ independent of $\varepsilon$ and $\gamma=\tfrac23+\lambda_*$.
\end{Theorem}

\begin{Theorem} \label{1.6}
There exists $\delta>0$ such that for every $|a-\tfrac12|<\delta$, one may choose odd functions 
$F_a, f_a\in C_c^\infty(\R)$ and $\varepsilon_0(a)>0$ such that for every $0<\varepsilon<\varepsilon_0(a)$, the data
\[
(\omega,J)(\cdot,0)=(F_a,\varepsilon f_a),\qquad
(u,B)(\cdot,0)=(\mathcal P F_a,\varepsilon\mathcal P f_a)
\]
generate a smooth solution $(u(t),B(t))\in L^2(\R)\times L^2(\R)$ of \eqref{1.1} on $[0,T)$ with $T<\infty$.
Moreover, there exist $\kappa_a, \gamma_a>0$ such that for $t$ sufficiently close to $T$
\begin{equation} \label{1.7}
c_a(T-t)^{-\kappa_a}\le\|u(t)\|_\infty\le C_a(T-t)^{-\kappa_a},
\qquad \|B(t)\|_\infty\ge c_a\varepsilon(T-t)^{-\gamma_a}
\end{equation}
with $c_a,C_a>0$ independent of $\varepsilon$. 
\end{Theorem}

\subsection{Dynamical rescaling}
We first fix $a = \frac{1}{2}$ to prove Theorem~\ref{1.4}. Let $A, L > 0$ and define the rescaled fields by
\begin{equation} \label{1.8}
\widetilde{\omega} (y, \tau) = A (\tau) \omega \left ( L(\tau) y, t(\tau) \right), \, \, \widetilde{J} (y, \tau) = A (\tau) J \left ( L(\tau) y, t(\tau) \right), 
\end{equation}
where
\begin{equation} \label{1.9}
t_\tau = A, \, \, A_\tau = c_A A, \, \, L_\tau = - c_\ell L, \, \, A(0) = L(0) = 1, \, \, t(0) = 0, 
\end{equation}
and the corresponding potentials satisfy
\[
(\widetilde u, \widetilde B) (y, \tau) = \frac{A}{L} (u, B) (L y, t).
\]
We drop the tildes in the normalized evolution. Physical fields will be identified explicitly in the proofs of local existence, where tildes are restored for the normalized fields. The normalized equations are 
\begin{equation} \label{1.10}
\begin{split}
\omega_\tau + \left ( c_\ell y + \frac{u}{2} \right) \omega_y = c_A \omega + \omega \mathcal H \omega - J \mathcal H J + \frac{1}{2} B J_y, \\
J_\tau + \left ( c_\ell y + \frac{u}{2} \right) J_y = c_A J + \omega \mathcal H J - J \mathcal H \omega + \frac{1}{2} B \omega_y.
\end{split}
\end{equation} 
For $(c_A, c_\ell) = ( - 1, \frac{1}{3})$, the pair $(F, 0)$ is stationary. Define
\begin{equation} \label{1.11}
V(y) = \frac{y}{3} + \frac{U(y)}{2} = \frac{y(y^2 + 5)}{3(1 + y^2)}, \qquad
 F'(0) = - \frac{16}{3}.
\end{equation}
Write $\omega = F + s$ and $u = U + \mathcal P s$, and choose the modulation rates by
\begin{equation} \label{1.12}
q_J = \frac{\mathcal H J(0) J'(0)}{F'(0)}, \, \, c_A = - 1 - \mathcal H s(0) + q_J, \, \, c_\ell = \frac{1}{3} - \frac{1}{2} \mathcal H s(0) + \frac{1}{2} q_J. 
\end{equation}
Differentiating \eqref{1.10} at $0$ yields
\begin{equation} \label{1.13}
\partial_\tau s'(0) = \frac{q_J}{2} s'(0),
\end{equation}
and
\[
\omega (0, \tau) = J(0, \tau) = u(0, \tau) = B(0, \tau) = 0,
\qquad u_y(0, \tau) = \frac83 + \mathcal H s(0).
\]
For $s_y(0, 0) = 0$, integration of equation~\eqref{1.13} shows that
\[
s_y(0, \tau) = s_y(0, 0)
\exp \left (\frac12 \int_0^\tau q_J(r) \, d r \right) = 0.
\]
Define the scalar and magnetic linearized operators by
\begin{equation} \label{1.14}
\mathcal L_s s = - V s' + (- 1 + Q ) s + F \mathcal H s - \frac{1}{2} F' \mathcal P s + \mathcal H s(0) \left (- F + \frac{1}{2} y F' \right), 
\end{equation}
\begin{equation} \label{1.15}
\mathcal L h = - V h' - (1 + Q ) h + F \mathcal H h + \frac{1}{2} F' \mathcal P h .
\end{equation} 
The scalar perturbation equation takes the form
\begin{equation} \label{1.16}
\begin{split}
s_\tau & = \mathcal L_s s - \frac{1}{2} \left (\mathcal P s - \mathcal H s(0) y \right) s' + \left ( \mathcal H s - \mathcal H s(0) \right) s  
 + q_J \left (  s - \frac{1}{2} y s' \right) + G_J,
\end{split}
\end{equation}
where
\begin{equation} \label{1.17}
G_J = - J \mathcal H J + \frac{1}{2} B J' + q_J \left (F - \frac{1}{2} y F' \right).
\end{equation}
The magnetic perturbation equation is
\begin{equation} \label{1.18}
J_\tau = \mathcal L J + \mathcal N,
\end{equation}
\begin{equation} \label{1.19}
\mathcal N = - \left (\delta_\ell y + \frac{\mathcal P s}{2} \right) J' + \delta_A J + s \mathcal H J - J \mathcal H s + \frac{1}{2} B s',
\end{equation}
where $\delta_A = c_A + 1$ and $\delta_\ell = c_\ell - \frac{1}{3}$.

Throughout, $\langle y \rangle = (1 + y^2)^{1/2}$ and $\| \cdot \|_p$ denotes the $L^p(\R)$ norm unless a domain is specified.

\begin{Lemma}\label{1.20}
Assume that $\| \langle y \rangle^\eta g \|_2 < \infty$. Then
\[
\| \langle y \rangle^\eta \mathcal H g \|_2
\leq C_\eta \| \langle y \rangle^\eta g \|_2
\]
in each of the following cases.
\begin{enumerate}
\item $g$ is odd and $- \frac12 < \eta < \frac32$.
\item $g$ is even, $\int g = 0$, and $\frac12 < \eta < \frac52$.
\item $g$ is odd, $\int y g(y) \, d y = 0$, and $\frac32 < \eta < \frac72$.
\end{enumerate}
In the second and third cases, respectively, $g$ and $y g$ are integrable.
\end{Lemma}

\begin{proof}
For $m = 1, 2$ and $\eta > m - \frac12$, we have
\[
\| y^{m - 1} g \|_1
\leq \| \langle y \rangle^\eta g \|_2
\| y^{m - 1} \langle y \rangle^{- \eta} \|_2 < \infty.
\]
Parity and the moment cancellations reduce the transform, for $x > 0$, to
\[
\mathcal H g(x) = \frac{2}{\pi x^m} \pv \int_0^\infty
\frac{y^{m + 1} g(y)}{x^2 - y^2} \, d y,
\qquad m = 0, 1, 2
\]
in the three respective cases. Away from $y = x$, the two tails are bounded by constant multiples of
\[
A_m g(x) = x^{- m - 2} \int_0^{x/2} y^{m + 1}|g(y)| \, d y,
\qquad
B_m g(x) = x^{- m} \int_{2 x}^\infty y^{m - 1}|g(y)| \, d y.
\]
With $G(s) = e^{(\eta + 1/2) s}|g(e^s)|$, logarithmic coordinates turn the tail operators into convolutions.
\[
\begin{split}
e^{(\eta + 1/2) s} A_m g(e^s)
& = \int_{\log 2}^\infty e^{- (m + 3/2 - \eta) t} G(s - t) \, d t, \\
e^{(\eta + 1/2) s} B_m g(e^s)
& = \int_{\log 2}^\infty e^{- (\eta - m + 1/2) t} G(s + t) \, d t.
\end{split}
\]
By Young's inequality,
\[
\| x^\eta (A_m g + B_m g) \|_{L^2(0, \infty)}
\leq C_{m, \eta} \| x^\eta g \|_{L^2(0, \infty)},
\qquad m - \frac12 < \eta < m + \frac32.
\]
Apply this estimate to $g \mathbf1_{[1/2, \infty)}$. The remaining contribution satisfies
\[
\left | \frac{2}{\pi x^m} \int_0^{1/2}
\frac{y^{m + 1} g(y)}{x^2 - y^2} \, d y \right |
\leq C_m x^{- m - 2} \| g \|_{L^2(0, 1/2)}, \qquad x \geq 2,
\]
and
\[
\int_2^\infty x^{2 \eta - 2 m - 4} \, d x < \infty,
\qquad \eta < m + \frac32.
\]

For the principal-value part, let $I_k = (2^k, 2^{k + 1})$ and $I_k^* = (2^{k - 1}, 2^{k + 2})$. On $I_k \times I_k^*$,
\[
\frac{2 y^{m + 1}}{x^m(x + y)(x - y)}
= \frac1{x - y} + O(2^{- k}).
\]
Define, for $x \in I_k$,
\[
T_{m, k} g(x) = \frac{2}{\pi x^m} \pv \int_{I_k^*}
\frac{y^{m + 1} g(y)}{x^2 - y^2} \, d y.
\]
The Hilbert bound and the remainder estimate give
\[
\| T_{m, k} g \|_{L^2(I_k)} \leq C_m \| g \|_{L^2(I_k^*)}.
\]
Since $2^{- 2| \eta |} 2^{k \eta} \leq x^\eta \leq 2^{2| \eta |} 2^{k \eta}$ on $I_k^*$,
\[
\sum_{k \geq 0} \| x^\eta T_{m, k} g \|_{L^2(I_k)}^2
\leq C_{m, \eta} \sum_{k \geq 0} \| x^\eta g \|_{L^2(I_k^*)}^2
\leq C_{m, \eta} \| \langle x \rangle^\eta g \|_2^2.
\]

For $g_0 = g \mathbf1_{\{|y| < 4 \}}$ and $g_\infty = g - g_0$,
\[
\begin{aligned}
\| \mathcal H g_0 \|_{L^2(- 2, 2)}
& \leq \| g_0 \|_2 \leq C_\eta \| \langle y \rangle^\eta g \|_2, \\
\sup_{|x| < 2}| \mathcal H g_\infty (x)|
& \leq C \| \langle y \rangle^\eta g \|_2
\left (\int_{|y| \geq 4} \frac{\langle y \rangle^{- 2 \eta}}{y^2} \, d y \right)^{1/2}
\leq C_\eta \| \langle y \rangle^\eta g \|_2, \qquad \eta > - \frac12.
\end{aligned}
\]
Combining these bounds proves the estimate. Weighted density preserves the stated parity and moments.
\end{proof}

We use the weights and space
\[
\sigma = \frac{7}{12}, \, \, w = \langle y \rangle^{2 \sigma}, \, \, \phi = \frac{(1 + y^2)^3}{y^4}, \, \, \psi = (1 + y^2)^2, \, \, X = L^2_{\text{odd}} (w d y).
\]

\begin{Lemma} \label{1.21}
Let $f \in X$, and let $\mathcal P f$ be its normalized potential, so that $(\mathcal P f)' = \mathcal H f$ and $\mathcal P f(0) = 0$. Then
\begin{equation} \label{1.22}
\| \mathcal P f \|_\infty + \sup_{|y| \geq 1} |y |^{\frac{1}{12}} | \mathcal P f(y) | \leq C \|  f \|_X.
\end{equation}
If, in addition, $\| f^{(4)} \|_{L^2(\psi)} < \infty$, then 
\begin{equation} \label{1.23}
\|   \langle  y \rangle^{\mathfrak a_j} f^{(j)} \|_2 \leq C \left ( \| f \|_X + \| f^{(4)} \|_{L^2 (\psi)} \right), \, \, \text{where} \, \, \mathfrak a_j = \frac{28 + 17 j}{48}, \, \, 0 \leq j \leq 4, 
\end{equation}
\begin{equation} \label{1.24}
\|   \langle  y \rangle^{\mathfrak b_j} f^{(j)} \|_\infty \leq C \left ( \| f \|_X + \| f^{(4)} \|_{L^2 (\psi)} \right), \, \, \text{where} \, \, \mathfrak b_j = \frac{73 + 34 j}{96}, \, \, 0 \leq j \leq 3. 
\end{equation}
The weighted derivative estimates \eqref{1.23} and \eqref{1.24} also hold with $f^{(j)}$ replaced by $\mathcal H f^{(j)}$, with the same right-hand side. 
\end{Lemma}

\begin{proof}
\emph{Potential bound.} By Cauchy--Schwarz and Lemma~\ref{1.20},
\[
\| f \|_1 + \| \mathcal H f \|_1
\leq \| \langle y \rangle^{- \sigma} \|_2
\bigl (\| f \|_X + \| \langle y \rangle^\sigma \mathcal H f \|_2 \bigr)
\leq C \| f \|_X.
\]
For the convention $\widehat f(\xi) = \int e^{- i y \xi} f(y) \, d y$, both Fourier transforms are continuous at zero. Oddness and the Hilbert multiplier give
\[
\widehat f(0) = 0, \qquad
\int \mathcal H f
= \lim_{\xi \to 0} \bigl [- i \operatorname{sgn}(\xi) \widehat f(\xi) \bigr] = 0.
\]
Since $\mathcal H f$ is even, for $y > 0$,
\[
\mathcal P f(y) = - \int_y^\infty \mathcal H f(x) d x. 
\]
Consequently,
\[
\begin{aligned}
| \mathcal P f(y)|  & \leq \| \mathcal H f \|_1 \leq C \| f \|_X, \qquad y > 0, \\
| \mathcal P f(y)|  & \leq \left (\int_y^\infty \langle x \rangle^{- 2 \sigma} \, d x \right)^{1/2}
\| \langle x \rangle^\sigma \mathcal H f \|_2
\leq C y^{- 1/12} \| f \|_X, \qquad y \geq 1.
\end{aligned}
\]
Oddness of $\mathcal P f$ extends \eqref{1.22} to $\R$.

\emph{Weighted interpolation.} For $R \geq 1$, let $I_R = \{ R < |y| < 2 R \}$ and $I_R^* = \{ R/2 < |y| < 4 R \}$. The rescaled local interpolation estimates are
\[
\| f^{(j)} \|_{L^2 (I_R)} \leq C \| f \|_{L^2 (I_R^*)}^{1 - \frac{j}{4}} \| f^{(4)} \|^{\frac{j}{4}}_{L^2(I_R^*)} + C R^{- j} \| f \|_{L^2 (I_R^*)},
\]
\[
\| f^{(j)} \|_{L^\infty (I_R)} \leq C \| f \|_{L^2 (I_R^*)}^{1 - \frac{2 j + 1}{8}} \| f^{(4)} \|^{\frac{2 j + 1}{8}}_{L^2(I_R^*)} + C R^{- j - \frac{1}{2}} \| f \|_{L^2 (I_R^*)}.
\]
The exponents satisfy
\[
\mathfrak a_j = \sigma (1 - \frac{j}{4}) + \frac{j}{2}, \, \, \mathfrak a_j - j - \sigma = - \frac{31 j}{48} \leq 0, 
\]
and
\[
\mathfrak b_j = \sigma (1 - \frac{2 j + 1}{8}) + \frac{2 j + 1}{4}, \, \, \mathfrak b_j - j - \frac{1}{2} - \sigma = - \frac{31 (j + \frac{1}{2})}{48} \leq 0. 
\]
For $k \geq 0$, set
\[
A_k = 2^{k \sigma} \| f \|_{L^2(I_{2^k}^*)}, \qquad
B_k = 2^{2 k} \| f^{(4)} \|_{L^2(I_{2^k}^*)}.
\]
The annuli have finite overlap, so
\[
\sum_k A_k^2 \leq C \| f \|_X^2, \qquad
\sum_k B_k^2 \leq C \| f^{(4)} \|_{L^2(\psi)}^2.
\]
For $0 \leq \theta \leq 1$, apply H\"older's inequality.
\[
\begin{aligned}
\sum_k A_k^{2(1 - \theta)} B_k^{2 \theta}
& \leq \left (\sum_k A_k^2 \right)^{1 - \theta}
\left (\sum_k B_k^2 \right)^\theta, \\
\sup_k A_k^{1 - \theta} B_k^\theta
& \leq \left (\sum_k A_k^2 \right)^{(1 - \theta)/2}
\left (\sum_k B_k^2 \right)^{\theta /2}.
\end{aligned}
\]
Use $\theta = j/4$ in the first bound and $\theta = (2 j + 1)/8$ in the second. The remainder terms are bounded by $C A_k$ because
\[
\mathfrak a_j - j - \sigma \leq 0, \qquad
\mathfrak b_j - j - \frac12 - \sigma \leq 0.
\]
On $|y| \leq 1$, unweighted interpolation completes \eqref{1.23} and \eqref{1.24}.

\emph{Hilbert transforms.} Choose an even $\chi \in C_c^\infty$ equal to one near zero. After integration by parts,
\[
\left | \int y \chi \left ( \frac{y}{R} \right) f^{(4)} d y \right | = \left | \int f(y) \partial_y^4 \biggl ( y \chi \left ( \frac{y}{R} \right) \biggr) d y \right |
\leq C R^{- \frac{5}{2}} \| f \|_2 \, \, \to \, \, 0,
\]
and
\[
\int | y f^{(4)} (y)| d y \leq \| \langle y \rangle^2 f^{(4)} \|_2 \|   y \langle y \rangle^{- 2} \|_2 < \infty. 
\]
We may therefore apply Lemma~\ref{1.20}. 
\[
\|   \langle y \rangle^{\sigma} \mathcal H f \|_2 + \|   \langle y \rangle^{2} \mathcal H f^{(4)} \|_2 \leq C \left ( \| f \|_X + \| f^{(4)} \|_{L^2 (\psi)} \right). 
\]
Interpolating $\mathcal H f$ with the same estimates completes the proof. 
\end{proof}

\bigskip

\section{Linear analysis}
In the spectral arguments, $X$ and the operators are complexified, and $X^*$ is the continuous complex-linear dual. We use
\[
\langle f, g \rangle = \int f \overline g,
\qquad \ell (h) = \int \zeta h.
\]
Nonlinear solutions and initial data are real-valued.

\begin{Proposition} \label{2.1}
The operator \eqref{1.15} generates strongly continuous semigroups on $L^2_{\mathrm{odd}}(\R)$ and $X$, with respective domains
\[
D(\mathcal L) = \{ h \in L^2_{\mathrm{odd}}(\R) : V h' \in L^2(\R) \}, \qquad
D_X(\mathcal L) = \{ h \in X : V h' \in X \}.
\]
Here derivatives are distributional. Moreover,
\begin{equation} \label{2.2}
\operatorname{Re} \langle h, \mathcal L h \rangle \leq - \frac{1}{4} \| h \|_2^2, \, \, \| e^{\tau \mathcal L} h \|_2 \leq e^{- \frac{\tau}{4}} \| h \|_2. 
\end{equation}
\end{Proposition}

\begin{proof}
\emph{1. Generation and domains.} For $\eta \in \{ 0, \sigma \}$, set
\[
w_\eta (y) = \langle y \rangle^{2 \eta}, \qquad
X_\eta = L^2_{\mathrm{odd}}(w_\eta \, d y), \qquad
\mathcal D_\eta = \{ h \in X_\eta : V h' \in X_\eta \}.
\]
Thus $X_0 = L^2_{\mathrm{odd}}(\R)$ and $X_\sigma = X$. On $y > 0$, the coordinate
\[
\xi (y) = \frac35 \log y + \frac65 \log (y^2 + 5), \qquad \xi'(y) = \frac1{V(y)},
\]
is a smooth increasing bijection from $(0, \infty)$ onto $\R$. Identifying an odd function with its restriction to the positive half line, define
\[
(\mathcal U_\eta h)(\xi) = \sqrt{2 w_\eta (y) V(y)} \, h(y).
\]
This map is unitary from $X_\eta$ onto $L^2(\R, d \xi)$. For smooth functions supported away from $y = 0$,
\[
\mathcal U_\eta (V h') = \partial_\xi (\mathcal U_\eta h)
- \left (\frac{V'} 2 + \eta \frac{V y}{1 + y^2} \right) \mathcal U_\eta h.
\]
The coefficient in parentheses is bounded. The identity holds distributionally on each open half-line. If $\mathcal U_\eta h \in H^1$, then
\[
h, V h' \in L^2(0, 1), \qquad
(V h)' = V' h + V h' \in L^2(0, 1).
\]
Thus $V h$ has a trace $b_+$ at $0^+$. If $b_+ \neq 0$, then, for some $\delta > 0$,
\[
|V(y) h(y)| \geq \frac{|b_+ |}{2}, \qquad
|h(y)| \geq \frac{3|b_+ |}{10 y} \quad (0 < y < \delta),
\]
contradicting $h \in L^2(0, 1)$. The trace at $0^-$ also vanishes by oddness. Consequently,
\[
[V h]_{0^- }^{0^+ } = 0, \qquad
(V h)' = V' h + V h' \quad \text{in } \mathcal D'(\R).
\]
Restriction to the half-lines proves the converse. Hence
\[
\mathcal U_\eta \mathcal D_\eta = H^1(\R_\xi),
\]
with no condition on $h(0)$.

Write $\mathcal L_0 = - V \partial_y - (1 + Q)$. Since $V' = 1/3 + Q/2$,
\[
\mathcal U_\eta \mathcal L_0 \mathcal U_\eta^{- 1}
= - \partial_\xi - \frac56 - \frac34 Q + \eta \frac{V y}{1 + y^2}.
\]
The group $g(\xi) \mapsto g(\xi - \tau)$ has generator $- \partial_\xi$ on $H^1(\R)$. By the bounded-perturbation theorem, $\mathcal L_0$ generates a semigroup with domain $\mathcal D_\eta$. For constants $0 < c_\eta \leq C_\eta < \infty$,
\[
c_\eta \| \mathcal U_\eta h \|_{H^1}
\leq \| h \|_{X_\eta} + \| \mathcal L_0 h \|_{X_\eta}
\leq C_\eta \| \mathcal U_\eta h \|_{H^1}.
\]
The odd extensions of $\mathcal U_\eta^{- 1} C_c^\infty (\R_\xi)$ form a smooth core supported away from zero and infinity.

The remaining operator $K = F \mathcal H + \frac12 F' \mathcal P$ is bounded on both $X_\eta$. Indeed, by Lemma~\ref{1.20} and Cauchy--Schwarz,
\[
\begin{gathered}
\| F \mathcal H h \|_{X_\eta} \leq C_\eta \| h \|_{X_\eta}, \qquad
| \mathcal P h(y)| \leq \sqrt{|y|} \| \mathcal H h \|_2
= \sqrt{|y|} \| h \|_2, \\
\| F' \mathcal P h \|_{X_\eta}
\leq \| \sqrt{|y|} \, w_\eta^{1/2} F' \|_2 \| h \|_2
\leq C_\eta \| h \|_{X_\eta}.
\end{gathered}
\]
The coefficient norm is finite because $F'$ is bounded near zero and $F'(y) = O(|y|^{- 4})$ at infinity. Thus $\mathcal L = \mathcal L_0 + K$ generates a semigroup with domain $\mathcal D_\eta$ and the same core.

\emph{2. The quadratic form gap.} We first work with real finite sine polynomials. The complex estimate follows by treating the real and imaginary parts separately. Set $\theta = 2 \arctan y$ and write
\[
v (\theta ) = \sum_{n \geq 1} a_n \sin n \theta, \, \, b_n = n a_n, \, \, f(\theta) = \sum b_n \sin n \theta, \, \, g(\theta) = \sum b_n \cos n \theta. 
\]
Here $f(\theta)$ and $g(\theta)$ denote function values. In these coordinates,
\[
h = - (1 + \cos \theta) f(\theta), \mathcal H h = ( 1 + \cos \theta ) g(\theta) = v_y, \, \, d y = \frac{d \theta}{1 + \cos \theta}.
\]
Thus, 
\[
\pi^{- 1} \| h \|_2^2 = \sum {b_n}^2 + \sum b_n b_{n + 1}. 
\]
Let $\mathsf{C}_j$ denote multiplication by $\cos j \theta$ in the sine basis. Then
\[
(\mathsf{C}_j)_{m n} = \frac{1}{2} (\mathbf 1_{|m - n | = j} - \mathbf 1_{m + n  = j}), \, \, \mathsf{C} = \mathsf{C}_1, \, \, \mathsf{C}_2 = 2 \mathsf{C}^2 - I, \, \, \mathsf{C}_3 = 4 \mathsf{C}^3 - 3 \mathsf{C}. 
\]
Let $\mathsf{N} = \operatorname{diag}(1, 2, 3, \ldots)$ and let $\mathsf{R}$ be the matrix whose only nonzero entries are $\mathsf{R}_{11} = \frac23$ and $\mathsf{R}_{12} = \mathsf{R}_{21} = \frac16$. Then we have
\[
\begin{split}
\pi^{- 1} \langle h, \mathcal L h \rangle
= {}  & b^T \Biggl [- \frac{11}{6} I - \frac{31}{12} \mathsf C - \mathsf C_2 - \frac14 \mathsf C_3 + \mathsf R \\
& + \frac12 \mathsf N^{- 1}
\left (\frac13 I + \mathsf C + \mathsf C_2 + \frac13 \mathsf C_3 \right)
+ \frac12
\left (\frac13 I + \mathsf C + \mathsf C_2 + \frac13 \mathsf C_3 \right) \mathsf N^{- 1}
\Biggr] b.
\end{split}
\]
The multiplication terms satisfy
\[
(1 + \cos \theta) (- \frac{5}{6} - \cos \theta - \cos^2 \theta ) = - \frac{11}{6} - \frac{31}{12} \cos \theta - \cos 2 \theta - \frac{1}{4} \cos 3 \theta. 
\]
For the Hilbert transform term,
\[
F = - \frac{4}{3} \sin \theta (1 + \cos \theta), \, \, \text{and} \, \, \pi^{- 1} \int F h \mathcal H h = b^T \mathsf{R} b
\]
after cancellation of the difference frequencies. For the potential term, we use
\[
\partial_\theta F = - \frac{4}{3} (\cos \theta + \cos 2 \theta),
\]
\[
\frac12 \int F' h \mathcal P h
= \int \left (\frac13 + \cos \theta + \cos 2 \theta + \frac13 \cos 3 \theta \right)
v f(\theta) \, d \theta,
\]
\[
\begin{split}
& \pi^{- 1} \int
\left (\frac13 + \cos \theta + \cos 2 \theta + \frac13 \cos 3 \theta \right)
v f(\theta) \, d \theta \\
& \quad = \frac12 b^T \Biggl [
\mathsf N^{- 1} \left (\frac13 I + \mathsf C + \mathsf C_2 + \frac13 \mathsf C_3 \right)
+
\left (\frac13 I + \mathsf C + \mathsf C_2 + \frac13 \mathsf C_3 \right) \mathsf N^{- 1}
\Biggr] b.
\end{split}
\]
Set $u_n = (- 1)^n b_n$ and $e_n = u_{n + 1} - u_n$, and extend the finite sequence by zero. Define 
\[
\mathcal Q (u) = \pi^{- 1} ( - \langle h, \mathcal L h \rangle - \frac{1}{4} \| h \|^2_2). 
\]
Then its diagonal terms are $q_{11} = \frac7{12}$ and, for $n \geq 2$, $q_{n n} = \frac{19}{12} - \frac1{3 n}$. The three upper bands are $- A_n$, $B_n$, and $- D_n$, respectively, where
\[
A_1 = \frac{5}{8}, \, \, A_n  = \frac{7}{6} - \frac{1}{4} \left (  \frac{1}{n} + \frac{1}{n + 1}   \right), \, \, \text{for} \, \, n \geq 2, 
\]
\[
B_n = \frac{1}{2} - \frac{1}{4} \left (   \frac{1}{n} + \frac{1}{n + 2} \right),
\qquad
D_n = \frac{1}{8} - \frac{1}{12} \left (    \frac{1}{n} + \frac{1}{n + 3} \right). 
\]
The row sums are
\[
\rho_1 = \frac{5}{48}, \, \, \rho_2 = \frac{19}{240}, \, \, \rho_3 = - \frac{41}{720},
\qquad
\rho_n =  \frac{6(n^2 + 1)}{ n (n^2 - 1) (n^2 - 4) (n^2 - 9 )}, \, \, \text{for} \, \, n \geq 4. 
\]
Define $s_0 = 0$, $s_n = - B_n + 2 D_n + 2 D_{n - 1}$, $t_n = A_n - B_n - B_{n - 1} + s_{n - 1} + s_n$, where $B_0 = 0 $ and $D_0 = 0 $. They satisfy
\[
s_1 = - \frac{1}{8}, \, \, s_n = \frac{1}{12} \left ( \frac{1}{n} + \frac{1}{n + 2} \right) - \frac{1}{6} \left (  \frac{1}{n - 1} + \frac{1}{n + 3}  \right) < 0, \, \, \text{for} \, \, n \geq 2, 
\]
\[
t_1 = \frac{1}{3}, \, \, t_2 = \frac{13}{60}, 
\]
for $n \geq 3$,
\[
t_n = \frac{1}{6} - \frac{1}{6(n - 2)(n - 1)} + \frac{1}{6(n + 2)(n + 3)} + \frac{1}{12} \left (   \frac{1}{n} + \frac{1}{n + 1}  \right) \geq \frac{1}{12},
\]
and
\[
D_n \geq D_1 = \frac{1}{48}. 
\]
The quadratic form decomposes as
\begin{align*}
\mathcal Q (u) = \frac{5}{48} u_1^2 + \sum_{n \geq 4} \rho_n u_n^2 + \frac{1}{45} \left ( u_2 - \frac{41}{16} e_2 \right)^2 + \frac{53}{3840} e_2^2 + \frac{1}{3} e_1^2 \\
+ \sum_{n \geq 3} t_n e_n^2 + \sum_{n \geq 1} D_n (e_n - e_{n + 1} + e_{n + 2})^2 + \sum_{n \geq 1 } (- s_n) ( e_n - e_{n + 1})^2.
\end{align*}
Every term in this identity is nonnegative, so $\mathcal Q(u) \geq 0$.

\emph{3. Extension to the generator domain.} Let $h \in L^2_{\mathrm{odd}}(\R)$ be orthogonal to $- (1 + \cos \theta) \sin n \theta$ for every $n \geq 1$. Then
\[
\int_{- \pi}^\pi \left |h \left (\tan \frac \theta 2 \right) \right |^2 d \theta
\leq 2 \| h \|_2^2,
\]
and
\[
\int_{- \pi}^\pi h \Bigg ( \tan \left ( \frac{\theta}{2} \right) \Bigg) \sin n \theta d \theta = 0,
\]
By completeness of the sine basis, $h = 0$. Thus the span is dense. On smooth real functions, integration by parts represents the quadratic form by
\[
\int \left (- \frac56 - \frac{3 Q}{4} \right) h^2
+ \int F h \mathcal H h + \frac12 \int F' h \mathcal P h.
\]
The first two terms satisfy
\[
\left | \int \left (- \frac56 - \frac{3 Q}{4} \right) h^2 \right |
+ \left | \int F h \mathcal H h \right |
\leq \left (\frac56 + \frac34 \| Q \|_\infty + \| F \|_\infty \right) \| h \|_2^2.
\]
For the last term,
\[
| \mathcal P h(y)| \leq \sqrt{|y|} \| h \|_2, \qquad
\left | \int F' h \mathcal P h \right |
\leq \| \sqrt{|y|} F' \|_2 \| h \|_2^2.
\]
The form is therefore continuous on $L^2$. Density and integration by parts extend the inequality to the smooth core, first for real functions and then for complex functions. For $h \in D(\mathcal L)$, choose core elements $h_j$ with
\[
\| h_j - h \|_2 + \| \mathcal L h_j - \mathcal L h \|_2 \longrightarrow 0.
\]
Then
\[
\operatorname{Re} \langle h, \mathcal L h \rangle
= \lim_{j \to \infty} \operatorname{Re} \langle h_j, \mathcal L h_j \rangle
\leq - \frac14 \lim_{j \to \infty} \| h_j \|_2^2
= - \frac14 \| h \|_2^2.
\]
For such $h$, the semigroup satisfies
\[
\frac12 \frac{d}{d \tau} \| e^{\tau \mathcal L} h \|_2^2
= \operatorname{Re} \langle e^{\tau \mathcal L} h, \mathcal L e^{\tau \mathcal L} h \rangle
\leq - \frac14 \| e^{\tau \mathcal L} h \|_2^2.
\]
By Gronwall's inequality and density of the domain, \eqref{2.2} holds for every $h \in L^2_{\mathrm{odd}}(\R)$.

\end{proof}

\begin{Proposition}\label{2.3}
There exist $\lambda_* \in (- \frac{23}{36}, - \frac12)$ and a nonzero smooth odd function $h \in D_X(\mathcal L)$ such that $\mathcal L h = \lambda_*h$. If $v = \mathcal P h$ and $\beta = 2 + 3 \lambda_*$, then, as $y \to \infty$,
\[
h(y) = O \left ( |y|^{- 1 - \beta} \right), \, \, v (y) = O \left ( |y|^{- \beta} \right). 
\]
\end{Proposition}

\begin{proof}
\emph{1. Analytic reduction.} Set $\theta = 2 \arctan y$ and write $v(\theta) = \operatorname{Im} W(e^{i \theta})$, where $W(z) = \sum_{n \geq 1} a_n z^n$ and $a_n \in \R$. Define
\begin{equation} \label{2.4}
\mathscr C = - \frac{(1 + z)^2}{2} W' (z) + \frac{a_1}{2}
\end{equation}
The endpoint behavior proved below justifies the identities
\[
h = \operatorname{Im} \mathscr C(e^{i \theta}), \qquad
\mathscr C(- 1) = \frac{a_1}{2}, \qquad
\mathcal H h = - \operatorname{Re} \mathscr C(e^{i \theta}) + \frac{a_1}{2} = v_y.
\]
The eigenvalue equation becomes
\begin{equation} \label{2.5}
\begin{split}
z^2 (z^2 - 1 ) ( z^2 + 3 z + 1) W'' + 2 z [ z^4 + 2 z^3 + 3 (1 + \lambda) z^2 + 3 z + 2] W' \\
- 2 (z + 1)^2 (z^2 - z + 1) W = 2 z [a_1 (z + 1)^2 + 2 a_2 z].
\end{split}
\end{equation}
The positive sine generating function of $\mathcal L h$ is
\[
\begin{split}
& - \frac{(z^2 - 1)(z^2 + 3 z + 1)}{6 z} \mathscr C'
- \left (\frac53 + \frac4{3 z} + \frac2{3 z^2} \right) \mathscr C
- \frac{(z + 1)^4(z^2 - z + 1)}{6 z^3} W \\
& \quad + \frac{a_1}{6 z^2} - \frac{a_2}{3 z}
- \left (\frac{a_1}{2} + a_2 + \frac{a_3}{3} \right)
- \frac{2 a_1 + a_2}{3} z - \frac{a_1}{6} z^2.
\end{split}
\]
Substitute \eqref{2.4} and set this expression equal to $\lambda \mathscr C$, using
\begin{equation} \label{2.6}
a_3 = - a_2 - \left (   1  + \frac{3 \lambda }{2} \right) a_1. 
\end{equation}
Multiplication by $\frac{12 z^3}{(z + 1)^2}$ implies \eqref{2.5}. Conversely, comparing the third Taylor coefficient of \eqref{2.5} recovers \eqref{2.6}. For $n \geq 3$, the recurrence is
\begin{align}
& ( - n^2 + 5 n - 2) a_n = 2 a_1 \mathbf 1_{n = 3} - (- 3 n^2 + 15 n - 14) a_{n - 1} \\
& \quad - 6 ( 1 + \lambda) (n - 2) a_{n - 2}
- (3 n^2 - 17 n + 22) a_{n - 3} - (n^2 - 7 n + 10) a_{n - 4}.
\end{align}
For $n \leq 0$, set $a_n = 0$. Since $- n^2 + 5 n - 2 \neq 0$, every choice of $(a_1, a_2) \in \R^2$ determines a solution such that, for some $\rho > 0 $, the series
\[
W(z) = \sum_{ n \geq 1 } a_n z^n
\]
converges for $|z| < \rho$. The singularities are $0$, $\pm 1$ and
\[
r = \frac{- 3 + \sqrt{5}}{2}, \, \, \frac{- 3 - \sqrt{5}}{2}
\]
At $1$, $r$, and $- 1$, the corresponding homogeneous exponents are
\[
\left ( 0, - \frac{6 + 3 \lambda}{5} \right), \, \, (0, \frac{- 7 - 2 \sqrt{5} - 6 \lambda }{5}), \, \, (0, 2 + 3 \lambda). 
\]
For $\lambda \in [- \frac{23}{36}, - \frac12]$, the exponents at $1$ and $r$ satisfy
\[
- \frac{6 + 3 \lambda}{5} \in (- 1, 0),
\qquad
\frac{- 7 - 2 \sqrt5 - 6 \lambda}{5} \in (- 2, - 1).
\]
\emph{2. Endpoint shooting.} We seek $q$ and $\lambda$ such that $W(z) = z + q z^2 + \cdots$ extends analytically through both $r$ and $1$. Normalize $a_1 = 1$ and $a_2 = q$, and write $W = z Y$. Then
\[
\begin{split}
& \mathcal T_\lambda Y
= {} z^2(z^2 - 1)(z^2 + 3 z + 1) Y''
+ 2 z[2 z^4 + 5 z^3 + 3(1 + \lambda) z^2 + 1] Y' \\
& \quad + 2[z^3 + 3(1 + \lambda) z^2 + 2 z + 1] Y
= {} 2(1 + z)^2 + 4 q z.
\end{split}
\]
For $z \in (r, 0) \cup (0, 1)$,
\[
z^2(z^2 - 1)(z^2 + 3 z + 1) < 0,
\qquad z^3 + 3(1 + \lambda) z^2 + 2 z + 1 > 0.
\]
On $[r, 0]$, the second inequality follows from
\[
z^3 + 3(1 + \lambda) z^2 + 2 z + 1
\geq z^2(1 + z) + 1 + 2 z \geq \sqrt5 - 2 > 0.
\]
At a positive interior maximum of $Z$,
\[
Z' = 0, \qquad Z'' \leq 0, \qquad
\mathcal T_\lambda Z
\geq 2[z^3 + 3(1 + \lambda) z^2 + 2 z + 1] Z > 0.
\]
For $e \in \{ r, 1 \}$, set
\[
b_e = \left |2 e[2 e^4 + 5 e^3 + 3(1 + \lambda) e^2 + 1] \right | > 0, \quad 
d_e = 2[e^3 + 3(1 + \lambda) e^2 + 2 e + 1] > 0.
\]
The coefficient of $Z'$ is negative at $r$ and positive at $1$. Thus $\mathcal T_\lambda Z(e) \leq 0$, for $Z$ analytic near $e$, is
\[
b_e \partial_n Z + d_e Z \leq 0. 
\]
A positive endpoint maximum has $\partial_n Z \geq 0$ and hence $b_e \partial_n Z + d_e Z > 0$. Therefore $\mathcal T_\lambda Z \leq 0$ and $Z(0) = 0$ imply $Z \leq 0$.
Set
\[
e_- = r, \qquad e_+ = 1, \qquad p_\pm = \frac{3 \pm \sqrt{17}} 2.
\]
Let $Y_{ \pm , b}$ be the solution of the preceding differential equation that is analytic near $e_\pm$ and satisfies $Y_{ \pm , b}(e_\pm ) = b$. Then
\[
Y'_{ \pm , b}(e_\pm )
= \frac{2(1 + e_\pm )^2 + 4 q e_\pm
- 2[e_\pm^3 + 3(1 + \lambda) e_\pm^2 + 2 e_\pm + 1] b}
{2 e_\pm [2 e_\pm^4 + 5 e_\pm^3 + 3(1 + \lambda) e_\pm^2 + 1]}.
\]
The coefficient of the divergent solution at $0$ is
\[
d_{ \pm } (b) = \lim_{z \to 0^\pm } |z |^{- p_{- }} Y_{ \pm , b } (z) = d_{ \pm , 0} + b d_{ \pm , 1}.
\]
If $d_{ \pm , 1} = 0$, the difference $Z = Y_{ \pm , 1} - Y_{ \pm , 0}$ is analytic at $e_\pm$ and satisfies
\[
\mathcal T_\lambda Z = 0, \qquad Z(e_\pm ) = 1, \qquad
\lim_{z \to 0^\pm } Z(z) = 0.
\]
The maximum principle applied to $\pm Z$ forces $Z = 0$, contradicting $Z(e_\pm ) = 1$. Thus
\[
d_{ \pm , 1} \neq 0, \qquad
Y_\pm = Y_{ \pm , - d_{ \pm , 0}/d_{ \pm , 1}}.
\]
This is the unique solution bounded at $0$ and analytic across $e_\pm$. Since $p_- < 0 < p_+ < 4$ and $p_+ > 3$,
\[
Y_\pm = Y_{\text{an}} + c_\pm (\lambda , q) |z|^{p_+ } H_{ \pm } (z)
\]
where $H_\pm (0) = 1$ and $Y_{\text{an}}$ is the analytic solution of the inhomogeneous equation near $0$, with the following expansion as $z \to 0$
\[
Y_{\text{an}} (z) = 1 + q z + a_3 z^2 + a_4 z^3 + O (z^4).
\]
Its Taylor polynomial is
\begin{equation} \label{2.7}
P_3 = 1 + q z + a_3 z^2 + a_4 z^3, \, \, a_3 = - q - 1 - \frac{3}{2} \lambda, \, \, a_4 = - (7 + 6 \lambda) q - 2 - \frac{3}{2} \lambda. 
\end{equation}
For $- 23/36 \leq \lambda \leq - 1/2$, the exponent differences satisfy
\[
\left \{ \sqrt{17}, - \frac{6 + 3 \lambda}{5},
\frac{- 7 - 2 \sqrt5 - 6 \lambda}{5}, 2 + 3 \lambda \right \} \cap \mathbb Z = \varnothing.
\]
The continued coefficients have the form
\[
c_\pm (\lambda , q) = A_\pm (\lambda) + q B_\pm (\lambda )
\]
where $A_\pm , B_\pm \in C([- \frac{23}{36}, - \frac12])$.

\begin{Lemma}\label{2.8}
Fix $\lambda \in [- \frac{23}{36}, - \frac{1}{2} ]$ and $q \in \R$. Let $I_{- } = [r, 0 ]$, $I_{+ } = [0 , 1]$, $ r = \frac{- 3 + \sqrt{5}}{2}$. Let $P$ be real analytic on a neighborhood of $I_\pm$. With $P_3$ defined by \eqref{2.7}, assume that 
\[
P(0) = 1, \, \, P'(0) = q, \, \, P''(0) = 2 a_3 , \, \, P'''(0) = 6 a_4
\]
Set $R_P (z) = \mathcal T_\lambda P(z) - 2(1 + z)^2 - 4 q z $. If $\min_{z \in I_{ \pm }} R_P (z) \geq 0$ and $\max_{z \in I_\pm } R_P (z) > 0$, then
\[
c_{ \pm } (\lambda, q ) < 0,
\]
if $\max_{z \in I_\pm } R_P(z) \leq 0$ and $\min_{z \in I_\pm } R_P(z) < 0$, then
\[
c_{ \pm } (\lambda, q ) > 0. 
\]
\end{Lemma}

\begin{proof}
Assume first that $R_P \geq 0$ and is not identically zero. With $Z = Y_\pm - P$, we have $Z(0) = 0$ and $\mathcal T_\lambda Z = - R_P \leq 0$, so $Z \leq 0$ by the maximum principle. The matching Taylor coefficients and $p_+ < 4$ imply
\[
Z = c_\pm |z|^{p_+ }(1 + o(1)) + O(|z|^4),
\]
hence $c_\pm \leq 0$.

To exclude equality, suppose $c_\pm = 0$. Then $Z$ is analytic at zero. In the coordinate $x = |z|$, write the first nonzero residual term as
\[
R_P = r_n x^n + O(x^{n + 1}), \qquad r_n > 0, \quad n \geq 4.
\]
The differential equation determines the leading term.
\[
Z = z_n x^n + O(x^{n + 1}), \qquad
(- n^2 + 3 n + 2) z_n = - r_n.
\]
Since $- n^2 + 3 n + 2 < 0$ for $n \geq 4$, we obtain $z_n > 0$, contradicting $Z \leq 0$. Thus $c_\pm < 0$. Reversing all signs proves the second assertion.
\end{proof}

\emph{3. Residual comparisons and a common zero.}

For every $n \geq 0$,
\begin{align*}
\mathcal T_\lambda z^n = ( - n^2 + 3 n + 2) z^n + (- 3 n^2 + 3 n + 4 ) z^{n + 1} + 6 (1 + \lambda ) (n + 1) z^{n + 2} \\
+ (3 n^2 + 7 n + 2) z^{n + 3} + (n^2 + 3 n) z^{n + 4}. 
\end{align*}
Consequently,
\begin{equation} \label{2.9}
\begin{split}
\frac{R_{P_3}}{z^4} = (92 + 66 \lambda)  q + 10 - 24 \lambda - 27 \lambda^2 + [24 ( 1 + \lambda) a_4 + 28 a_3 + 4 q] z \\
+ (50 a_4 + 10 a_3 ) z^2 + 18 a_4 z^3. 
\end{split}
\end{equation}
Set $t = \lambda + \frac{2}{3} \in [\frac{1}{36}, \frac{1}{6}]$ and $s = - z \in [0, - r]$, so that $s < \frac{2}{5}$. For $P = P_3$ and $q = - \frac{1}{3}$, \eqref{2.9} becomes
\[
\frac{R_P}{z^4} = - 2 - 10 t - 27 t^2 + (8 - 38 t + 12 t^2 ) z + \left ( \frac{10}{3} + 10 t \right) z^2 + 9 t z^3
\]
Since $(8 - 38 t + 12 t^2 ) \geq 2 $ and $  \left ( \frac{10}{3} + 10 t \right) \leq 5$ for $ t \in [\frac{1}{36}, \frac{1}{6}]$,
\[
\frac{R_P}{z^4} \leq - \frac{6}{5}. 
\]
Therefore, $c_{- } \left ( \lambda, - \frac{1}{3} \right) > 0$. For $q = - \frac{1}{4}$, take $P = P_3 + \frac{z^4 }{12}$. Its residual divided by $z^4$ is 
\[
\frac{11}{6} - \frac{9}{2} t - 27 t^2 + \left ( \frac{4}{3} - 48 t \right) z - \left ( \frac{55}{6} + \frac{25}{2} t \right) z^2  + 2 z^3 + \frac{7}{3} z^4. 
\]
For $s = - z \in [0, - r]$, the quotient $R_P(- s)/s^4$ extends continuously to $s = 0$ and is concave in $t$. Since $0 \leq s < \frac25$, the endpoint bounds are
\[
\left. \frac{R_P(- s)}{s^4} \right |_{t = 1/36}
\geq \frac53 - \frac{48}{5} s^2 - 2 s^3 > 0,
\]
\[
\left. \frac{R_P(- s)}{s^4} \right |_{t = 1/6}
\geq \frac13 + s \left (\frac{20}{3} - \frac{45}{4} s - 2 s^2 \right) + \frac73 s^4 > 0.
\]
Thus $c_{- } \left ( \lambda, - \frac{1}{4} \right) < 0 $. By affine dependence on the shooting parameter, there is a unique continuous function $q_{- } (\lambda )$ such that 
\begin{equation} \label{2.10}
c_{- } \left (  \lambda, q_{- } (\lambda ) \right) = 0, \, \, - \frac{1}{3} < q_{- } (\lambda ) < - \frac{1}{4}. 
\end{equation}
For this comparison, fix $\lambda = - \frac{23}{36}$ and write $P_3(z; q)$ for the polynomial in \eqref{2.7} at this value of $\lambda$. For $q = - \frac{1}{3}$ and $- \frac{1}{4}$, define
\[
P_{- \frac{1}{3}} (z) = P_3 \left (z ; - \frac{1}{3} \right) - \frac{ z^4 (5 + 2 z + 2 z^2)}{4 ( 1 + z )^4},
\]
\[
P_{- \frac{1}{4}} (z) = P_3 \left (z ; - \frac{1}{4} \right) + \frac{ z^4 (4 + z )}{8 ( 1 + z )^2}.
\]
Both functions are analytic near $[0, 1 ]$ and agree with $P_3$ through Taylor order $3$. Define
\begin{equation} \label{2.11}
\begin{split}
H_{- 1/3}(z)  & = \frac{(1 + z)^6}{z^4}
\left [\mathcal T_{- 23/36} P_{- 1/3}(z) - 2(1 + z)^2 + \frac43 z \right], \\
H_{- 1/4}(z)  & = \frac{(1 + z)^4}{z^4}
\left [\mathcal T_{- 23/36} P_{- 1/4}(z) - 2(1 + z)^2 + z \right].
\end{split}
\end{equation}
These quotients extend to polynomials of degrees $n_{- 1/3} = 9$ and $n_{- 1/4} = 7$, respectively. Their Bernstein expansions on $[0, 1]$ are 
\begin{equation} \label{2.12}
\begin{split}
H_q (z) = \sum_{j = 0 }^{n_q} b_{q, j } \binom{n_q}{j} z^j (1 - z )^{n_q - j} \\
b_{q, j} = \sum_{i = 0 }^j \frac{ \binom{j}{i} }{ \binom{n_q}{i} } \cdot \frac{H_q ^{(i)} (0) }{i ! }, \, \, \text{where} \, \, 0 \leq j \leq n_q. 
\end{split}
\end{equation}
Define
\[
(c_0, \ldots, c_9) = \left (\tfrac15, \tfrac32, 2, \tfrac52, 2, 2, 3, 5, 11, 17 \right), \qquad 
(d_0, \ldots, d_7) = \left (\tfrac45, \tfrac25, \tfrac15, \tfrac23, \tfrac52, 6, 10, 4 \right).
\]
Substituting into \eqref{2.11} and \eqref{2.12}, we find
\[
b_{- \frac13, j} \geq c_j \quad (0 \leq j \leq 9), \qquad
b_{- \frac14, j} \geq d_j \quad (0 \leq j \leq 7).
\]
All lower bounds are at least $\frac{1}{5 }$. For $0 \leq z \leq 1 $,
\[
\binom{n_q}{j} z^j (1 - z)^{n_q - j} \geq 0, \, \, \sum_{j = 0 }^{n_q } \binom{n_q}{j} z^j (1 - z)^{n_q - j} = 1,
\]
Thus $H_q \geq \frac15$ on $[0, 1]$, and, for $0 < z \leq 1$,
\[
\begin{split}
R_{P_{- 1/3}}(z)  & = \frac{z^4 H_{- 1/3}(z)}{(1 + z)^6}
\geq \frac{z^4}{5(1 + z)^6} > 0, \\
R_{P_{- 1/4}}(z)  & = \frac{z^4 H_{- 1/4}(z)}{(1 + z)^4}
\geq \frac{z^4}{5(1 + z)^4} > 0.
\end{split}
\]
Lemma~\ref{2.8} implies $c_+ (- \frac{23}{36}, - \frac13) < 0$ and $c_+ (- \frac{23}{36}, - \frac14) < 0$. Since the dependence on $q$ is affine,
\begin{equation} \label{2.13}
c_+ \Bigg (  - \frac{23}{36}, q_- \left (  - \frac{23}{36} \right)  \Bigg) < 0. 
\end{equation}
For $\lambda = - \frac{1}{2}$ and $q = - \frac{5}{16} $, the residual of $P_3 $ satisfies
\[
\frac{R_{P_3}}{z^4} = - \frac{51}{16} + \frac{z}{2} + \frac{ 5 z^2 }{8}
\]
The displayed polynomial is negative on both intervals, so
\[
c_\pm \left ( - \frac{1}{2}, - \frac{5}{16} \right) > 0. 
\]
Equation~\eqref{2.10} then locates the shooting root.
\[
- \frac{5}{16} < q_- \left ( - \frac{1}{2} \right) < - \frac{1}{4}. 
\]
For $q = - \frac14$, take $P = P_3 + \frac{z^4}{4(1 + z)}$. On the positive interval $0 < z \leq 1$, its residual satisfies 
\[
R_P = - \frac{ z^5 ( 40 + 115 z + 123 z^2 + 58 z^3 + 10 z^4 )}{4 (1 + z)^3 } < 0. 
\]
Therefore $ c_+ \left ( - \frac{1}{2} , - \frac{1}{4} \right) > 0 $. The affine dependence preserves this sign at the shooting root.
\begin{equation} \label{2.14}
c_+ \Bigg ( - \frac{1}{2}, q_- \left ( - \frac{1}{2} \right) \Bigg) > 0 .
\end{equation}
By continuity, \eqref{2.13} and \eqref{2.14} imply the existence of $\lambda_* \in (- \frac{23}{36}, - \frac12)$ such that $c_\pm (\lambda_*, q_- (\lambda_*)) = 0$. Thus $W = z Y_{\text{an}}$ extends across $r$ and $1$ and is analytic on $|z| < 1$.

\emph{4. Return to the physical eigenfunction.} At the remaining singularity,
\[
W (z) = A(z) + (1 + z)^\beta B(z) 
\]
where
\[
\beta = 2 + 3 \lambda_* \in \left ( \frac{1}{12}, \frac{1}{2} \right)
\]
and both $A$ and $B$ are analytic near $- 1$. In particular, $W$ is bounded. Returning to $\theta = 2 \arctan y $, we obtain from \eqref{2.4}
\[
h(y) = \operatorname{Im} \mathscr C (e^{i \theta }), \, \, v(y) = \operatorname{Im} W (e^{i \theta })
\]
with $\mathscr C (- 1) = \frac{1}{2} $, $ \mathcal H h = v_y $, and $\mathcal L h = \lambda_* h $. These functions satisfy
\[
| h(y) | + | V(y) h' (y) | = O ( |y |^{- 1 - \beta }) + O ( |y|^{- 3}), \, \, v(y) = O (|y|^{- \beta} ). 
\]
Both fields are odd. Analyticity at $z = 1$ implies smoothness at $y = 0$, and
\[
\int_1^\infty y^{2 \sigma - 2 - 2 \beta} \, d y < \infty,
\qquad \beta > \sigma - \frac12 
\]
yields $h, V h' \in X$.
If $h = 0$, then $\operatorname{Im} \mathscr C = 0$ on the unit circle. Bounded analyticity and $\mathscr C(0) = 0$ force $\mathscr C = 0$. Equation~\eqref{2.4} then implies
\[
W'(z) = \frac1{(1 + z)^2}, \qquad W(0) = 0, \qquad
W(z) = \frac z{1 + z},
\]
contradicting boundedness at $- 1$. Thus $h \neq 0$.
\end{proof}

\begin{Proposition}\label{2.15}
There exist a smooth even function $\varrho$ and constants $c, C > 0$ such that $c w(y) \leq \varrho (y) \leq C w(y)$ for every $y \in \R$ and, for every $h \in D_X(\mathcal L)$,
\begin{equation} \label{2.16}
\operatorname{Re} \langle \mathcal L h, h \varrho \rangle \leq - \frac{6}{25} \| h \|_{L^2 (\varrho)}^2. 
\end{equation}
Consequently,
\[
\| e^{\tau \mathcal L} h \|_{L^2(\varrho)}
\leq e^{- 6 \tau /25} \| h \|_{L^2(\varrho)}, \qquad h \in X.
\]
The spectrum on $X$ in $\operatorname{Re} z > - 23/36$ consists of isolated eigenvalues of finite multiplicity. There exists $0 \neq \ell \in X^*$ such that
\[
\ell (\mathcal L h) = \lambda_* \ell (h), \qquad h \in D_X(\mathcal L).
\]
\end{Proposition}

\begin{proof}
\emph{1. A coercive equivalent weight.} The transport coefficient for $w = \langle y \rangle^{2 \sigma}$ is
\[
d_\sigma = - \frac56 - \frac34 Q + \sigma \frac{V y}{1 + y^2},
\qquad d_\sigma (\infty) = - \frac56 + \frac \sigma 3 = - \frac{23}{36}.
\]
Choose an even $\chi \in C^\infty (\R)$ such that
\[
0 \leq \chi \leq 1, \qquad
\chi (y) = 0 \quad (|y| \leq 1), \qquad
\chi (y) = 1 \quad (|y| \geq 2).
\]
For $R \geq 1$, set
\[
\chi_R(y) = \chi \left (\frac y R \right), \qquad
E_R^2 = \int \chi_R w|h|^2.
\]
Then
\[
\operatorname{supp} \chi_R' \subset \{ R \leq |y| \leq 2 R \}, \qquad
\| \chi_R' \|_\infty \leq \frac{\| \chi' \|_\infty}{R}.
\]
For $h$ in the smooth core of Proposition~\ref{2.1}, integrate by parts.
\[
\begin{split}
\operatorname{Re} \int \chi_R w(\mathcal L h) \bar h
= {}  & \int \chi_R w d_\sigma |h|^2
+ \frac12 \int \chi_R' w V|h|^2 \\
& + \operatorname{Re} \int \chi_R w F \bar h \mathcal H h
+ \frac12 \operatorname{Re} \int \chi_R w F' \bar h \mathcal P h.
\end{split}
\]

On $R \leq |y| \leq 2 R$, the bounds $|V| \leq C R$, $w \leq C R^{2 \sigma}$, and $| \chi_R'| \leq C/R$ yield
\[
\left | \frac12 \int \chi_R' w V|h|^2 \right |
\leq C R^{2 \sigma} \| h \|_2^2.
\]
The explicit profile satisfies
\[
\sup_{|y| \geq R}|F(y)| \leq C R^{- 3}, \qquad
|F'(y)| \leq C|y|^{- 4} \quad (|y| \geq 1).
\]

Lemma~\ref{1.20}, with $\sigma = 7/12$, yields
\[
\| \chi_R^{\frac12} w^{\frac12} F \mathcal H h \|_2
\leq \sup_{|y| \geq R}|F(y)| \, \| w^{\frac12} \mathcal H h \|_2
\leq C R^{- 3} \| h \|_X,
\]
By \eqref{1.22}, $| \mathcal P h(y)| \leq C|y|^{- 1/12} \| h \|_X$ for $|y| \geq 1$. Hence
\[
\begin{split}
& \| \chi_R^{\frac12} w^{\frac12} F' \mathcal P h \|_2^2
= \int \chi_R w|F'|^2| \mathcal P h|^2 \\
& \leq C \| h \|_X^2 \int_{|y| \geq R}|y|^{\frac76 - 8 - \frac16} \, d y
= C \| h \|_X^2 \int_{|y| \geq R}|y|^{- 7} \, d y
\leq C R^{- 6} \| h \|_X^2.
\end{split}
\]
Cauchy--Schwarz bounds these pairings.
\[
\left | \int \chi_R w F \bar h \mathcal H h
+ \frac12 \int \chi_R w F' \bar h \mathcal P h \right |
\leq C R^{- 3} E_R \| h \|_X.
\]
Hence
\[
\operatorname{Re} \int \chi_R w(\mathcal L h) \bar h
\leq \left (\sup_{|y| \geq R} d_\sigma (y) \right) E_R^2
+ C R^{2 \sigma} \| h \|_2^2 + C R^{- 3} E_R \| h \|_X,
\]
where $C$ depends only on the fixed profiles, weight, and cutoff, and is independent of $R \geq 1$, $h$, and $\tau$.
Since $\chi_R = 1$ for $|y| \geq 2 R$,
\[
\| h \|_X^2 \leq (1 + 4 R^2)^\sigma \| h \|_2^2 + E_R^2,
\]
and, by Young's inequality,
\[
C R^{- 3} E_R \| h \|_X
\leq \left (C R^{- 3} + \frac1{16} \right) E_R^2 + C_R \| h \|_2^2.
\]
Choose $R$ so that
\[
\sup_{|y| \geq R} d_\sigma (y) \leq - \frac12,
\qquad C R^{- 3} \leq \frac1{16},
\qquad \nu = \frac38 > \frac14.
\]
Then
\[
\operatorname{Re} \int \chi_R w(\mathcal L h) \bar h
\leq - \nu E_R^2 + C_R \| h \|_2^2.
\]
Fix this $R$ and choose $M \geq \max \{ 1, 100 C_R \}$. Define
\[
\varrho = M + \chi_R w.
\]
The weight is smooth and even, with
\[
c w \leq \varrho \leq C w, \qquad
\left \| \frac{y \varrho'}{\varrho} \right \|_\infty < \infty.
\]
Adding $M$ times the inequality in Proposition~\ref{2.1} completes the coercive bound.
\[
\operatorname{Re} \int (\mathcal L h) \bar h \varrho
\leq - \left (\frac M 4 - C_R \right) \| h \|_2^2 - \nu E_R^2 
\leq - \frac6{25} \left (M \| h \|_2^2 + E_R^2 \right)
= - \frac6{25} \| h \|_{L^2(\varrho)}^2.
\]
The core is graph dense, so the form inequality extends to $D_X(\mathcal L)$. The form inequality holds for strong solutions. Density extends the semigroup estimate to all data.
For $h \in X$, define
\[
q_\varrho (h) = \int \left (- 1 - Q + \frac{(V \varrho)'}{2 \varrho} \right) |h|^2 \varrho
+ \operatorname{Re} \int (K h) \bar h \varrho,
\]
where $K = F \mathcal H + \frac12 F' \mathcal P$. The coefficient is bounded because $Q$, $V'$, $V/y$, and $y \varrho'/ \varrho$ are bounded. Since $K : X \to X$ is bounded, $q_\varrho$ is continuous on $X$. It agrees with the preceding pairing on the core. Density in $X$ therefore gives
\[
q_\varrho (h) \leq - \frac6{25} \| h \|_{L^2(\varrho)}^2, \qquad h \in X.
\]
For $h \in X \cap H^1_{\mathrm{loc}}(\R)$, let $\eta_S(y) = 1 - \chi (y/S)$, $S \geq 1$. The extra transport term in the cutoff integration by parts satisfies
\[
\left | \int \eta_S \eta_S' V |h|^2 \varrho \right |
\leq C \int_{|y| \geq S} |h|^2 \varrho \longrightarrow 0,
\]
with $C$ independent of $S$. The bounded coefficient and $K h \in X$ then imply
\[
\lim_{S \to \infty} \operatorname{Re} \int \eta_S^2 (\mathcal L h) \bar h \varrho = q_\varrho (h).
\]

\emph{2. The transport resolvent.} Write $\mathcal L = \mathcal L_0 + K$, with $\mathcal L_0 = - V \partial_y - (1 + Q)$ and $K = F \mathcal H + \frac12 F' \mathcal P$. In the coordinate $\xi$, the unitary map $\mathcal U_\sigma$ from Proposition~\ref{2.1} satisfies
\[
\mathcal U_\sigma \mathcal L_0 \mathcal U_\sigma^{- 1}
= - \partial_\xi + d_\sigma (\xi),
\]

where $d_\sigma (- \infty) = - \frac{17}{6}$ and $d_\sigma (\infty) = - \frac{23}{36}$. The positive excess of this coefficient is integrable.
\[
\int_\R \left (d_\sigma + \frac{23}{36} \right)_+ \, d \xi < \infty.
\]
Indeed,
\[
d_\sigma + \frac{23}{36}
= \frac{93 y^2 - 79}{36(1 + y^2)^2}, \quad 
\int_\R \left (d_\sigma + \frac{23}{36} \right)_+ d \xi
= \int_{\sqrt{79/93}}^\infty
\frac{93 y^2 - 79}{36(1 + y^2)^2 V(y)} \, d y < \infty.
\]
The integrand is $O(y^{- 3})$ at infinity. The transport semigroup is
\[
( e^{ \tau \mathcal U_\sigma \mathcal L_0 \mathcal U_\sigma^{- 1}} g ) (\xi ) = \exp \left ( \int_{\xi - \tau }^\xi d_\sigma (r) d r \right) g (\xi - \tau),
\]
and
\[
\| e^{\tau \mathcal L_0} \|_{X \to X}
\leq \exp \left (\int_\R \left (d_\sigma + \frac{23}{36} \right)_+ \, d \xi \right)
e^{- \frac{23 \tau}{36}}.
\]
Also,
\[
(z - \mathcal L_0)^{- 1} = \int_0^\infty e^{- z \tau} e^{\tau \mathcal L_0} d \tau.
\]
\emph{3. Relative compactness.} Fix $\operatorname{Re} z > - \frac{23}{36}$, let $\mathcal R_0(z) = (z - \mathcal L_0)^{- 1}$, and take $0 < \delta < R < \infty$. For
\[
A_{\delta, R} = \{ y \in \R : \delta < |y| < R \},
\]
the resolvent equation reads
\[
V \partial_y(\mathcal R_0(z) f) = f - (z + 1 + Q) \mathcal R_0(z) f, \qquad f \in X,
\]
\[
\| \mathcal R_0(z) f \|_X
\leq \frac{\displaystyle \exp \left (\int_\R
\left (d_\sigma + \frac{23}{36} \right)_+ \, d \xi \right)}
{\operatorname{Re} z + \frac{23}{36}} \| f \|_X.
\]
Since $V$ does not vanish on $\overline{A_{\delta /2, 2 R}}$, the resolvent equation and bound show that
\[
\| \mathcal R_0(z) f \|_{H^1(A_{\delta /2, 2 R})}
\leq C_{z, \delta, R} \| f \|_X.
\]
By Rellich compactness, the restriction
\[
\mathbf1_{A_{\delta /2, 2 R}} \mathcal R_0(z) : X \longrightarrow L^2(A_{\delta /2, 2 R})
\]
is compact. For the zero extension $\widetilde g$ of $g \in L^2(A_{\delta /2, 2 R})$,
\[
\| \mathbf1_{A_{\delta, R}} F \mathcal H \widetilde g \|_{L^2(w \, d y)}
\leq \| \mathbf1_{A_{\delta, R}} w^{1/2} F \|_\infty
\| g \|_{L^2(A_{\delta /2, 2 R})}.
\]
The composition preserves oddness. Thus
\[
\mathbf1_{A_{\delta, R}} F \mathcal H
\mathbf1_{A_{\delta /2, 2 R}} \mathcal R_0(z) : X \longrightarrow X
\]
is compact.

For the complementary integral,
\[
|y - x| \geq \frac \delta 2 \quad \text{if } y \in A_{\delta, R}, \ |x| \leq \frac \delta 2,
\qquad
|y - x| \geq \frac{|x|}{2} \quad \text{if } y \in A_{\delta, R}, \ |x| \geq 2 R.
\]
After multiplication of the input and output by $w^{\frac12}$, the integral kernel satisfies
\[
\begin{split}
& \iint_{\R^2} \left |
\frac{\mathbf1_{A_{\delta, R}}(y) \mathbf1_{\R \setminus A_{\delta /2, 2 R}}(x)
w(y)^{\frac12} F(y)}{\pi (y - x) w(x)^{\frac12}}
\right |^2 \, d x \, d y \\
& \quad \leq \frac4{\pi^2} \int_{A_{\delta, R}} w(y)|F(y)|^2 \, d y
\left (\frac1 \delta + \int_{|x| \geq 2 R} \frac{d x}{x^2 w(x)} \right) < \infty.
\end{split}
\]
Hence $\mathbf1_{A_{\delta, R}} F \mathcal H \mathbf1_{\R \setminus A_{\delta /2, 2 R}}$ is Hilbert--Schmidt on $L^2(w \, d y)$. Its composition with the bounded operator $\mathcal R_0(z)$ is compact. Both terms preserve oddness, and
\[
\mathbf1_{A_{\delta, R}} F \mathcal H \mathcal R_0(z)
= {} \mathbf1_{A_{\delta, R}} F \mathcal H \mathbf1_{A_{\delta /2, 2 R}} \mathcal R_0(z) 
+ \mathbf1_{A_{\delta, R}} F \mathcal H \mathbf1_{\R \setminus A_{\delta /2, 2 R}} \mathcal R_0(z)
\]
is therefore compact from $X$ to $X$. The constants above depend only on $z, \delta, R$ and the fixed profiles and weight.
As $\delta \to 0$ and $R \to \infty$, 
\[
\|  \mathbf 1_{\R \setminus A_{\delta, R}}  F \mathcal H \mathcal R_0 (z) \|_{X \to X } \leq C_z \sup_{\R \setminus A_{\delta, R }} |F | \to 0. 
\]
For the potential term, $(\mathcal P h )' = \mathcal H h$ and $ \|    \mathcal P h \|_{H^1 (- R, R )} \leq C_R \| h \|_X $, while 
\[
\|    \mathbf 1_{|y| > R } F' \mathcal P h  \|_X \leq C  \|  \mathbf 1_{|y| > R } w^{\frac{1}{2}} F' |y|^{- \frac{1}{12}}            \|_2 \| h \|_X
= o_R (1) \| h \|_X. 
\]
Local Rellich compactness and this tail bound show that $F' \mathcal P : X \to X$ is compact. Consequently $K \mathcal R_0(z)$ is compact.

\emph{4. Spectral isolation and the left eigenfunctional.} The factorization
\[
z - \mathcal L = (I - K \mathcal R_0(z))(z - \mathcal L_0)
\]
and the analytic Fredholm theorem, starting from a sufficiently large positive $z$, show that $\sigma (\mathcal L ) \cap \{ \operatorname{Re} z > - \frac{23}{36}  \}$ consists of isolated eigenvalues of finite multiplicity. Together with Proposition~\ref{2.3}, this establishes
\[
\lambda_* > - \frac{23}{36}, \, \, 0 \neq \ell \in X^*, \, \, \ell \mathcal L = \lambda_* \ell. 
\]
Since $\mathcal L$ commutes with conjugation and $\lambda_* \in \R$, for real $h \in D_X(\mathcal L)$,
\[
\operatorname{Re} \ell (\mathcal L h) = \lambda_* \operatorname{Re} \ell (h), \qquad
\operatorname{Im} \ell (\mathcal L h) = \lambda_* \operatorname{Im} \ell (h).
\]
At least one of these real functionals is nonzero. Complexifying it and retaining the notation $\ell$, we have
\[
\ell (\overline h) = \overline{\ell (h)}, \qquad h \in X.
\]

\end{proof}

\begin{Proposition} \label{2.17}
There exists a real odd $\zeta $ such that $\ell (h) = \int \zeta h$, and $\zeta \in H^1 (\R) $, $y \zeta ' \in L^2 (\R)$, $ (\mathcal H \zeta )' \in L^1 (\R) $. For every $h \in X$, we have
\begin{equation} \label{2.18}
\ell (h) = - \int (\mathcal H \zeta)' \mathcal P h, \, \, | \ell (h) | \leq \| (\mathcal H \zeta)' \|_1 \| \mathcal P h \|_\infty. 
\end{equation}
\end{Proposition}

\begin{proof}
\emph{1. The adjoint source.} By Riesz representation, there is a real odd $\zeta$ with
\[
\ell (h) = \int \zeta h, \qquad
\| w^{- 1/2} \zeta \|_2 = \| \ell \|_{X^*}.
\]
For odd functions,
\[
\widehat{\mathcal P h}(\eta) = - \frac{\widehat h(\eta)}{| \eta |}.
\]
Testing against compactly supported odd functions identifies the adjoint equation.
\begin{equation} \label{2.19}
V \zeta' - \left (    \frac{2}{3} + \frac{Q}{2} + \lambda_* \right) \zeta = \mathcal H (F \zeta) - \frac{1}{2} \mathcal P ( F' \zeta ).
\end{equation}
Set
\[
R_\zeta = \mathcal H(F \zeta) - \frac12 \mathcal P(F' \zeta).
\]

\[
\| F \zeta \|_2 + \| F' \zeta \|_2 \leq \left (  \| w^{\frac{1}{2} } F  \|_\infty   +   \| w^{\frac{1}{2} } F'  \|_\infty \right) \| w^{- \frac{1}{2} } \zeta  \|_2 < \infty 
\]
and 
\[
\| y F' \zeta \|_1 \leq \|  y F' w^{\frac{1}{2} }      \|_2 \|   w^{- \frac{1}{2}} \zeta   \|_2 < \infty
\]
The bounds imply $F' \zeta \in L^1$. Since this product is odd, $\int F' \zeta = 0$, so
\[
| \widehat{F' \zeta}(\eta)|
= \left | \int (e^{- i y \eta} - 1) F'(y) \zeta (y) \, d y \right |
\leq | \eta | \| y F' \zeta \|_1.
\]

\[
\| \mathcal P(F' \zeta) \|_{H^1}
\leq C \bigl (\| y F' \zeta \|_1 + \| F' \zeta \|_2 \bigr), \qquad
\| R_\zeta \|_2
\leq \| F \zeta \|_2 + \frac12 \| \mathcal P(F' \zeta) \|_2 < \infty.
\]
Since $V(y) \neq 0$ for $y \neq 0$, equation~\eqref{2.19} shows that $\zeta \in H^1_{\mathrm{loc}}(\R \setminus \{ 0 \})$.

\emph{2. Regularity at the origin.} Extend $F/y$ continuously by $F'(0)$. By oddness, 
\[
\int \frac F y \zeta = 0, \qquad \text{which yields} \qquad
\mathcal H (F \zeta) = y \mathcal H \left (\frac F y \zeta \right)
\]
and
\[
\left \| \frac{\mathcal P(F' \zeta)}{y} \right \|_{L^2(0, 1)} \leq 2 \| \mathcal H(F' \zeta) \|_{L^2(0, 1)},
\]
Consequently,
\[
\frac{R_\zeta}{V} \in L^2(0, 1), \qquad
\frac{\frac23 + \frac Q 2 + \lambda_*}{V} = \frac \alpha y + O(y),
\qquad \alpha = \frac{6 + 3 \lambda_*}{5} > \frac{49}{60} > \frac12.
\]
Thus, with
\[
g_0 = \frac{R_\zeta}{V}
+ \left (\frac{\frac23 + \frac Q 2 + \lambda_*}{V} - \frac \alpha y \right) \zeta
\in L^2(0, 1),
\]
the origin equation and its solution are
\[
\zeta' - \alpha \frac{\zeta}{y} = g_0 \in L^2 (0, 1)
\qquad \text{and} \qquad
\zeta (y) = c_{\mathrm{orig}} y^\alpha - y^\alpha \int_y^1 t^{- \alpha} g_0 (t) d t.
\]
It follows from Hardy's inequality 
\[
\left \| y^{\alpha - 1} \int_y^1 t^{- \alpha} g_0(t) \, d t \right \|_{L^2(0, 1)}
\leq \frac{\| g_0 \|_{L^2(0, 1)}}{\alpha - \frac12},
\]
and consequently 
\[
\left \| \frac{\zeta}{y} \right \|_{L^2(0, 1)} + \| \zeta' \|_{L^2(0, 1)} \leq C_\alpha ( |c_{\mathrm{orig}}| + \| g_0 \|_{L^2(0, 1)}). 
\]
Hence $\zeta \in H^1(- 1, 1)$ and $\zeta (0) = 0$.

\emph{3. Regularity at infinity.} Let $\beta = 2 + 3 \lambda_*$. Then
\[
\frac{\frac23 + \frac Q 2 + \lambda_*}{V} - \frac{\beta}{y} = O(y^{- 3}),
\]
and
\[
r_\zeta = \frac{y R_\zeta}{V}
+ y \left (\frac{\frac23 + \frac Q 2 + \lambda_*}{V} - \frac{\beta}{y} \right) \zeta
\in L^2(1, \infty)
\]
because $y^{- 2 } w^{ \frac{1}{2}} $ is bounded. The tail equation becomes
\[
\zeta ' - \beta \frac{\zeta}{y} = \frac{r_\zeta}{y},
\]
where $\zeta (y) = c_\infty y^\beta - y^\beta \int_y^\infty t^{- \beta - 1} r_\zeta (t) d t $, and 
\[
\left | y^\beta \int_y^\infty t^{- \beta - 1 } r_\zeta (t) d t \right | \leq \frac{ y^{ - \frac{1}{2} }}{\sqrt{ 2 \beta + 1}} \| r_\zeta \|_{L^2 (y, \infty )} = o (y^{- \frac{1}{2}}). 
\]
Since $ \beta > \sigma - \frac{1}{2} $ and $ \zeta \in L^2 (w^{- 1})$, we must have $c_\infty = 0$. Hardy's inequality bounds the tail integral, and the differential equation bounds the derivative.
\[
\| \zeta \|_{L^2 (1, \infty )} \leq \frac{ \| r_\zeta \|_{L^2(1, \infty)} }{\beta + \frac{1}{2} }, \, \, \| y \zeta' \|_{L^2(1, \infty)} \leq | \beta | \| \zeta \|_{L^2(1, \infty)} + \| r_\zeta \|_{L^2(1, \infty)}.
\]
\emph{4. The potential observable.} From the bounds $\zeta \in H^1$ and $y \zeta' \in L^2$,
\[
\| \zeta' \|_1
\leq \| \langle y \rangle^{- 1} \|_2 \| \langle y \rangle \zeta' \|_2 < \infty, \quad 
\lim_{y \to \pm \infty} \zeta (y) = 0,
\quad \int \zeta' = 0,
\quad \zeta'(- y) = \zeta'(y).
\]
Apply Lemma~\ref{1.20}, with $1/2 < \eta < 1$.
\[
\| \mathcal H \zeta ' \|_1 \leq \| \langle y \rangle^{- \eta } \|_2  \| \langle y \rangle^\eta \mathcal H \zeta' \|_2 \leq C_\eta \|  \langle y \rangle \zeta '  \|_2 < \infty. 
\]
Finally, $h = - \mathcal H(\mathcal P h)'$, $\mathcal H \zeta \in H^1$, and $\mathcal P h \in L^\infty$. We then have
\[
\ell (h) = \int \zeta h
= \int \mathcal H \zeta (\mathcal P h)'
= - \int (\mathcal H \zeta)' \mathcal P h.
\]
Here
\[
\lim_{y \to \pm \infty}(\mathcal H \zeta)(y) \mathcal P h(y) = 0,
\qquad
| \ell (h)| \leq \|(\mathcal H \zeta)' \|_1 \| \mathcal P h \|_\infty,
\]
which proves \eqref{2.18}.
\end{proof}

\bigskip

\section{Nonlinear analysis and blowup}

\subsection{Energy estimates}
Let $s$ and $J$ be smooth odd functions satisfying $s'(0) = 0 $, and assume that the norms below are finite. Define
\[
E_4^2 = \| s^{(4)} \|_{L^2 (\psi)}^2 + \| J^{(4)} \|^2_{L^2 (\psi )}
\]
and 
\[
E_{\text{raw}} = \| s \|_{L^2 (\phi)} + \| J \|_{L^2 (\varrho)} + E_4. 
\]
Since $w \leq C \phi$, it follows from Lemma~\ref{1.21} 
\[
\| s \|_{H^4 } + \| J \|_{H^4} + \| \mathcal P s \|_\infty + \| \mathcal P J \|_\infty \leq C E_{\text{raw}}
\]
\begin{equation} \label{3.3}
| \mathcal H s(0) | \leq C \| s \|_{L^2 (\phi)}, \, \, | q_J | \leq C \| J \|^2_{H^2} \leq C E^2_{\text{raw}},
\end{equation}
Indeed, by Cauchy--Schwarz,
\[
| \mathcal H s(0)|
\leq \pi^{- 1} \| s \|_{L^2(\phi)}
\left (\int \frac{d y}{y^2 \phi (y)} \right)^{1/2},
\]
and the integral is finite. The bound for $q_J$ follows from the $H^2$ trace estimates and $F'(0) \neq 0$.

\begin{Lemma}\label{3.1}
Let $s$ be smooth, real and odd, with $s'(0) = 0$ and $\| s \|_{L^2(\phi)} + \| s^{(4)} \|_{L^2(\psi)} < \infty$. Then the cutoff quadratic form satisfies $\langle \mathcal L_s s, s \phi \rangle \leq - \frac12 \| s \|_{L^2(\phi)}^2$. 
\end{Lemma}

\begin{proof}
First let $s \in C_c^\infty (\R)$ satisfy these conditions. We transfer the inviscid estimate in \cite[Equation 2.20]{Chen2020} to the present normalization. For $x = b y$, $b > 0$, $b^2 = \frac38$, the profile and weight become
\[
- \frac{2 b x }{ ( x^2 + b^2)^2 } \bigg |_{x  = b y} = - \frac{2 y}{ b^2 (1 + y^2)^2 } = F(y),
\qquad
\frac{ (x^2 + b^2 )^3}{2 b x^4} \bigg |_{x  = b y} = \frac{b}{2} \phi (y ).
\]
The scalar potential in \cite{Chen2020} satisfies
\[
\frac1 b \frac{x}{x^2 + b^2} \bigg |_{x = b y}
= \frac{y}{b^2(1 + y^2)} = U(y).
\]
For $q_J = 0$, the normalization corrections are
$- \mathcal H s(0)$ and $- \frac12 \mathcal H s(0)$. Hence
\[
c_A = - 1 - \mathcal H s(0), \qquad
c_\ell = \frac13 - \frac12 \mathcal H s(0).
\]
Thus the normalized linear operator in \cite{Chen2020} becomes \eqref{1.14}. The preserved condition $s'(0) = 0$ is precisely the hypothesis of that estimate. After cancelling the common Jacobian factor $b$,
\[
\frac{b}{2} \langle \mathcal L_s s, s \phi \rangle \leq - \frac{1}{2} \frac{b}{2} \| s \|^2_{L^2 (\phi)}. 
\]
Cancelling $b/2 > 0$ proves the estimate for compactly supported $s$.

For the stated class, define the pairing after integration by parts by
\[
\begin{split}
\langle \mathcal L_s s, s \phi \rangle
= {}  & \int \left (- 1 + Q + \frac{(V \phi)'}{2 \phi} \right) s^2 \phi \\
& + \int \left [F(\mathcal H s - \mathcal H s(0))
- \frac12 F'(\mathcal P s - y \mathcal H s(0)) \right] s \phi.
\end{split}
\]
Lemma~\ref{1.21} ensures $s \in H^1$. Hardy's inequalities and the rational coefficients imply
\[
\begin{gathered}
\left \| \frac{\mathcal H s - \mathcal H s(0)}{y} \right \|_2
+ \left \| \frac{\mathcal P s - y \mathcal H s(0)}{y^2} \right \|_2
\leq C \| s' \|_2, \\
\| Q \|_\infty + \left \| \frac{(V \phi)'}{\phi} \right \|_\infty
+ \| y \phi^{1/2} F \|_\infty
+ \| y^2 \phi^{1/2} F' \|_\infty < \infty.
\end{gathered}
\]
Thus the form is finite and continuous in $H^1 \cap L^2(\phi)$. Choose an even $\chi \in C_c^\infty (\R)$ equal to one near zero. Then
\[
s_R = \chi (y/R) s, \qquad
\| s_R - s \|_{H^1} + \| s_R - s \|_{L^2(\phi)} \longrightarrow 0.
\]
The estimates for $s_R$ pass to the limit. The transport cutoff flux vanishes since $V/y$ is bounded and $s^2 \phi$ is integrable.
\end{proof}

\begin{Lemma} \label{3.2}
The forcing in \eqref{1.17} satisfies $G_J (0) = G'_J (0) = 0 $ and $ \| G_J \|_{L^2 (\phi )} \leq C E^2_{\text{raw}}$. 
\end{Lemma}

\begin{proof}
Oddness implies $G_J(0) = 0$, and 
\[
G_J' (0) = - \frac{1}{2} \mathcal H J(0) J'(0) + \frac{1}{2} F'(0) q_J = 0 . 
\]
On $(- 1, 1)$, the product rule, \eqref{3.3} and Lemma~\ref{1.21} bound the forcing.
\[
\| G_J \|_{H^2(- 1, 1)}
\leq C \left (\| J \|_{H^3}^2
+ \| B \|_\infty \| J^{(3)} \|_2 + |q_J| \right)
\leq C E_{\text{raw}}^2.
\]
The two vanishing traces allow the second-order Hardy inequality.
\[
\left \| \frac{G_J}{y^2} \right \|_{L^2(- 1, 1)}
\leq C \| G_J'' \|_{L^2(- 1, 1)} \leq C E_{\text{raw}}^2.
\]
On $|y| \geq 1$, $|y| \leq \phi^{\frac12}(y) \leq 2 \sqrt2 \,|y|$ and $\mathfrak a_0 + \mathfrak b_0 = \frac{43}{32}$, so
\[
\begin{split}
\| |y|J \mathcal H J \|_{L^2(|y| \geq 1)}
& \leq \| \langle y \rangle^{\mathfrak b_0} J \|_\infty
\| \langle y \rangle^{\mathfrak a_0} \mathcal H J \|_2 
 \leq C \| \langle y \rangle^{\mathfrak b_0} J \|_\infty
\| \langle y \rangle^{\mathfrak a_0} J \|_2 \leq C E_{\text{raw}}^2.
\end{split}
\]
Similarly, since $\frac{1}{12} + \mathfrak a_1 = \frac{49}{48}$,
\[
\|    |y|  B J' \|_{L^2( |y| \geq 1 )} \leq \| |y |^{\frac{1}{12}} B \|_{L^\infty (|y| \geq 1 )} \| \langle y \rangle^{\mathfrak a_1} J' \|_2 \leq  C E^2_{\text{raw}}.
\]
The rational-profile term decays as $O (|y|^{- 3})$ and has coefficient $O ( E^2_{\text{raw}})$. Hence
\[
\| G_J \|_{L^2 (\phi )} \leq C E^2_{\text{raw}}. 
\]
\end{proof} 

\emph{Zeroth order energies.}
All integrals below are over $\R$. Assume that $s, J$ are real and $E_{\text{raw}} \leq 1$. The integrations by parts are understood with the spatial cutoffs justified in the proof of Proposition~\ref{3.21}. Since
\[
\frac{y \phi'}{\phi} = \frac{6 y^2}{1 + y^2} - 4, \qquad \mathcal P s(y) - y \mathcal H s(0) = \int_0^y \bigl (\mathcal H s(x) - \mathcal H s(0) \bigr) \, d x,
\]
we have
\[
\begin{split}
\left \| \frac{y \phi'}{\phi} \right \|_\infty  & \leq 4, \qquad
\left \| \frac{\mathcal P s - y \mathcal H s(0)}{y} \right \|_\infty
\leq \| \mathcal H s - \mathcal H s(0) \|_\infty \leq C E_{\text{raw}}, \\
\left \| \frac{(\mathcal P s - y \mathcal H s(0)) \phi'}{\phi} \right \|_\infty
& \leq 4 \| \mathcal H s - \mathcal H s(0) \|_\infty \leq C E_{\text{raw}}.
\end{split}
\]
Multiplying \eqref{1.16} by $s \phi$ and integrating by parts,
\[
\begin{split}
\frac12 \frac{d}{d \tau} \| s \|_{L^2( \phi )}^2
= {}  & \int (\mathcal L_s s) s \phi
+ \frac14 \int \bigl [(\mathcal P s - y \mathcal H s(0)) \phi \bigr]' s^2 \\
& + \int \bigl (\mathcal H s - \mathcal H s (0) \bigr) s^2 \phi
+ q_J \int \left (\phi + \frac14(y \phi)' \right) s^2
+ \int G_J s \phi.
\end{split}
\]
By \eqref{3.3} and Lemmas \ref{3.1}, \ref{3.2},
\[
\begin{split}
\left | \frac14 \int \bigl [(\mathcal P s - y \mathcal H s(0)) \phi \bigr]' s^2 \right |
& \leq C E_{\text{raw}} \| s \|_{L^2(\phi)}^2, \\
\left | \int \bigl (\mathcal H s - \mathcal H s(0) \bigr) s^2 \phi \right |
& \leq C E_{\text{raw}} \| s \|_{L^2(\phi)}^2, \\
\left |q_J \int \left (\phi + \frac14(y \phi)' \right) s^2 \right |
& \leq C E_{\text{raw}}^2 \| s \|_{L^2(\phi)}^2, \\
\left | \int G_J s \phi \right |
& \leq \| G_J \|_{L^2(\phi)} \| s \|_{L^2(\phi)}
\leq C E_{\text{raw}}^2 \| s \|_{L^2(\phi)}.
\end{split}
\]
Therefore,
\begin{equation} \label{3.4}
\begin{split}
\frac{1}{2} (\| s \|_{L^2 (\phi)}^2)'
& \leq - \frac{1}{2} \| s \|_{L^2 (\phi)}^2
+ C E_{\text{raw}} \| s \|_{L^2 (\phi)}^2
+ C E_{\text{raw}}^2 \| s \|_{L^2 (\phi)} \\
& \leq - \frac{1}{2} \| s \|_{L^2 (\phi)}^2 + C E^3_{\text{raw}}.
\end{split}
\end{equation}
From \eqref{1.18} and \eqref{1.19},
\[
\begin{split}
\frac12 \frac{d}{d \tau} \| J \|_{L^2(\varrho)}^2
= {}  & \int (\mathcal L J) J \varrho
+ \frac12 \int \left [\left (\delta_\ell y + \frac12 \mathcal P s \right) \varrho \right]' J^2 \\
& + \int (\delta_A - \mathcal H s) J^2 \varrho
+ \int s \mathcal H J \, J \varrho
+ \frac12 \int B s' J \varrho.
\end{split}
\]
The weight in Proposition~\ref{2.15} satisfies
\[
\left \| \frac{y \varrho'}{\varrho} \right \|_\infty < \infty,
\qquad
\left \| \frac{\mathcal P s}{y} \right \|_\infty \leq \| \mathcal H s \|_\infty.
\]
Thus, by \eqref{3.3},
\[
| \delta_A| + | \delta_\ell | + \| \mathcal H s \|_\infty
\leq C(E_{\text{raw}} + E_{\text{raw}}^2) \leq C E_{\text{raw}},
\]
\[
\left \| \frac{\left [\left (\delta_\ell y + \frac12 \mathcal P s \right) \varrho \right]'}{\varrho} \right \|_\infty
\leq | \delta_\ell | + \frac12 \| \mathcal H s \|_\infty
+ \left (| \delta_\ell | + \frac12 \left \| \frac{\mathcal P s}{y} \right \|_\infty \right)
\left \| \frac{y \varrho'}{\varrho} \right \|_\infty
\leq C E_{\text{raw}}.
\]
By Lemmas \ref{1.20}, \ref{1.21} and the bounds $c \langle y \rangle^{2 \sigma} \leq \varrho (y) \leq C \langle y \rangle^{2 \sigma}$, with fixed constants $c, C > 0$,
\[
\left | \int s \mathcal H J \, J \varrho \right |
\leq \| s \|_\infty \| \mathcal H J \|_{L^2(\varrho)} \| J \|_{L^2(\varrho)}
\leq C E_{\text{raw}} \| J \|_{L^2(\varrho)}^2,
\]
\[
\begin{split}
\| \varrho^{\frac12} B s' \|_2
& \leq C \| J \|_X \| \langle y \rangle^{\sigma - \frac1{12}} s' \|_2 = C \| J \|_X \| \langle y \rangle^{\frac12} s' \|_2
\leq C E_{\text{raw}} \| J \|_{L^2(\varrho)}, \\
\left | \int B s' J \varrho \right |
& \leq \| \varrho^{\frac12} B s' \|_2 \| J \|_{L^2(\varrho)}
\leq C E_{\text{raw}} \| J \|_{L^2(\varrho)}^2.
\end{split}
\]
Combining the preceding bounds with Proposition~\ref{2.15} yields
\begin{equation} \label{3.5}
\frac{1}{2} (\| J \|_{L^2 (\varrho)}^2)'
\leq - \left (\frac{6}{25} - C E_{\text{raw}} \right) \| J \|_{L^2(\varrho)}^2.
\end{equation}
\emph{Fourth-order linear damping.}
For $f = s, J$, differentiation and one integration by parts give
\[
\partial_y^4(- V f')
= - V f^{(5)} - 4 V' f^{(4)}
- \sum_{j = 2}^4 \binom4 j V^{(j)} f^{(5 - j)},
\]
\[
- \int \psi V f^{(4)} f^{(5)}
= \frac12 \int (\psi V)'|f^{(4)}|^2.
\]
Since
\[
V' = \frac13 + \frac Q 2,
\qquad \frac{\psi'}{\psi} = \frac{4 y}{1 + y^2},
\]
the coefficients of $|s^{(4)}|^2 \psi$ and $|J^{(4)}|^2 \psi$ are
\begin{equation} \label{3.6}
d_{s, 4}
= - 1 + Q - \frac72 V' + \frac{V \psi'}{2 \psi}
= - \frac{13}{6} - \frac34 Q + \frac{2 V y}{1 + y^2}
= - \frac9{10} - \frac{(3 y^2 - 7)^2}{15(1 + y^2)^2},
\end{equation}
\begin{equation} \label{3.7}
d_{J, 4}
= - 1 - Q - \frac72 V' + \frac{V \psi'}{2 \psi}
= - \frac{13}{6} - \frac{11}{4} Q + \frac{2 V y}{1 + y^2}
= - \frac{19}{78} - \frac{(7 y^2 - 19)^2}{39(1 + y^2)^2}.
\end{equation}
For real $g \in C_c^\infty (\R)$,
\[
\int y g(y) \mathcal H g(y) \, d y
= \lim_{\epsilon \to 0} \frac1{2 \pi}
\iint_{|y - z| > \epsilon} \frac{y - z}{y - z} g(y) g(z) \, d z \, d y
= \frac1{2 \pi} \left (\int g \right)^2.
\]
If $\langle y \rangle g \in L^2$, choose real $g_n \in C_c^\infty$ with
$\| \langle y \rangle (g_n - g) \|_2 \to 0$. Then
\[
\begin{split}
\left | \int y g_n \mathcal H g_n - \int y g \mathcal H g \right |
& \leq \| y(g_n - g) \|_2 \| g_n \|_2 + \| y g \|_2 \| g_n - g \|_2 \to 0, \\
\left | \int (g_n - g) \right |
& \leq \| \langle y \rangle^{- 1} \|_2 \| \langle y \rangle (g_n - g) \|_2 \to 0.
\end{split}
\]
Consequently,
\begin{equation} \label{3.8}
\int y g(y) \mathcal H g(y) \, d y = \frac1{2 \pi} \left (\int g \right)^2.
\end{equation}
Since $f^{(4)}$ is odd and $\psi F = - 16 y/3$,
\[
\int \psi F f^{(4)} \mathcal H f^{(4)}
= - \frac{16}{3} \int y f^{(4)} \mathcal H f^{(4)}
= - \frac8{3 \pi} \left (\int f^{(4)} \right)^2 = 0.
\]
To estimate the remaining linear terms, we use
\[
\partial_y^4(Q f) = Q f^{(4)} + \sum_{j = 1}^4 \binom4 j Q^{(j)} f^{(4 - j)},
\]
\[
\partial_y^4(F \mathcal H f)
= F \mathcal H f^{(4)} + \sum_{j = 1}^4 \binom4 j F^{(j)} \mathcal H f^{(4 - j)},
\]
\[
\partial_y^4(F' \mathcal P f)
= F^{(5)} \mathcal P f + \sum_{j = 0}^3 \binom4 j F^{(j + 1)} \mathcal H f^{(3 - j)}.
\]
The coefficients satisfy
\[
\begin{split}
& \sum_{j = 2}^4 \| \psi^{\frac12} V^{(j)} \|_\infty
+ \sum_{j = 1}^4 \left (\| \psi^{\frac12} Q^{(j)} \|_\infty
+ \| \psi^{\frac12} F^{(j)} \|_\infty \right) \leq C, \\
& \| \sqrt{|y|} \psi^{\frac12} F^{(5)} \|_2
+ \left \| \psi^{\frac12} \partial_y^4 \left (- F + \frac12 y F' \right) \right \|_2 \leq C.
\end{split}
\]
For $0 \leq j \leq 3$ and $\eta > 0$,
\[
\| \mathcal H f^{(j)} \|_2 = \| f^{(j)} \|_2
\leq C \| f \|_2^{1 - \frac j 4} \| f^{(4)} \|_2^{\frac j 4},
\]
\[
\begin{split}
\| \psi^{\frac12} F^{(5)} \mathcal P f \|_2
& \leq \| \sqrt{|y|} \psi^{\frac12} F^{(5)} \|_2 \| f \|_2, \\
\| f^{(4)} \|_{L^2(\psi)} \| f \|_2^{1 - \frac j 4} \| f^{(4)} \|_2^{\frac j 4}
& \leq \eta \| f^{(4)} \|_{L^2(\psi)}^2 + C_\eta \| f \|_2^2, \\
| \mathcal H s(0)| \| s^{(4)} \|_{L^2(\psi)}
& \leq \eta \| s^{(4)} \|_{L^2(\psi)}^2 + C_\eta \| s \|_{L^2(\phi)}^2.
\end{split}
\]
Thus, for every $\delta > 0$,
\begin{equation} \label{3.9}
\begin{split}
& \left | \int \psi s^{(4)} \partial_y^4(\mathcal L_s s)
- \int d_{s, 4} \psi |s^{(4)}|^2 \right |
+ \left | \int \psi J^{(4)} \partial_y^4(\mathcal L J)
- \int d_{J, 4} \psi |J^{(4)}|^2 \right | \\
& \leq \delta E_4^2 + C_\delta \left (\| s \|_{L^2(\phi)}^2 + \| J \|_{L^2(\varrho)}^2 \right).
\end{split}
\end{equation}
\emph{Fourth-order nonlinear terms.}
The product terms follow from the weighted bounds in Lemma~\ref{1.21}. For $f, g \in \{ s, J \}$ and $1 \leq j \leq 3$,
\[
\begin{split}
\| f^{(j)} \mathcal H g^{(4 - j)} \|_{L^2(\psi)}
& \leq \| \langle y \rangle^{\mathfrak b_j} f^{(j)} \|_\infty
\| \langle y \rangle^{\mathfrak a_{4 - j}} \mathcal H g^{(4 - j)} \|_2
\leq C E_{\text{raw}}^2, \\
\mathfrak b_j + \mathfrak a_{4 - j}  & = \frac{265}{96} > 2.
\end{split}
\]
The endpoint terms satisfy
\[
\begin{aligned}
\| f \mathcal H g^{(4)} \|_{L^2(\psi)}
+ \| f^{(4)} \mathcal H g \|_{L^2(\psi)}
& \leq \| f \|_\infty \| \mathcal H g^{(4)} \|_{L^2(\psi)}
\\
& \quad + \| f^{(4)} \|_{L^2(\psi)} \| \mathcal H g \|_\infty
\leq C E_{\text{raw}}^2.
\end{aligned}
\]
Hence
\[
\| \partial_y^4(f \mathcal H g) \|_{L^2(\psi)}
\leq \sum_{j = 0}^4 \binom4 j \| f^{(j)} \mathcal H g^{(4 - j)} \|_{L^2(\psi)}
\leq C E_{\text{raw}}^2.
\]
For the diagonal transport terms,
\[
\begin{split}
& \int \psi f^{(4)} \partial_y^4 \left [- \left (\delta_\ell y + \frac12 \mathcal P s \right) f' + \delta_A f \right] \\
& = \int \left [\delta_A - \frac72 \left (\delta_\ell + \frac12 \mathcal H s \right)
+ \frac{\psi'}{2 \psi} \left (\delta_\ell y + \frac12 \mathcal P s \right) \right] \psi |f^{(4)}|^2 \\
& \quad - \frac12 \sum_{j = 2}^4 \binom4 j \int \psi f^{(4)}(\mathcal H s^{(j - 1)}) f^{(5 - j)}.
\end{split}
\]
Here
\[
\begin{split}
& \left \| \delta_A - \frac72 \left (\delta_\ell + \frac12 \mathcal H s \right)
+ \frac{\psi'}{2 \psi} \left (\delta_\ell y + \frac12 \mathcal P s \right) \right \|_\infty
\leq C E_{\text{raw}}, \\
& \| (\mathcal H g^{(j - 1)}) f^{(5 - j)} \|_{L^2(\psi)}
\leq \| \langle y \rangle^{\mathfrak b_{j - 1}} \mathcal H g^{(j - 1)} \|_\infty
\| \langle y \rangle^{\mathfrak a_{5 - j}} f^{(5 - j)} \|_2
\leq C E_{\text{raw}}^2, 2 \leq j \leq 4.
\end{split}
\]
For magnetic transport, use $B^{(j)} = \mathcal H J^{(j - 1)}$. The two fifth-derivative terms combine before integration by parts.
\[
\frac12 \int \psi B \left (s^{(4)} J^{(5)} + J^{(4)} s^{(5)} \right)
= \frac12 \int \psi B \left (s^{(4)} J^{(4)} \right)'
= - \frac12 \int (\psi B)' s^{(4)} J^{(4)}.
\]
Consequently,
\[
\begin{split}
& \frac12 \int \psi \left [s^{(4)} \partial_y^4(B J') + J^{(4)} \partial_y^4(B s') \right]
= \int \left (\frac72 \psi B' - \frac12 \psi' B \right) s^{(4)} J^{(4)} \\
& \quad + \frac12 \sum_{j = 2}^4 \binom4 j
\int \psi B^{(j)} \left (s^{(4)} J^{(5 - j)} + J^{(4)} s^{(5 - j)} \right).
\end{split}
\]
Since
\[
\left \| \frac72 B' - \frac12 \frac{\psi'}{\psi} B \right \|_\infty
\leq C \left (\| \mathcal H J \|_\infty + \| B \|_\infty \right)
\leq C E_{\text{raw}},
\]
we obtain
\[
\left | \frac12 \int \psi \left [s^{(4)} \partial_y^4(B J') + J^{(4)} \partial_y^4(B s') \right] \right |
\leq C E_{\text{raw}} E_4^2 + C E_{\text{raw}}^2 E_4
\leq C E_{\text{raw}}^3.
\]
The remaining normalization term satisfies
\[
\begin{split}
& \left |q_J \int \psi s^{(4)} \partial_y^4 \left (F - \frac12 y F' \right) \right | \\
& \qquad \leq |q_J| \| s^{(4)} \|_{L^2(\psi)}
\left \| \partial_y^4 \left (F - \frac12 y F' \right) \right \|_{L^2(\psi)}
\leq C E_{\text{raw}}^2 E_4 \leq C E_{\text{raw}}^3.
\end{split}
\]
Together with \eqref{3.6}, \eqref{3.7}, and \eqref{3.9},
\[
\begin{split}
\frac12(E_4^2)'
& \leq - \frac9{10} \| s^{(4)} \|_{L^2(\psi)}^2
- \frac{19}{78} \| J^{(4)} \|_{L^2(\psi)}^2 + \delta E_4^2 \\
& \quad + C_\delta \left (\| s \|_{L^2(\phi)}^2 + \| J \|_{L^2(\varrho)}^2 \right)
+ C E_{\text{raw}}^3.
\end{split}
\]
Take $\delta = \frac{53}{7800} < \frac{19}{78} - \frac{23}{100}$ and $C_1 = C_{\frac{53}{7800}}$. Then
\begin{equation} \label{3.10}
\frac12(E_4^2)'
\leq - \frac{23}{100} E_4^2
+ C_1 \left (\| s \|_{L^2(\phi)}^2 + \| J \|_{L^2(\varrho)}^2 \right)
+ C E_{\text{raw}}^3.
\end{equation}
The constants $C$ and $C_1$ depend only on the fixed profiles and weights, and are independent of the solution and of $\tau$.

\begin{Proposition}\label{3.11}
There exist $\mu > 0$ and $E_* > 0$ with the following property. Define
\[
E^2 = \| s \|_{L^2(\phi)}^2 + \| J \|_{L^2(\varrho)}^2 + \mu E_4^2.
\]
Every smooth normalized solution considered above with $E(0) < E_*$ satisfies, throughout its interval of existence,
\begin{equation} \label{3.12}
E (\tau ) \leq E (0) e^{- \frac{9 \tau }{40}}.
\end{equation}
Therefore, 
\[
\| J \|_X + \| J \|_{H^2} \leq C E (0) e^{- \frac{9 \tau}{40}}. 
\]
\end{Proposition}

\begin{proof}
Fix $\mu > 0$ with $\mu C_1 < \frac1{100}$. Since
\[
\sqrt{\min \{ 1, \mu \}/3} \, E_{\text{raw}}
\leq E \leq \sqrt{\max \{ 1, \mu \}} \, E_{\text{raw}},
\]
by \eqref{3.4}, \eqref{3.5} and \eqref{3.10}, there is a constant $C$, fixed after $\mu$, such that
\[
\frac12(E^2)' \leq - \frac{23}{100} E^2 + C E^3
\leq - \frac9{40} E^2
\]
whenever $E \leq E_*$ and $C E_* \leq \frac1{200}$.
Decrease $E_*$ so that all preceding estimates hold for $E \leq E_*$. On every such interval,
\[
E(\tau) \leq E(0) e^{- 9 \tau /40} < E_*.
\]
Continuity closes the bootstrap and proves \eqref{3.12}. For the magnetic norms, it follows from Lemma~\ref{1.21} 
\[
\| J \|_X + \| J \|_{H^2} \leq C E \leq C E(0) e^{- 9 \tau /40}.
\]
\end{proof}
Fix a sufficiently small $\mu_3 > 0 $ and define $M_3^2 = \| s \|_{L^2 (\phi)}^2 + \mu_3 \| s^{(3)} \|_{L^2 (\psi)}^2 $.

\begin{Proposition}\label{3.21}
For the solutions considered above, if $s(y, 0) \equiv 0$ and $E(0) \leq C_0 \varepsilon$, then
\begin{equation} \label{3.17}
M_3 (\tau) \leq C \varepsilon^2 e^{- \frac{9 \tau}{20}}, \qquad
\| s \|_{H^1} + \| s \|_{L^2 (\phi)} \leq C M_3.
\end{equation}
\end{Proposition}

\begin{proof}
\emph{1. Weighted regularity before damping.} Fix a compact normalized interval of existence $I = [0, \tau_0]$ for the constructed solution and its physical interval $I_t = t(I)$. In this step, $\omega, J, u, B$ denote the physical fields. By Lemma~\ref{3.15},
\[
(\omega, J) \in C(I_t; Z \times Z).
\]
Differentiating \eqref{1.2} three times,
\[
\| \partial_t \omega^{(3)} \|_{L^2(\psi)} \leq C \|(\omega, J) \|_Z^2.
\]
Indeed, the principal terms satisfy
\[
\| u \omega^{(4)} \|_{L^2(\psi)} + \| B J^{(4)} \|_{L^2(\psi)}
\leq C \|(\omega, J) \|_Z^2.
\]
For $f, g \in \{ \omega, J \}$ and $0 \leq j \leq 3$, the remaining terms satisfy
\[
\begin{aligned}
\| f^{(j)} \mathcal H g^{(3 - j)} \|_{L^2(\psi)}
& \leq \| \langle y \rangle^{9/8 + j/4} f^{(j)} \|_\infty
\| \langle y \rangle^{1 + (3 - j)/4} \mathcal H g^{(3 - j)} \|_2
\leq C \| f \|_Z \| g \|_Z, \\
\frac98 + \frac j 4 + 1 + \frac{3 - j}{4}  & = \frac{23}{8} > 2.
\end{aligned}
\]
These bilinear estimates imply that
\[
(\omega, J) \longmapsto \partial_x^3 \left [- \frac12(\mathcal P \omega) \omega_x + \frac12(\mathcal P J) J_x + \omega \mathcal H \omega - J \mathcal H J \right]
\]
is continuous from the closed odd subspace of $Z \times Z$ to $L^2(\psi)$. Since $\omega^{(3)}( \cdot , 0) = F^{(3)} \in L^2(\psi)$,
\[
F^{(3)} + \int_0^t \partial_x^3 \left [- \frac12 u \omega_x + \frac12 B J_x + \omega \mathcal H \omega - J \mathcal H J \right]( \cdot , r) \, d r
\]
belongs to $C(I_t; L^2(\psi))$ and equals $\omega^{(3)}$ in distributions. On $I$,
\[
0 < \min_I \{ A, L \} \leq \max_I \{ A, L \} < \infty.
\]
Under dilation,
\[
\| \partial_y^3(A f(L \cdot )) \|_{L^2(\psi)}^2
\leq A^2 L^5 \max \{ 1, L^{- 4} \} \| f^{(3)} \|_{L^2(\psi)}^2
\]
Dilation is continuous in $L^2(\psi)$ by compactly supported approximation, hence $s^{(3)} \in C(I; L^2(\psi))$. Using the preserved traces,
\[
\begin{gathered}
s(0, \tau) = s'(0, \tau) = 0, \qquad
\| s/y^2 \|_{L^2(- 1, 1)} \leq C \| s'' \|_{L^2(- 1, 1)}, \\
\langle y \rangle^2 \leq \phi (y) \leq 4 \langle y \rangle^2
\quad (|y| \geq 1).
\end{gathered}
\]
Applying these bounds to $s(\tau) - s(r)$,
\[
\lim_{\tau \to r} \| s(\tau) - s(r) \|_{L^2(\phi)} = 0,
\qquad M_3 \in C(I).
\]

\emph{2. Cutoff justification.} We justify the weighted identities before using their signs. Choose an even $\chi \in C_c^\infty (\R)$ with $0 \leq \chi \leq 1$, $\chi = 1$ on $[- 1, 1]$ and $\chi = 0$ outside $[- 2, 2]$, and set
\[
\eta_{\epsilon, R}(y) = \bigl (1 - \chi (y/ \epsilon) \bigr) \chi (y/R),
\qquad 0 < \epsilon < \tfrac12, \quad R > 2.
\]
Test the scalar equation and its third derivative with
\[
\eta_{\epsilon, R} s \phi, \qquad
\eta_{\epsilon, R} s^{(3)} \psi.
\]
Integrate on $I$ before removing the cutoffs. The arguments of $\mathcal H$ and $\mathcal P$ remain $s$ and $J$.

The transport velocities satisfy
\[
\begin{gathered}
\sup_{\tau \in I} \| v/y \|_\infty \leq C_I, \qquad
v \in \left \{ V, \  \frac12(\mathcal P s - \mathcal H s(0) y + q_J y), \,
\delta_\ell y + \frac12 \mathcal P s, \  B \right \}, \\
\| \mathcal P s/y \|_\infty \leq \| \mathcal H s \|_\infty,
\qquad \| B/y \|_\infty \leq \| \mathcal H J \|_\infty.
\end{gathered}
\]
Consequently,
\[
| \eta_{\epsilon, R}'(y) v(y, \tau)|
\leq C_I \left (\mathbf1_{\{ \epsilon \leq |y| \leq 2 \epsilon \}}
+ \mathbf1_{\{ R \leq |y| \leq 2 R \}} \right).
\]
Set
\[
\mathcal A_{\epsilon, R}
= \{ \epsilon \leq |y| \leq 2 \epsilon \} \cup \{ R \leq |y| \leq 2 R \}.
\]
The scalar boundary fluxes satisfy
\[
\int_I \! \int | \eta_{\epsilon, R}' v|
\bigl (s^2 \phi + |s^{(3)}|^2 \psi \bigr) \, d y \, d \tau \leq C_I \int_I \! \int_{\mathcal A_{\epsilon, R}}
\bigl (s^2 \phi + |s^{(3)}|^2 \psi \bigr) \, d y \, d \tau
\longrightarrow 0
\]
as $\epsilon \to 0$ and $R \to \infty$. For the Hilbert term,
\[
\int_I \! \int | \psi F s^{(3)} \mathcal H s^{(3)}| \, d y \, d \tau
\leq \frac{16}{3} \int_I \| y s^{(3)} \|_2 \| s^{(3)} \|_2 \, d \tau < \infty.
\]
Thus its cutoff limit is the original Hilbert pairing. For the remaining terms, the weighted bounds below reduce convergence to
\[
f, g \in L^2(I \times \R), \qquad \int_I \! \int |(1 - \eta_{\epsilon, R}) f g|
\leq \|(1 - \eta_{\epsilon, R}) f \|_{L^2(I \times \R)}
\| g \|_{L^2(I \times \R)} \longrightarrow 0.
\]

For the fourth order estimates, use the same $\eta_{\epsilon, R}$ in both equations. The boundary fluxes satisfy
\[
\begin{aligned}
& \int_I \! \int | \eta_{\epsilon, R}' \psi B s^{(4)} J^{(4)}| \, d y \, d \tau
+ \sum_{f \in \{ s, J \}} \int_I \! \int
| \eta_{\epsilon, R}' v| \psi |f^{(4)}|^2 \, d y \, d \tau \\
& \qquad \leq C_I \int_I \! \int_{\mathcal A_{\epsilon, R}}
\psi \bigl (|s^{(4)}|^2 + |J^{(4)}|^2 \bigr) \, d y \, d \tau
\longrightarrow 0.
\end{aligned}
\]
The limits follow from $C(I; Z \times Z)$ and the scalar Hardy estimate in Step~1. For the constructed solution,
\[
\| s \|_{L^2(\phi)}^2, \quad
\| J \|_{L^2(\varrho)}^2, \quad E_4^2, \quad M_3^2
\in W^{1, 1}(I).
\]
The energy inequalities therefore hold almost everywhere. A pairing containing an unbounded transport derivative denotes its cutoff limit after integration by parts. This requires no $L^2(\psi)$ estimate for $\partial_\tau s^{(3)}$.

\emph{3. Linear damping.} The third-order coefficient satisfies
\[
d_{s, 3} = - \frac{11}{6} - \frac14 Q + \frac{2 V y}{1 + y^2}
= - \frac{13}{18} - \frac{4(y^2 - 2)^2}{9(1 + y^2)^2}.
\]
Equation \eqref{3.8} applies since $s^{(3)} (- y) = s^{(3)} (y) $ and $ \int s^{(3) } = 0 $.
\[
\int \psi F s^{(3)} \mathcal H s^{(3)} = - \frac{16}{3} \int y s^{(3)} \mathcal H s^{(3)} = 0 
\]
and 
\[
\| \langle y \rangle^2 \mathcal H s^{(3)} \|_2 \leq C \| \langle y \rangle^2 s^{(3)} \|_2.
\]
The third derivative of \eqref{1.14} is
\[
\begin{split}
\partial_y^3(\mathcal L_s s)
= {}  & - V s^{(4)} + (- 1 + Q - 3 V') s^{(3)} + F \mathcal H s^{(3)} \\
& + (3 Q' - 3 V'') s'' + (3 Q'' - V^{(3)}) s' + Q^{(3)} s
+ \frac52 F' \mathcal H s'' \\
& + \frac32 F'' \mathcal H s'
- \frac12 F^{(3)} \mathcal H s - \frac12 F^{(4)} \mathcal P s
+ \frac12 \mathcal H s(0) \left (F^{(3)} + y F^{(4)} \right).
\end{split}
\]
Moreover,
\[
- \int \psi V s^{(3)} s^{(4)}
= \frac12 \int (\psi V)'|s^{(3)}|^2,
\qquad
- 1 + Q - \frac52 V' + \frac{V \psi'}{2 \psi} = d_{s, 3}.
\]
The coefficients satisfy
\[
\begin{split}
& \sum_{j = 2}^3 \| \psi^{\frac12} V^{(j)} \|_\infty
+ \sum_{j = 1}^3 \left (\| \psi^{\frac12} Q^{(j)} \|_\infty
+ \| \psi^{\frac12} F^{(j)} \|_\infty \right) \leq C, \\
& \| \sqrt{|y|} \psi^{\frac12} F^{(4)} \|_2
+ \left \| \psi^{\frac12} \left (F^{(3)} + y F^{(4)} \right) \right \|_2 \leq C.
\end{split}
\]
Since
\[
\| \mathcal H s^{(j)} \|_2 = \| s^{(j)} \|_2, \qquad
| \mathcal P s(y)| \leq \sqrt{|y|} \| s \|_2, \qquad
| \mathcal H s(0)| \leq C \| s \|_{L^2(\phi)},
\]
the remaining linear terms, after the cancellation $\int \psi F s^{(3)} \mathcal H s^{(3)} = 0$, satisfy
\[
\left | \int \psi s^{(3)} \partial_y^3(\mathcal L_s s)
- \int d_{s, 3} \psi |s^{(3)}|^2 \right |
\leq C \| s^{(3)} \|_{L^2(\psi)}
\left (\sum_{j = 0}^2 \| s^{(j)} \|_2 + \| s \|_{L^2(\phi)} \right).
\]
For $0 \leq j \leq 2$,
\[
\| s^{(j)} \|_2
\leq C \| s \|_2^{1 - \frac j 3} \| s^{(3)} \|_2^{\frac j 3}
\leq C \| s \|_{L^2(\phi)}^{1 - \frac j 3}
\| s^{(3)} \|_{L^2(\psi)}^{\frac j 3}.
\]
For every $\delta > 0$, Young's inequality bounds the preceding linear error.
\[
C \| s^{(3)} \|_{L^2( \psi  )}
\left (\sum_{j = 0}^2 \| s^{(j)} \|_2 + \| s \|_{L^2( \phi )}   \right)
\leq \delta \| s^{(3)} \|_{L^2(\psi)}^2
+ C_\delta \| s \|_{L^2(\phi)}^2,
\]
where $C_\delta$ depends only on $\delta$ and the fixed profiles and weights. Since $d_{s, 3} \leq - \frac{13}{18}$,
\[
\int  \psi s^{(3)} \partial_y^3(\mathcal L_s s)
\leq - \left (\frac{13}{18} - \delta \right) \| s^{(3)} \|_{L^2(\psi)}^2 +    C_\delta \| s \|_{L^2(\phi)}^2.
\]
Since $\frac{49}{100} \in (\frac{12}{25}, \frac12)$, take
\[
\delta = \frac12 \left (\frac{13}{18} - \frac{49}{100} \right),
\qquad
0 < \mu_3 \leq \min \left \{ 1, \frac{\frac12 - \frac{49}{100}}{1 + C_\delta} \right \}.
\]
Then
\[
\frac{13}{18} - \delta \geq \frac{49}{100},
\qquad
\frac12 - \mu_3 C_\delta \geq \frac{49}{100}.
\]
By Lemma~\ref{3.1},
\[
\begin{split}
\frac12(M_3^2)'_{\text{linear}}
& = \int (\mathcal L_s s) s \phi
+ \mu_3 \int \psi s^{(3)} \partial_y^3(\mathcal L_s s) \\
& \leq - \left (\frac12 - \mu_3 C_\delta \right) \| s \|_{L^2(\phi)}^2
- \mu_3 \left (\frac{13}{18} - \delta \right) \| s^{(3)} \|_{L^2(\psi)}^2 \\
& \leq - \frac{49}{100} \left (\| s \|_{L^2(\phi)}^2 + \mu_3 \| s^{(3)} \|_{L^2(\psi)}^2 \right)
= - \frac{49}{100} M_3^2.
\end{split}
\]
\emph{4. Scalar nonlinear terms.} Define the scalar perturbation velocity by
\[
v_s = \frac12 \bigl (\mathcal P s - \mathcal H s(0) y + q_J y \bigr), \qquad
v_s' = \frac12 \bigl (\mathcal H s - \mathcal H s(0) + q_J \bigr).
\]
Then \eqref{1.16} becomes $s_\tau = \mathcal L_s s + N_s + G_J$, where $N_s = - v_s s' + (\mathcal H s - \mathcal H s(0) + q_J) s$. For $E \leq 1$, Lemma~\ref{1.21} and \eqref{3.3} imply
\[
\| \mathcal H s - \mathcal H s(0) + q_J \|_\infty + \| v_s' \|_\infty
+ \| v_s/y \|_\infty \leq C E.
\]
The quotient is extended continuously at zero. Since $y \phi'/ \phi$ and $y \psi'/ \psi$ are bounded, the zeroth-order contribution satisfies
\[
\int N_s s \phi
= \int \left (\mathcal H s - \mathcal H s(0) + q_J + \frac12 v_s' + \frac{v_s \phi'}{2 \phi} \right) s^2 \phi, \qquad \left | \int N_s s \phi \right | \leq C E \| s \|_{L^2(\phi)}^2.
\]
For the third derivative, the exact expansion is
\[
\begin{aligned}
\partial_y^3 N_s  & = - v_s s^{(4)} + (\mathcal H s - \mathcal H s(0) + q_J - 3 v_s') s^{(3)} \\
& + \frac32(\mathcal H s') s'' + \frac52(\mathcal H s'') s' + (\mathcal H s^{(3)}) s.
\end{aligned}
\]
Since $\phi \geq \langle y \rangle^2$,
\[
\| s \|_{H^3} \leq C_{\mu_3} M_3, \qquad
y s \in L^2, \qquad (y s)^{(3)} = y s^{(3)} + 3 s'' \in L^2.
\]
By Fourier interpolation,
\[
\| y s \|_{H^3} \leq C \bigl (\| y s \|_2 + \|(y s)^{(3)} \|_2 \bigr)
\leq C_{\mu_3} M_3.
\]
Since $y s' = (y s)' - s$ and $y s'' = (y s)'' - 2 s'$,
\[
\| \langle y \rangle s' \|_2 + \| \langle y \rangle s'' \|_2
\leq C_{\mu_3} M_3.
\]
Moreover, $\mathfrak b_1 = 107/96 > 1$ and $\mathfrak b_2 = 141/96 > 1$ in Lemma~\ref{1.21}, hence
\[
\| \langle y \rangle \mathcal H s' \|_\infty
+ \| \langle y \rangle \mathcal H s'' \|_\infty + \| s \|_\infty \leq C E.
\]
Here $s^{(3)}$ is even with zero integral, so Lemma~\ref{1.20} applies.
\[
\| \mathcal H s^{(3)} \|_{L^2(\psi)} \leq C \| s^{(3)} \|_{L^2(\psi)}
\leq C_{\mu_3} M_3.
\]
Consequently,
\[
\begin{split}
& \left \| \frac32(\mathcal H s') s'' + \frac52(\mathcal H s'') s' + (\mathcal H s^{(3)}) s \right \|_{L^2(\psi)} \\
& \leq \frac32 \| \langle y \rangle \mathcal H s' \|_\infty \| \langle y \rangle s'' \|_2
+ \frac52 \| \langle y \rangle \mathcal H s'' \|_\infty \| \langle y \rangle s' \|_2
+ \| s \|_\infty \| \mathcal H s^{(3)} \|_{L^2(\psi)} \leq C_{\mu_3} E M_3.
\end{split}
\]
After integrating the principal transport term by parts, its coefficient together with the zeroth-order term is
\[
\mathcal H s - \mathcal H s(0) + q_J - \frac52 v_s' + \frac{v_s \psi'}{2 \psi},
\]
whose $L^\infty$ norm is bounded by $C E$. In the cutoff-justified identity,
\[
\left | \int \psi s^{(3)} \partial_y^3 N_s \right |
\leq C E \| s^{(3)} \|_{L^2(\psi)}^2
+ C_{\mu_3} E M_3 \| s^{(3)} \|_{L^2(\psi)}.
\]
With $\mu_3 > 0$ fixed as above, this and the zeroth-order bound imply
\[
\left | \int N_s s \phi
+ \mu_3 \int \psi s^{(3)} \partial_y^3 N_s \right |
\leq C_{\mu_3} E M_3^2.
\]
\emph{5. Magnetic forcing.} The transport product satisfies
\[
\partial_y^3 (B J' ) = B J^{(4)} + \sum_{j = 1}^3 \binom{3}{j} B^{(j)} J^{(4 - j)},
\]
where $B^{(j)} = \mathcal H J^{(j - 1 )}$. The endpoint and the remaining products satisfy
\[
\| B J^{(4)} \|_{L^2 (\psi )} \leq \| B \|_\infty \| J^{(4)} \|_{L^2 (\psi)} \leq C E^2,
\]
 and 
\[
\| ( \mathcal H J^{(j )} ) J^{(3 - j)} \|_{L^2 (\psi )} \leq \| \langle y \rangle^{\mathfrak b_j} \mathcal H J^{(j)} \|_\infty \|   \langle y \rangle^{\mathfrak a_{3 - j}} J^{(3 - j )}        \|_2 \leq C E^2 
\]
where $0 \leq j \leq 3$ and $\mathfrak b_j + \mathfrak a_{3 - j} = 77/32 > 2$. Also,
\[
|q_J| \left \| \partial_y^3 \left (F - \frac12 y F' \right) \right \|_{L^2(\psi)}
\leq C E^2.
\]
Together with Lemma~\ref{3.2},
\[
\| G_J \|_{L^2(\phi)} + \| G_J^{(3)} \|_{L^2(\psi)} \leq C E^2,
\qquad
\left | \int G_J s \phi + \mu_3 \int \psi s^{(3)} G_J^{(3)} \right |
\leq C E^2 M_3.
\]
\emph{6. Closing the estimate.} Combining Steps~3 to~5,
\[
\frac12(M_3^2)' \leq - \frac{49}{100} M_3^2
+ C_{\mu_3} E M_3^2 + C E^2 M_3.
\]
Decrease the common smallness threshold, after fixing $\mu_3$, so that $C_{\mu_3} E \leq 1/100$. Then
\[
\frac{1}{2} ( M_3^2)' \leq - \frac{12}{25} M_3^2 + C E^2 M_3. 
\]
Consequently,
\begin{equation} \label{3.13}
M_3' \leq - \frac{12}{25} M_3 + C E^2. 
\end{equation}
At points where $M_3 = 0 $, regularize by $( M_3^2 + \delta^2)^{\frac{1}{2}}$ and let $\delta \to 0$. Since $M_3(0) = 0$, Proposition~\ref{3.11} and variation of constants give
\[
\begin{split}
M_3(\tau)  & \leq C \varepsilon^2 \int_0^\tau e^{- 12(\tau - r)/25} e^{- 9 r/20} \, d r \\
& = \frac{100 C}{3} \varepsilon^2 \left (e^{- 9 \tau /20} - e^{- 12 \tau /25} \right)
\leq C' \varepsilon^2 e^{- 9 \tau /20}.
\end{split}
\]
and
\[
\| s \|_{H^1} \leq C ( \| s \|_{L^2 (\phi)} + \|   s^{(3)} \|_{L^2 (\psi)}) \leq C M_3. 
\]
\end{proof}

\subsection{Local existence and continuation}\label{3.14}
We first construct the physical solution and then solve for the rescaling factors. For arbitrary functions, define the weighted spaces by
\[
\| f \|_{Y_m}^2 = \sum_{j = 0 }^m \|  \langle y \rangle f^{(j)}      \|_2^2, \, \, \| f \|_Z^2 = \| f \|_{Y_4}^2 + \|   \langle y \rangle^2 f^{(4)} \|_2^2. 
\]
Oddness is imposed separately. For any of these spaces $Y$ and $\mathbf f = (f_1, f_2)$, write
\[
\| \mathbf f \|_Y^2 = \| f_1 \|_Y^2 + \| f_2 \|_Y^2,
\qquad \mathbf f \in Y \  \Longleftrightarrow \ (f_1, f_2) \in Y \times Y.
\]

\begin{Lemma}\label{3.15}
Let $a = \frac{1}{2}$. For real odd initial data $(\omega_0, J_0) \in Z \times Z$, there exist $t_0 > 0$ and a unique solution of \eqref{1.2} satisfying
\[
(\omega, J) \in C ( [0, t_0]; Z \times Z ), \, \, (\omega, J) ( \cdot , 0) = (\omega_0, J_0 ),
\]
and 
\[
(u, B) = (\mathcal P \omega, \mathcal P J), \, \, (\omega, J) (- x, t) = - (\omega, J) (x, t) .
\]
Higher regularity persists.
\[
(\omega_0, J_0) \in \bigcap_{m \geq 0}(Y_m \times Y_m)
\quad \Longrightarrow \quad
(\omega, J) \in \bigcap_{m \geq 0} C([0, t_0]; Y_m \times Y_m).
\]
Let $\tau_* \geq 0$ and let $(\omega_*, J_*) \in Z \times Z$ be real and odd, with $\omega_*'(0) = F'(0)$. There exist $\eta > 0$ and a unique solution of \eqref{1.10}, \eqref{1.12} satisfying
\[
(\omega, J) \in C ( [\tau_*, \tau_* + \eta ]; Z \times Z ) \, \, \text{and} \, \, (\omega, J) ( \cdot , \tau_* ) = (\omega_*, J_* ),
\]
where $\omega_y (0, \tau) = F'(0)$. The relative scales may be initialized by $A (\tau_*) = L (\tau_*) = 1$. 
\end{Lemma}

\begin{proof}
\emph{1. Weighted bounds and Picard construction.} Until Step~5, the fields are physical and $y = x$. For odd $f$,
\[
f^{(j)}(- y) = (- 1)^{j + 1} f^{(j)}(y).
\]
Integrating by parts,
\[
\int f^{(j)} = 0 \quad \text{when $j$ is odd}, \qquad
\int y f^{(4)} = 0.
\]
By Lemma~\ref{1.20} and interpolation, $\| \mathcal H f \|_{Y_m } \leq C_m \| f \|_{Y_m}$, $ \| \mathcal H f \|_Z \leq C \| f \|_Z $, and 
\[
\| \mathcal P f \|_\infty \leq C \| \langle y \rangle f \|_2,
\]
when $ j \geq 1 $, 
\[
(\mathcal P f )^{(j)} = \mathcal H f^{(j - 1 )}.
\]
Weighted interpolation also implies, for $0 \leq j \leq 4 $,
\[
\| \langle y \rangle^{1 + \frac{j}{4}} f^{(j)} \|_2 \leq C \| f \|_Z,
\]
 for $0 \leq j \leq 3 $,
\[
\| \langle y \rangle^{\frac{9}{8} + \frac{j}{4}} f^{(j)} \|_\infty \leq C \| f \|_Z. 
\]
For $\mathbf W = (\omega, J)^T$, write $\mathbf W_t + M(\mathbf W) \mathbf W_y = \mathcal R(\mathbf W)$, where
\[
M = \frac{1}{2} \begin{pmatrix}  u & - B \\ - B & u    \end{pmatrix}, \, \, \mathcal R = \begin{pmatrix}  \omega \mathcal H \omega - J \mathcal H J \\ \omega \mathcal H J - J \mathcal H \omega   \end{pmatrix}.
\]
The transport matrix has the constant orthogonal diagonalization
\[
\frac12 \begin{pmatrix} 1  & 1 \\ 1  & - 1 \end{pmatrix}
M \begin{pmatrix} 1  & 1 \\ 1  & - 1 \end{pmatrix}
= \frac12 \operatorname{diag}(u - B, u + B).
\]
Starting with $\mathbf W^0(t) = \mathbf W_0$, define the Picard iterates by
\[
\partial_t \mathbf W^{n + 1} + M (\mathbf W^n) \partial_y \mathbf W^{n + 1 } = \mathcal R (\mathbf W^n),
\qquad \mathbf W^{n + 1}(0) = \mathbf W_0.
\]
The velocities are bounded and Lipschitz, and the transport matrix is symmetric. Hence
\[
\frac{d}{d t} \| \mathbf W^{n + 1} \|_Z \leq C \| \mathbf W^n \|_Z \|  \mathbf W^{n + 1 }  \|_Z + C \| \mathbf W^n \|_Z^2. 
\]
If $\| \mathbf W_0 \|_Z \leq R$, choose $T_R > 0$ so that
\[
\sup_n \sup_{t \leq T_R } \|  \mathbf W^n (t) \|_Z \leq 2 R. 
\]
Let $\delta_n = \mathbf W^n - \mathbf W^{n - 1}$. Then
\[
\left ( \partial_t + M (\mathbf W ^n ) \partial_y \right) \delta_{n + 1 } = - \left (    M (\mathbf W^n) - M (\mathbf W^{n - 1}) \right) \partial_y \mathbf W^n + \mathcal R (\mathbf W^n ) - \mathcal R ( \mathbf W^{n - 1 }),
\]
and 
\[
\frac{d}{d t} \| \delta_{n + 1 } \|_{Y_3} \leq C_R \left (  \|   \delta_{n + 1 } \|_{Y_3} + \|  \delta_n  \|_{Y_3}  \right),
\]
where $\delta_{n + 1}(0) = 0$. Decrease $T_R$ so that $C_R T_R e^{C_R T_R} \leq \frac12$. Then the iterates converge in $C_t Y_3$, since
\[
\| \delta_{n + 1} \|_{C_t Y_3} \leq \frac{1}{2} \| \delta_n \|_{C_t Y_3}, \qquad
\sup_{t \leq T_R} \| \mathbf W(t) \|_Z \leq 2 R.
\]
\emph{2. Uniqueness in the low norm.} Let $\mathbf W_1, \mathbf W_2 \in L^\infty ([0, t_0]; Z) \cap C([0, t_0]; Y_3)$ be two real odd solutions. Their difference satisfies
\[
\begin{split}
& \partial_t(\mathbf W_1 - \mathbf W_2)
+ M(\mathbf W_1) \partial_y(\mathbf W_1 -  \mathbf W_2) \\
& \qquad = - \bigl (M(\mathbf W_1) - M(\mathbf W_2) \bigr) \partial_y \mathbf W_2
+ \mathcal R(\mathbf W_1) - \mathcal R(\mathbf W_2).
\end{split}
\]
The product terms satisfy
\[
\begin{split}
& \| \bigl (M(\mathbf W_1) - M(\mathbf W_2) \bigr) \partial_y \mathbf W_2 \|_{Y_3}
\leq C \| \mathbf W_2 \|_Z \| \mathbf W_1 - \mathbf W_2 \|_{Y_3}, \\
& \| \mathcal R(\mathbf W_1) - \mathcal R(\mathbf W_2) \|_{Y_3}
\leq C \bigl (\| \mathbf W_1 \|_Z + \| \mathbf W_2 \|_Z \bigr)
\| \mathbf W_1 - \mathbf W_2 \|_{Y_3}.
\end{split}
\]
Since $M(\mathbf W_1)^T = M(\mathbf W_1)$,
\[
\frac{d}{d t} \| \mathbf W_1 - \mathbf W_2 \|_{Y_3}^2
\leq C \bigl (\| \mathbf W_1 \|_Z + \| \mathbf W_2 \|_Z \bigr)
\| \mathbf W_1 - \mathbf W_2 \|_{Y_3}^2.
\]
Hence
\[
\begin{aligned}
\| \mathbf W_1(t) - \mathbf W_2(t) \|_{Y_3}^2
& \leq \| \mathbf W_1(0) - \mathbf W_2(0) \|_{Y_3}^2 \\
& \quad \times \exp \left (C \int_0^t
\bigl (\| \mathbf W_1(r) \|_Z + \| \mathbf W_2(r) \|_Z \bigr) \, d r \right).
\end{aligned}
\]
In particular, for equal initial data, $\mathbf W_1 = \mathbf W_2$.

\emph{3. Strong continuity by forward and backward restart.} For $t_n \to t_*$ in $[0, t_0]$, the convergence in $Y_3$ and the uniform $Z$ bound imply
\[
\mathbf W(t_n) \to \mathbf W(t_*) \quad \text{in } Y_3,
\qquad
\mathbf W(t_n) \rightharpoonup \mathbf W(t_*) \quad \text{in } Z,
\]
\[
\| \mathbf W(t_*) \|_Z \leq \liminf_{n \to \infty} \| \mathbf W(t_n) \|_Z.
\]
Fix $t_* \in [0, t_0]$ and assume $\sup_{0 \leq t \leq t_0} \| \mathbf W(t) \|_Z \leq 2 R$. For $h \geq 0$, consider
\[
\begin{split}
\partial_h \mathbf W_+ + M(\mathbf W_+ ) \partial_y \mathbf W_+
& = \mathcal R(\mathbf W_+ ), \\
\partial_h \mathbf W_- - M(\mathbf W_- ) \partial_y \mathbf W_-
& = - \mathcal R(\mathbf W_- ), \\
\mathbf W_+ (0) = \mathbf W_- (0)  & = \mathbf W(t_*).
\end{split}
\]
Start both Picard schemes with $\mathbf W_\pm^0(h) = \mathbf W(t_*)$. The norm estimates hold for either sign. There exist $\eta_R, C_R > 0$, independent of $t_*$ and $n$, such that
\[
\sup_n \sup_{0 \leq h \leq \eta_R} \| \mathbf W_\pm^n(h) \|_Z \leq 4 R,
\qquad
\frac{d}{d h} \| \mathbf W_\pm^{n + 1} \|_Z
\leq \frac{C_R}{2} \left (\| \mathbf W_\pm^{n + 1} \|_Z + \| \mathbf W_\pm^n \|_Z \right).
\]
After multiplying by $e^{- C_R h/2}$ and integrating,
\[
\| \mathbf W_\pm^{n + 1}(h) \|_Z
\leq {} e^{C_R h/2} \| \mathbf W(t_*) \|_Z
+ \frac{C_R}{2} \int_0^h e^{C_R(h - r)/2} \| \mathbf W_\pm^n(r) \|_Z \, d r.
\]
Since $\| \mathbf W_\pm^0(h) \|_Z = \| \mathbf W(t_*) \|_Z$, induction in the preceding integral inequality proves
\[
\| \mathbf W_\pm^{n + 1}(h) \|_Z
\leq \| \mathbf W(t_*) \|_Z e^{C_R h}.
\]
For every $h \in [0, \eta_R]$,
\[
\mathbf W_\pm^n(h) \to \mathbf W_\pm (h) \quad \text{in } Y_3,
\qquad
\mathbf W_\pm^n(h) \rightharpoonup \mathbf W_\pm (h) \quad \text{in } Z,
\]
\[
\| \mathbf W_\pm (h) \|_Z
\leq \liminf_{n \to \infty} \| \mathbf W_\pm^n(h) \|_Z
\leq \| \mathbf W(t_*) \|_Z e^{C_R h}.
\]
The functions $h \mapsto \mathbf W(t_* \pm h)$ solve the respective systems in $L^\infty Z \cap C_t Y_3$. By Step~2,
\[
\mathbf W_\pm (h) = \mathbf W(t_* \pm h)
\quad \text{if } 0 \leq h \leq \eta_R, \quad t_* \pm h \in [0, t_0].
\]
Thus, for $t \in [0, t_0]$ with $|t - t_*| \leq \eta_R$,
\[
\| \mathbf W(t) \|_Z \leq \| \mathbf W(t_*) \|_Z e^{C_R|t - t_*|},
\qquad
\limsup_{\substack{t \to t_* \\ t \in [0, t_0]}} \| \mathbf W(t) \|_Z
\leq \| \mathbf W(t_*) \|_Z.
\]
By weak lower semicontinuity and the preceding bound,
\[
\lim_{t \to t_*} \| \mathbf W(t) \|_Z = \| \mathbf W(t_*) \|_Z.
\]
Let $\langle \cdot , \cdot \rangle_Z$ denote the inner product of $Z \times Z$. Then
\[
\| \mathbf W(t) - \mathbf W(t_*) \|_Z^2
= \| \mathbf W(t) \|_Z^2 + \| \mathbf W(t_*) \|_Z^2
- 2 \operatorname{Re} \langle \mathbf W(t), \mathbf W(t_*) \rangle_Z
\longrightarrow 0.
\]
Thus $\mathbf W \in C([0, t_0]; Z)$.

\emph{4. Higher-order persistence.} For $m \geq 4$ and $0 \leq j \leq m$,
\[
\partial_t \partial_y^j \mathbf W + M(\mathbf W) \partial_y^{j + 1} \mathbf W
= - \sum_{k = 1}^j \binom j k \bigl (\partial_y^k M(\mathbf W) \bigr)
\partial_y^{j - k + 1} \mathbf W
+ \partial_y^j \mathcal R(\mathbf W),
\]
where the sum is zero for $j = 0$. For $f, g \in \{ \omega, J \}$, the weighted product estimates read
\[
\begin{split}
& \sum_{k = 1}^j \binom j k
\| \langle y \rangle (\mathcal H f^{(k - 1)}) g^{(j - k + 1)} \|_2
\leq C_m \left (\| f \|_{H^4} \| g \|_{Y_m}
+ \| g \|_{H^4} \| f \|_{Y_m} \right), 1 \leq j \leq m, \\
& \| \partial_y^j(f \mathcal H g) \|_{L^2(\langle y \rangle^2 d y)}
\leq C_m \left (\| f \|_{H^4} \| g \|_{Y_m}
+ \| g \|_{H^4} \| f \|_{Y_m} \right), 0 \leq j \leq m.
\end{split}
\]
Using $\partial_y^k \mathcal P f = \mathcal H f^{(k - 1)}$ for $k \geq 1$, we obtain
\[
\begin{split}
\frac12 \frac{d}{d t} \| \mathbf W \|_{Y_m}^2
= {}  & \frac12 \sum_{j = 0}^m \int
(\partial_y^j \mathbf W)^T \partial_y \bigl ((1 + y^2) M(\mathbf W) \bigr)
\partial_y^j \mathbf W \\
& - \sum_{j = 1}^m \sum_{k = 1}^j \binom j k \int (1 + y^2)
(\partial_y^j \mathbf W)^T(\partial_y^k M(\mathbf W))
\partial_y^{j - k + 1} \mathbf W \\
& + \sum_{j = 0}^m \int (1 + y^2)(\partial_y^j \mathbf W)^T
\partial_y^j \mathcal R(\mathbf W).
\end{split}
\]
Since
\[
\left \| \frac{\partial_y((1 + y^2) M(\mathbf W))}{1 + y^2} \right \|_\infty
\leq C \left (\| \mathbf W \|_{H^4} + \|(u, B) \|_\infty \right),
\]
we have
\[
\frac{d}{d t} \| \mathbf W \|_{Y_m}
\leq C_m \left (\| \mathbf W \|_{H^4} + \|(u, B) \|_\infty \right) \| \mathbf W \|_{Y_m},
\]
\[
\begin{split}
\| \mathbf W(t) \|_{Y_m}
& \leq \| \mathbf W(0) \|_{Y_m}
\exp \left (C_m \int_0^t \left (\| \mathbf W(r) \|_{H^4} + \|(u, B)(r) \|_\infty \right) \, d r \right) \\
& \leq \| \mathbf W(0) \|_{Y_m}
\exp \left (C_m \int_0^{t_0} \| \mathbf W(r) \|_Z \, d r \right) < \infty.
\end{split}
\]
These bounds pass from smooth approximations to the solution. Restarting at $t_*$,
\[
\| \mathbf W(t) \|_{Y_m}
\leq \| \mathbf W(t_*) \|_{Y_m}
\exp \left (C_m \left | \int_{t_*}^t
\left (\| \mathbf W(r) \|_{H^4} + \|(u, B)(r) \|_\infty \right) \, d r \right | \right).
\]
With weak continuity in $Y_m$, it follows that
\[
\lim_{t \to t_*} \| \mathbf W(t) - \mathbf W(t_*) \|_{Y_m} = 0,
\qquad
\mathbf W \in \bigcap_{m \geq 0} C([0, t_0]; Y_m)
\]
when $\mathbf W(0) \in \bigcap_{m \geq 0} Y_m$. For such data, repeated differentiation of \eqref{1.2} and these bounds at successively higher orders establish joint smoothness in space and time.

\emph{5. Normalized solutions and restart.} For the normalized initial data, take the physical solution with
$(\omega, J)( \cdot , 0) = (\omega_*, J_*)$ and set
\[
\begin{split}
\widetilde \omega (y, \tau (t))  & = A(t) \omega (L(t) y, t), \qquad
\widetilde J(y, \tau (t)) = A(t) J(L(t) y, t), \\
(\widetilde u, \widetilde B)(y, \tau (t))  & = \frac{A(t)}{L(t)}(u, B)(L(t) y, t).
\end{split}
\]
Then
\[
\begin{split}
\mathcal H(\widetilde \omega - F)(0, \tau (t))  & = A \mathcal H \omega (0, t) - \frac83, \\
\mathcal H \widetilde J(0, \tau (t)) \, \partial_y \widetilde J(0, \tau (t))
& = A^2 L \mathcal H J(0, t) J_x(0, t),
\end{split}
\]
and \eqref{1.12} becomes
\[
\begin{split}
c_A(t, A, L)  & = \frac53 - A \mathcal H \omega (0, t)
+ \frac{A^2 L}{F'(0)} \mathcal H J(0, t) J_x(0, t), \\
c_\ell (t, A, L)  & = \frac53 - \frac A 2 \mathcal H \omega (0, t)
+ \frac{A^2 L}{2 F'(0)} \mathcal H J(0, t) J_x(0, t).
\end{split}
\]
The traces satisfy
\[
| \mathcal H \omega (0, t)| + | \mathcal H J(0, t)| + |J_x(0, t)|
\leq C \| \mathbf W(t) \|_{H^2}.
\]
For $0 < r < R$ and $(A_i, L_i) \in [r, R]^2$, $i = 1, 2$,
\[
\begin{split}
& |c_A(t, A_1, L_1) - c_A(t, A_2, L_2)|
+ \left | \frac{L_1}{A_1} c_\ell (t, A_1, L_1)
- \frac{L_2}{A_2} c_\ell (t, A_2, L_2) \right |
+ \left | \frac1{A_1} - \frac1{A_2} \right | \\
& \leq C_{r, R} \left (1 + \| \mathbf W(t) \|_{H^2} + \| \mathbf W(t) \|_{H^2}^2 \right)
\left (|A_1 - A_2| + |L_1 - L_2| \right).
\end{split}
\]
Thus the system
\[
\begin{split}
A_t  & = c_A(t, A, L), \qquad
L_t = - \frac L A c_\ell (t, A, L), \qquad
\tau_t = \frac1 A, \\
A(0)  & = L(0) = 1, \qquad \tau (0) = \tau_*
\end{split}
\]
has a unique local solution with $A, L > 0$. In particular, $t_\tau = A$ and
\[
\begin{split}
\partial_\tau \widetilde \omega
& = c_A \widetilde \omega - c_\ell y \partial_y \widetilde \omega + A^2 \omega_t(L y, t), \\
\partial_\tau \widetilde J
& = c_A \widetilde J - c_\ell y \partial_y \widetilde J + A^2 J_t(L y, t).
\end{split}
\]
Substituting the physical equations \eqref{1.2} recovers the normalized equations \eqref{1.10} with the rates \eqref{1.12}.
For $k = 1$, $0 \leq j \leq 4$, and for $(k, j) = (2, 4)$,
\[
\begin{split}
\| \langle y \rangle^k \partial_y^j(A f(L y)) \|_2^2
& = A^2 L^{2 j - 1} \int \left (1 + \frac{x^2}{L^2} \right)^k|f^{(j)}(x)|^2 \, d x \\
& \leq A^2 L^{2 j - 1} \max \{ 1, L^{- 2 k} \}
\| \langle x \rangle^k f^{(j)} \|_2^2.
\end{split}
\]
By approximation with smooth compactly supported functions,
$(\widetilde \omega, \widetilde J) \in C([\tau_*, \tau_* + \eta]; Z \times Z)$ for some $\eta > 0$.
Solving \eqref{1.13},
\[
\partial_y \widetilde \omega (0, \tau) - F'(0)
= \bigl (\omega_*'(0) - F'(0) \bigr)
\exp \left (\int_{\tau_*}^{\tau}
\frac{\mathcal H \widetilde J(0, r) \, \partial_y \widetilde J(0, r)}{2 F'(0)} \, d r \right)
= 0.
\]
Conversely, solving
\[
A_\tau = c_A A, \qquad L_\tau = - c_\ell L, \qquad t_\tau = A,
\qquad (A, L, t)(\tau_*) = (1, 1, 0),
\]
and setting
\[
(\omega, J)(x, t(\tau)) = \frac1{A(\tau)}
(\widetilde \omega, \widetilde J) \left (\frac{x}{L(\tau)}, \tau \right)
\]
recovers \eqref{1.2}. The normalized solution is unique because the physical solution and the scale equations are unique.
\end{proof}

\begin{Proposition}\label{3.16}
For initial data $\omega ( \cdot , 0) = F$ and $J( \cdot , 0) = \varepsilon f$, where $f \in C_c^\infty (\R; \R)$ is fixed and odd, a smooth normalized solution exists for sufficiently small $\varepsilon$ and all $\tau \geq 0$. The bounds \eqref{3.12} and \eqref{3.17} hold globally. 

\end{Proposition} 

\begin{proof}
The initial pair $(F, \varepsilon f)$ belongs to $Z \cap \bigcap_m Y_m$. By Lemma~\ref{3.15}, a smooth local normalized solution exists. The preserved traces $s(0) = s'(0) = 0$ allow the second-order Hardy inequality.
\[
\| s/y^2 \|_{L^2(- 1, 1)} \leq C \| s'' \|_{L^2(- 1, 1)}.
\]
Together with
\[
\langle y \rangle^2 \leq \phi (y) \leq 4 \langle y \rangle^2
\quad (|y| \geq 1),
\]
we conclude $s \in C_\tau L^2(\phi)$. Choose $\varepsilon$ so that $E(0) < E_*$. It then follows from Proposition~\ref{3.11}
\[
E(\tau) \leq C \varepsilon e^{- 9 \tau /40}.
\]
In particular,
\[
\|(s, J) \|_{H^4} + \|(\mathcal P s, B) \|_\infty
+ |c_A + 1| + |c_\ell - \tfrac13| \leq C(E + E^2).
\]
The symmetric transport energy estimate takes the form
\begin{equation} \label{3.18}
\frac{d}{d \tau} \|(\omega, J) \|_{Y_4}^2 \leq C \|(\omega, J) \|_{Y_4}^2.
\end{equation}
Indeed, the preceding bound, $\| y \partial_y \log \langle y \rangle^2 \|_\infty \leq 2$, and
\[
\| f g \|_{Y_4} \leq C \left (\| f \|_{H^4} \| g \|_{Y_4}
+ \| g \|_{H^4} \| f \|_{Y_4} \right)
\]
control all coefficients and lower-order products. The constant in \eqref{3.18} is uniform as long as $E \leq E_*$. Hence
\[
\begin{split}
\|(\omega, J)(\tau) \|_{Y_4}
& \leq \|(\omega, J)(0) \|_{Y_4} e^{C \tau /2}, \\
\|(\omega, J)(\tau) \|_Z
& \leq C \left (\|(\omega, J)(\tau) \|_{Y_4} + E(\tau)
+ \| \langle y \rangle^2 F^{(4)} \|_2 \right).
\end{split}
\]
If the maximal normalized time $\tau_{\max}$ were finite, the preceding bounds would imply
\[
\sup_{0 \leq \tau < \tau_{\max}} \|(\omega, J)(\tau) \|_Z < \infty.
\]
The uniform local lifespan in Lemma~\ref{3.15} would extend the solution past $\tau_{\max}$. Hence $\tau_{\max} = \infty$. Its higher-order estimates also ensure, for every finite $\tau_0$,
\[
(\omega, J) \in \bigcap_{m \geq 0} C([0, \tau_0]; Y_m \times Y_m).
\]

For every finite $\tau_0$, Step~1 of Proposition~\ref{3.21} applies to
\[
(\omega, J) \in C([0, \tau_0]; Z \times Z), \qquad
\omega ( \cdot , 0) = F, \qquad s( \cdot , 0) = 0.
\]
Thus $M_3 \in C([0, \tau_0])$ before damping is used. That proposition therefore proves \eqref{3.17} for all $\tau \geq 0$.
\end{proof}

\subsection{Proof of finite-time blowup}

\begin{Lemma}\label{3.19}
Let $\ell$ be the left eigenfunctional in Proposition~\ref{2.17}. Choose even $\chi \in C_c^\infty$ equal to one near zero, and define
\[
\ell (\mathcal N) = \lim_{R \to \infty} \int \chi (y/R) \zeta \mathcal N,
\qquad \mathcal N \text{ as in \eqref{1.19}}.
\]
For $\mathcal N \in X$, this equals the ordinary dual pairing. Then
\begin{equation} \label{3.20}
| \ell ( \mathcal N ) | \leq C (  \| s \|_{L^2 (\phi)} + \| s \|_{H^1} + \| J \|_{H^2}^2    ) \| J \|_X. 
\end{equation}
\end{Lemma}

\begin{proof}
By Proposition~\ref{2.17}, $\zeta \in H^1(\R)$ and $y \zeta' \in L^2(\R)$. Integration by parts moves the transport derivatives onto this fixed adjoint function.
\[
- \delta_\ell \int y \zeta J'
= \delta_\ell \int (\zeta + y \zeta') J, \qquad - \frac12 \int (\mathcal P s) \zeta J'
= \frac12 \int \left ((\mathcal P s) \zeta' + (\mathcal H s) \zeta \right) J.
\]
These identities also prove existence of the cutoff limits. Indeed, for $\chi = 1$ on $[- 1, 1]$ and $\operatorname{supp} \chi \subset [- 2, 2]$, the extra term is bounded by
\[
\left | \frac1 R \int \chi' \left (\frac y R \right)
\left (\delta_\ell y + \frac12 \mathcal P s \right) \zeta J \right |
\leq C \left (| \delta_\ell | + \frac{\| \mathcal P s \|_\infty}{R} \right)
\| \zeta \|_{L^2(|y| \geq R)} \| J \|_2 \longrightarrow 0.
\]
By \eqref{3.3}, Lemma~\ref{1.21}, Sobolev embedding and the $H^1$ boundedness of $\mathcal H$,
\[
\begin{gathered}
| \mathcal H s(0)| + \| \mathcal P s \|_\infty \leq C \| s \|_{L^2(\phi)}, \qquad
\| s \|_\infty + \| \mathcal H s \|_\infty \leq C \| s \|_{H^1}, \\
|q_J| \leq C \| J \|_{H^2}^2, \qquad
| \delta_A| + | \delta_\ell |
\leq C \left (\| s \|_{L^2(\phi)} + \| J \|_{H^2}^2 \right).
\end{gathered}
\]
Cauchy--Schwarz bounds the integrated transport terms.
\[
\begin{split}
\left | - \delta_\ell \int y \zeta J' \right |
+ \left | \delta_A \int \zeta J \right |
& \leq C \left (\| s \|_{L^2(\phi)} + \| J \|_{H^2}^2 \right) \| J \|_X, \\
\left | - \frac12 \int (\mathcal P s) \zeta J' \right |
& \leq C \left (\| s \|_{L^2(\phi)} + \| s \|_{H^1} \right) \| J \|_X.
\end{split}
\]
The remaining products satisfy
\[
\begin{split}
\left | \int \zeta (s \mathcal H J - J \mathcal H s) \right |
& \leq \| \zeta \|_2
\left (\| s \|_\infty \| \mathcal H J \|_2 + \| \mathcal H s \|_\infty \| J \|_2 \right)
\leq C \| s \|_{H^1} \| J \|_X, \\
\left | \int \zeta B s' \right |
& \leq \| B \|_\infty \| \zeta \|_2 \| s' \|_2
\leq C \| J \|_X \| s \|_{H^1}.
\end{split}
\]
Adding the bounds proves \eqref{3.20}. The constant depends only on the fixed weights and $\| \zeta \|_{H^1} + \| y \zeta' \|_2$.
\end{proof}

\begin{proof}[Proof of Theorem~\ref{1.4}]
\emph{1. Choice of magnetic seed.}
The functional $\ell$ is nonzero and real on real functions. Since real odd functions in $C_c^\infty (\R)$ are dense in the real subspace of $X$, continuity implies that some such $h$ satisfies $\ell (h) \neq 0$. Set
\[
f = \frac{h}{\ell (h)} \in C_c^\infty (\R; \R), \qquad f(- y) = - f(y), \qquad \ell (f) = 1.
\]
For sufficiently small $\varepsilon > 0$, Proposition~\ref{3.16} provides a global normalized solution with $\omega ( \cdot , 0) = F$ and $J( \cdot , 0) = \varepsilon f$.

\emph{2. Magnetic lower bound.}
Test \eqref{1.18} against smooth approximations of $\chi (y/R) \zeta$. The cutoff argument in Lemma~\ref{3.19} treats the nonlinear transport term, and the linear transport error tends to zero because $V(y) = O(|y|)$ and
\[
\left | \int \frac y R \chi' \left (\frac y R \right) \zeta J \, d y \right |
\leq C \| \zeta \|_{L^2(|y| > R)} \| J \|_2 \longrightarrow 0.
\]
The adjoint relation from Proposition~\ref{2.17} now reads
\[
\frac d{d \tau} \ell (J) = \lambda_* \ell (J) + \ell (\mathcal N), \qquad
\ell (J( \cdot , 0)) = \varepsilon,
\]
with the cutoff pairing defined in Lemma~\ref{3.19}. Propositions~\ref{3.11} and~\ref{3.21}, together with that lemma, imply
\[
| \ell (\mathcal N(\tau))| \leq C \varepsilon^3 e^{- 27 \tau /40}, \qquad
\frac{27}{40} + \lambda_* > \frac{13}{360} > 0.
\]
By variation of constants, uniformly for $\tau \geq 0$,
\[
\left |e^{- \lambda_* \tau} \ell (J(\tau)) - \varepsilon \right |
\leq C \varepsilon^3 \int_0^\tau e^{- (27/40 + \lambda_*) r} \, d r
\leq \frac{C \varepsilon^3}{27/40 + \lambda_*}.
\]
Decrease $\varepsilon_0$ so that the last quantity is at most $\varepsilon /2$. Then $| \ell (J(\tau))| \geq \varepsilon e^{\lambda_* \tau}/2$. Using \eqref{2.18},
\begin{equation} \label{3.22}
\| B(\tau) \|_\infty \geq c \varepsilon e^{\lambda_* \tau}.
\end{equation}

\emph{3. Return to physical variables.}
The improved scalar bound and \eqref{3.3} imply
\[
|c_A + 1| + |c_\ell - \tfrac13| \leq C \varepsilon^2 e^{- 9 \tau /20}.
\]
These errors are integrable with integral at most $C \varepsilon^2$. Integrating \eqref{1.9}, we obtain positive limiting factors $A_\infty, L_\infty$ satisfying
\[
e^{- C \varepsilon_0^2} \leq A_\infty, L_\infty \leq e^{C \varepsilon_0^2}.
\]
In particular, $T = \int_0^\infty A(r) \, d r < \infty$, and
\begin{equation} \label{3.23}
A(\tau) \sim A_\infty e^{- \tau}, \qquad L(\tau) \sim L_\infty e^{- \tau /3}, \qquad T - t(\tau) = \int_\tau^\infty A(r) \, d r \sim A_\infty e^{- \tau}.
\end{equation}
Since $\| \mathcal P s \|_\infty \to 0$, the normalized velocity satisfies $\| U + \mathcal P s \|_\infty \to \| U \|_\infty > 0$. The physical potentials are
\[
(u, B)_{\mathrm{phys}}(L y, t) = \frac L A(u, B)_{\mathrm{res}}(y, \tau).
\]
Combining this identity with \eqref{3.22} and \eqref{3.23} proves \eqref{1.5} with $\gamma = \frac23 + \lambda_*$. The bounds on $A_\infty, L_\infty$ make the constants uniform in $0 < \varepsilon < \varepsilon_0$.

\emph{4. Finite pre-blowup energy.}
At each finite normalized time, for $f \in \{ \omega, J \}$,
\[
f \in Y_4, \qquad f, y f \in L^2, \qquad
\| f \|_1 \leq \| \langle y \rangle^{- 1} \|_2 \| \langle y \rangle f \|_2 < \infty.
\]
With $\widehat f(\eta) = \int e^{- i y \eta} f(y) \, d y$, oddness and Plancherel give
\[
\widehat f(0) = 0, \qquad
\partial_\eta \widehat f = - i \widehat{y f} \in L^2,
\qquad \widehat f \in H^1.
\]
Hardy's inequality bounds the potential.
\[
\| \mathcal P f \|_2^2
= \frac1{2 \pi} \int \frac{| \widehat f(\eta)|^2}{| \eta |^2} \, d \eta
\leq \frac4{2 \pi} \int | \partial_\eta \widehat f(\eta)|^2 \, d \eta
= 4 \| y f \|_2^2.
\]
Thus, for every $t < T$,
\[
\| u_{\mathrm{phys}}(t) \|_2^2 + \| B_{\mathrm{phys}}(t) \|_2^2
= \frac{L^3}{A^2} \left (\| \mathcal P \omega (\tau) \|_2^2 + \| \mathcal P J(\tau) \|_2^2 \right) < \infty.
\]
\end{proof}

\bigskip

\section{Parameters near $\frac{1}{2}$ }
We prove Theorem~\ref{1.6} using the profile in \eqref{1.3}.
\[
U = \mathcal P F, \qquad F'(0) = - \frac{16}{3}, \qquad \mathcal H F(0) = \frac83.
\]
\subsection{Nearby profiles and spectral persistence}
We use the inviscid profile construction in \cite[Sections 2.5 and 2.6]{Chen2020}, with the change of variables $x = b y$, $ b^2 = \frac{3}{8}$. 

\begin{Lemma}\label{4.1}
There exists $\delta_0 > 0$ such that, for every $ a \in \R$ satisfying $\left | a - \frac{1}{2} \right | < \delta_0$, there is a profile $F_a \in C^\infty (\R; \R )$ with
\[
F_a (- y) = - F_a (y), \, \, F'_a (0) = F'(0), \, \, U_a = \mathcal P F_a. 
\]
The profile also satisfies
\begin{equation} \label{4.2}
V_a F_a ' = ( \bar c_A + Q_a ) F_a, \, \, V_a = \bar c_\ell y + a U_a. 
\end{equation}
For every integer $m \geq 0 $, there is a constant $C_m > 0$, independent of $a$ on $ \left |   a - \frac{1}{2}  \right | < \delta_0$, such that
\begin{equation} \label{4.3}
\| F_a - F \|_{L^2 (\phi)} + \sum_{j = 0 }^m \|     \langle y \rangle^{j + 1} \left ( F_a^{(j)} - F^{(j)} \right)         \|_2 \leq C_m \left |   a - \frac{1}{2} \right |. 
\end{equation}
There are constants $c > 0$, $C > 0$, independent of $a$, such that  
\[
| \bar c_A + 1 | + \left | \bar c_\ell - \frac{1}{3} \right | \leq C \left | a - \frac{1}{2}    \right |,
\]
Moreover, for $y \neq 0$ we have $ 0 < c \leq \frac{V_a(y)}{y} \leq C$. At $y = 0 $, we define $\frac{V_a(y)}{y}  = V'_a (0) $ by continuity. 
\end{Lemma}

\begin{proof}
In the normalized scalar equation, set
\[
c_A = \frac53 - \mathcal H \omega (0), \qquad
c_\ell = \frac53 - a \mathcal H \omega (0), \qquad \omega ( \cdot , 0) = F.
\]
The residual is
\[
R_a = - \left ( a - \frac{1}{2} \right) (U - \frac{8}{3} y ) F',
\qquad R_a'(0) = 0.
\]
Moreover,
\[
\left |a - \frac12 \right | \Biggl (\left \| \left (U - \frac83 y \right) F' \right \|_{L^2(\psi)} + \left \| \langle y \rangle^3 \partial_y^4 \left [\left (U - \frac83 y \right) F' \right] \right \|_2 \Biggr) \leq C \left |a - \frac12 \right |.
\]
Let $d = |a - \tfrac12|$, and let $\omega_a(y, \tau)$ be the normalized scalar trajectory with parameter $a$ and initial profile $\omega_a( \cdot , 0) = F$. Chen's scalar stability and convergence estimates give
\[
\sup_{\tau \geq 0 } \left (  \| \omega_a (\tau) - F \|_{L^2 (\phi )}   +  \|    \langle y \rangle^3 \partial_y^4 ( \omega_a (\tau) - F )     \|_2   \right) \leq C d, 
\]
and
\[
\lim_{\tau \to \infty } \| \omega_a (\tau) - F_a \|_{L^2 (\phi )} = 0. 
\]
The fourth-order bound, convergence, and time-derivative decay below are from \cite[Sections 2.5-2.6]{Chen2020}. Set $g = F_a - F$. The subscript $*$ denotes $a = 1/2$.

\emph{1. Weak compactness and the stationary limit.}
For any sequence $\tau_n \to \infty$, set $g_n = \omega_a(\tau_n) - F$. After passing to a subsequence,
\[
\langle y \rangle^3 g_n^{(4)} \rightharpoonup z \quad \hbox{in } L^2,
\qquad g_n \to g \quad \hbox{in } L^2(\phi).
\]
For every $\varphi \in C_c^\infty$, integration by parts shows
\[
\langle z, \varphi \rangle
= \lim_n \langle g_n, \partial_y^4(\langle y \rangle^3 \varphi) \rangle
= \langle g, \partial_y^4(\langle y \rangle^3 \varphi) \rangle.
\]
Thus $z = \langle y \rangle^3 g^{(4)}$ distributionally. By weak lower semicontinuity and unweighted Fourier interpolation,
\[
\| g \|_{L^2(\phi)} + \| \langle y \rangle^3 g^{(4)} \|_2
+ \| g \|_{H^4} \leq C d.
\]
For $0 \leq r < 4$, interpolation proves strong convergence.
\[
\| g_n - g \|_{H^r}
\leq C \| g_n - g \|_2^{1 - r/4} \| g_n - g \|_{H^4}^{r/4}
\longrightarrow 0.
\]
Thus
\[
g(- y) = - g(y), \qquad g'(0) = \lim_n g_n'(0) = 0,
\qquad \| \mathcal H(g_n - g) \|_{H^3} \longrightarrow 0.
\]
We use Lemma~\ref{1.21} for the potentials and Chen's estimate for the time derivative.
\[
\| \mathcal P(g_n - g) \|_\infty
\leq C \| g_n - g \|_X
\leq C \| g_n - g \|_{L^2(\phi)} \longrightarrow 0, \qquad \| \partial_\tau \omega_a(\tau_n) \|_{L^2(\phi)} \longrightarrow 0.
\]
Passing to the scalar equation on each compact interval proves \eqref{4.2}, with
\[
Q_a = \mathcal H F_a, \qquad
\bar c_A = \frac53 - Q_a(0), \qquad
\bar c_\ell = \frac53 - a Q_a(0).
\]
The $H^4$ bound controls the modulation and transport coefficients.
\[
| \bar c_A + 1| + | \bar c_\ell - \tfrac13|
+ \left \| \frac{V_a - V}{y} \right \|_\infty \leq C d.
\]
Indeed,
\[
\left \| \frac{\mathcal P g}{y} \right \|_\infty
\leq \| \mathcal H g \|_\infty \leq C \| g \|_{H^4} \leq C d.
\]
Reduce $\delta_0$ so that $C \delta_0 \leq 1/6$. Then
\[
\frac16 \leq \frac{V_a(y)} y \leq C, \qquad
V_a'(0) = \bar c_A + Q_a(0) = \frac53.
\]

\emph{2. The source trace and the origin bootstrap.}
Suppose inductively that $\| g \|_{H^j} \leq C_j d$ for some $j \geq 4$. Differentiate the stationary equation in distributions.
\[
V_a F_a^{(j + 1)} + A_{a, j} F_a^{(j)} = F_a \mathcal H F_a^{(j)} + G_{a, j},
\qquad A_{a, j} = j V_a' - \bar c_A - Q_a,
\]
with the lower-order source
\[
G_{a, j} = \sum_{k = 1}^{j - 1}
\left [\binom j k - a \binom j{k + 1} \right]
Q_a^{(k)} F_a^{(j - k)}.
\]
For $1 \leq k \leq j - 1$,
\[
Q_a^{(k)}, F_a^{(j - k)} \in H^1(\R), \qquad
\| v w \|_{H^1} \leq C \| v \|_{H^1} \| w \|_{H^1}.
\]
Under the induction hypothesis,
\[
\| G_{a, j} \|_{H^1} \leq C_j, \qquad
\| G_{a, j} - G_{*, j} \|_{H^1} \leq C_j d.
\]
The source need not vanish at zero. For example, $G_{*, 5}(0) = - 12800$. Define, only on $(- 1, 1)$,
\[
c_{a, j} = \frac{G_{a, j}(0)}{(j - 1) \frac53}, \qquad
\widetilde f_{a, j} = F_a^{(j)} - c_{a, j}, \qquad
T_{a, j} = G_{a, j} - c_{a, j} A_{a, j}.
\]
These definitions use only the traces of $G_{a, j}$ and $A_{a, j}$.
\[
A_{a, j}(0) = (j - 1) \frac53, \qquad T_{a, j}(0) = 0.
\]
Moreover,
\[
|c_{a, j} - c_{*, j}| +
\| T_{a, j} - T_{*, j} \|_{H^1(- 1, 1)} \leq C_j d.
\]
Set
\[
\beta_{a, j} = \frac{A_{a, j}}{V_a} - \frac{j - 1}{y}.
\]
The coefficient identities are
\[
\begin{split}
V_a  & = \frac53 y + a \int_0^y \bigl (Q_a(t) - Q_a(0) \bigr) \, d t, \\
A_{a, j}  & = (j - 1) \frac53 + (j a - 1) \bigl (Q_a - Q_a(0) \bigr).
\end{split}
\]
Thus
\[
\beta_{a, j}
= \frac1{V_a} \left [(j a - 1) \bigl (Q_a - Q_a(0) \bigr)
- \frac{(j - 1) a}{y} \int_0^y \bigl (Q_a(t) - Q_a(0) \bigr) \, d t \right].
\]
Evenness and the $H^4$ bounds give the Taylor estimates on $(- 1, 1)$.
\[
Q_a(y) - Q_a(0) = O(y^2), \qquad
Q_a(y) - Q_a(0) - Q(y) + Q(0) = O(d y^2).
\]
Consequently,
\[
| \beta_{a, j}(y)| \leq C_j|y|, \qquad
| \beta_{a, j}(y) - \beta_{*, j}(y)| \leq C_j d|y|,
\qquad
\left \| \frac{F_a}{V_a} - \frac F V \right \|_{L^\infty (- 1, 1)} \leq C d.
\]
On $(- 1, 0) \cup (0, 1)$ the corrected equation is
\[
\widetilde f_{a, j}' + \frac{j - 1}{y} \widetilde f_{a, j} = r_{a, j},
\qquad
r_{a, j} = \frac{F_a}{V_a} \mathcal H F_a^{(j)}
+ \frac{T_{a, j}}{V_a} - \beta_{a, j} \widetilde f_{a, j}.
\]
Hardy's inequality implies
\[
\left \| \frac{T_{a, j}}{V_a} \right \|_{L^2(- 1, 1)}
\leq 6 \left \| \frac{T_{a, j}} y \right \|_{L^2(- 1, 1)}
\leq 12 \| T_{a, j}' \|_{L^2(- 1, 1)}, \qquad
r_{a, j} \in L^2(- 1, 1).
\]
Since $y^{- (j - 1)} \notin L^2(0, 1)$, the distributional equation has the representation
\[
\widetilde f_{a, j}(y)
= y^{- (j - 1)} \int_0^y t^{j - 1} r_{a, j}(t) \, d t.
\]
In particular,
\[
\frac{\widetilde f_{a, j}(y)} y
= \int_0^1 s^{j - 1} r_{a, j}(s y) \, d s, \qquad
\| r_{a, j}(s \, \cdot ) \|_{L^2(0, 1)}
\leq s^{- 1/2} \| r_{a, j} \|_{L^2(0, 1)}.
\]
By Minkowski's inequality,
\[
\left \| \frac{\widetilde f_{a, j}} y \right \|_{L^2(0, 1)}
\leq \frac{\| r_{a, j} \|_{L^2(0, 1)}}{j - \frac12}, \qquad
\| \widetilde f_{a, j}' \|_{L^2(0, 1)}
\leq \left (1 + \frac{j - 1}{j - \frac12} \right) \| r_{a, j} \|_{L^2(0, 1)}.
\]
On $(- 1, 0)$, replace $y$ by $- y$. For $0 < |y| < 1$,
\[
| \widetilde f_{a, j}(y)|
\leq \frac{|y|^{1/2}}{\sqrt{2 j - 1}}
\| r_{a, j} \|_{L^2(- |y|, |y|)} \longrightarrow 0.
\]
Hence $F_a^{(j)} \in H^1(- 1, 1)$ and $F_a^{(j)}(0) = c_{a, j}$.

The difference satisfies
\[
\frac{T_{a, j}}{V_a} - \frac{T_{*, j}} V
= \frac{T_{a, j} - T_{*, j}}{V_a}
+ \left (\frac1{V_a} - \frac1 V \right) T_{*, j}, \qquad
\left | \frac1{V_a} - \frac1 V \right | \leq \frac{C d}{|y|}.
\]
Both $T_{a, j} - T_{*, j}$ and $T_{*, j}$ vanish at zero. Hardy's inequality bounds the preceding difference.
\[
\left \| \frac{T_{a, j}}{V_a} - \frac{T_{*, j}} V \right \|_{L^2(- 1, 1)} \leq C_j d.
\]
Together with $\| \mathcal H(F_a^{(j)} - F^{(j)}) \|_2 \leq C_j d$, this controls the source difference.
\[
\| r_{a, j} - r_{*, j} \|_{L^2(- 1, 1)} \leq C_j d.
\]
The same Hardy bounds apply to the difference of the two representations.
\[
\| g^{(j + 1)} \|_{L^2(- 1, 1)} \leq C_j d.
\]

\emph{3. Exterior regularity and unweighted differences.}
The exact difference equation is
\[
V_a g' = (\bar c_A + Q_a) g + F \mathcal H g + (\bar c_A + 1) F - (V_a - V) F'.
\]
Its $j$th derivative reads
\[
\begin{aligned}
V_a g^{(j + 1)}
& = \partial_y^j \left [(\bar c_A + Q_a) g + F \mathcal H g + (\bar c_A + 1) F - (V_a - V) F' \right] \\
& \quad - \sum_{r = 1}^j \binom j r V_a^{(r)} g^{(j + 1 - r)}.
\end{aligned}
\]
For $|y| \geq 1$, $1/|V_a| \leq C$. By Sobolev embedding and the decay of $F$,
\[
\begin{gathered}
V_a - V = O(d|y|), \qquad
\| V_a^{(r)} - V^{(r)} \|_\infty \leq C_j d
\quad (1 \leq r \leq j), \\
\left \|(\bar c_A + 1) F - (V_a - V) F' \right \|_{H^j(|y| \geq 1)}
\leq C_j d.
\end{gathered}
\]
At the highest coefficient orders,
\[
\| Q_a^{(j)} g \|_2
\leq \| Q_a^{(j)} \|_2 \| g \|_\infty \leq C_j d, \qquad \| V_a^{(j)} g' \|_2
\leq |a| \| Q_a^{(j - 1)} \|_\infty \| g' \|_2 \leq C_j d.
\]
For $1 \leq r \leq j - 1$, the other products satisfy
\[
\begin{aligned}
\| Q_a^{(r)} g^{(j - r)} \|_2
& \leq \| Q_a^{(r)} \|_\infty \| g^{(j - r)} \|_2 \leq C_j d, \\
\| V_a^{(r)} g^{(j + 1 - r)} \|_2
& \leq \| V_a^{(r)} \|_\infty \| g^{(j + 1 - r)} \|_2 \leq C_j d.
\end{aligned}
\]
Also,
\[
\| \partial_y^j(F \mathcal H g) \|_2
\leq C_j \sum_{r = 0}^j \| F^{(r)} \|_\infty
\| \mathcal H g^{(j - r)} \|_2 \leq C_j d.
\]
Thus $\| g^{(j + 1)} \|_{L^2(|y| \geq 1)} \leq C_j d$. Together with the origin estimate,
\[
\| g \|_{H^j} \leq C_j d \quad \Longrightarrow \quad
\| g \|_{H^{j + 1}} \leq C_{j + 1} d, \qquad j \geq 4.
\]
Step~1 starts the induction, so $F_a \in C^\infty (\R)$.

\emph{4. Moments, Hilbert transforms and all-order weighted differences.}
Define
\[
W_j(g) = \sum_{r = 0}^j \| \langle y \rangle^{r + 1} g^{(r)} \|_2.
\]
By Step~1, $W_0(g) \leq C d$. Assume $W_j(g) \leq C_j d$. For $0 \leq k \leq j$,
\[
\int |y^k g^{(j)}(y)| \, d y
\leq \| \langle y \rangle^{j + 1} g^{(j)} \|_2
\| y^k \langle y \rangle^{- j - 1} \|_2 < \infty.
\]
Choose an even $\chi \in C_c^\infty$ with
\[
\chi = 1 \  \hbox{on }[- 1, 1], \qquad
\operatorname{supp} \chi \subset [- 2, 2], \qquad \chi_R(y) = \chi (y/R).
\]
It follows from integration by parts
\[
\int \chi_R y^k g^{(j)} = (- 1)^j \int g \, \partial_y^j(y^k \chi_R),
\]
The terms without cutoff derivatives vanish.
\[
\partial_y^j y^k = 0 \quad (k < j), \qquad
\int \chi_R g = 0 \quad (k = j).
\]
The terms with cutoff derivatives are bounded by
\[
C_j R^{k - j} \| g \|_{L^1(R \leq |y| \leq 2 R)}
\leq C_j R^{k - j - \frac12}
\| \langle y \rangle g \|_{L^2(R \leq |y| \leq 2 R)} \longrightarrow 0.
\]
Consequently,
\[
\int y^k g^{(j)}(y) \, d y = 0, \qquad 0 \leq k \leq j,
\]
including $\int g = 0$ for $j = 0$. The moment cancellations eliminate the polynomial remainder in the Hilbert kernel.
\[
y^n \mathcal H(g^{(j)}) = \mathcal H(y^n g^{(j)}),
\qquad 0 \leq n \leq j + 1.
\]
The cutoff identities pass to distributions by the absolute moment bounds, and $y^n g^{(j)} \in L^2$. The Hilbert isometry then implies
\[
\| \langle y \rangle^{j + 1} \mathcal H g^{(j)} \|_2^2
= \sum_{n = 0}^{j + 1} \binom{j + 1} n \| \mathcal H(y^n g^{(j)}) \|_2^2
= \| \langle y \rangle^{j + 1} g^{(j)} \|_2^2.
\]
For $0 \leq k \leq j - 1$,
\[
\| \langle y \rangle^{k + 1} g^{(k)} \|_{H^1}
\leq C_j \left (
\| \langle y \rangle^{k + 1} g^{(k)} \|_2 +
\| \langle y \rangle^{k + 2} g^{(k + 1)} \|_2 \right)
\leq C_j W_j(g).
\]
The supremum bound follows by Sobolev embedding.
\[
\| \langle y \rangle^{k + 1} g^{(k)} \|_\infty \leq C_j W_j(g),
\qquad 0 \leq k \leq j - 1.
\]
The corresponding coefficient estimates, using $V_a^{(r)} = a Q_a^{(r - 1)}$ for $r \geq 2$, are
\[
\begin{split}
& \| \langle y \rangle^{r + 1} Q_a^{(r)} \|_2 \leq C_j, \qquad
\| \langle y \rangle^{r + 1}(Q_a^{(r)} - Q^{(r)}) \|_2 \leq C_j d,
\quad 0 \leq r \leq j, \\
& \| V_a' \|_\infty \leq C, \qquad
\|(V_a' - V') \|_\infty + \|(V_a - V)/y \|_\infty \leq C d, \\
& \| \langle y \rangle^r V_a^{(r)} \|_2 \leq C_j, \qquad
\| \langle y \rangle^r(V_a^{(r)} - V^{(r)}) \|_2 \leq C_j d,
\quad 2 \leq r \leq j.
\end{split}
\]
Here $C_j$ depends only on the preceding weighted and unweighted bounds. Since
\[
|F^{(r)}(y)| \leq C_r \langle y \rangle^{- r - 3},
\]
Therefore we get
\[
\left \| \langle y \rangle^{j + 1} \partial_y^j \left [(\bar c_A + 1) F - (V_a - V) F' \right] \right \|_2 \leq C_j d.
\]
For the differentiated difference equation in Step~3, the remaining products satisfy
\[
\begin{split}
& \| \langle y \rangle^{j + 1} \partial_y^j((\bar c_A + Q_a) g) \|_2
+ \| \langle y \rangle^{j + 1} \partial_y^j(F \mathcal H g) \|_2 \\
& \quad + \sum_{r = 1}^j \| \langle y \rangle^{j + 1} V_a^{(r)} g^{(j + 1 - r)} \|_2
\leq C_j W_j(g).
\end{split}
\]
Indeed, the endpoint terms satisfy
\[
\| \langle y \rangle^{j + 1}(\bar c_A + Q_a) g^{(j)} \|_2 + \| \langle y \rangle^{j + 1} F \mathcal H g^{(j)} \|_2 + \| \langle y \rangle^{j + 1} V_a' g^{(j)} \|_2 \leq C_j W_j(g).
\]
For $1 \leq r \leq j$ and $2 \leq r \leq j$, respectively,
\[
\begin{aligned}
& \| \langle y \rangle^{j + 1} Q_a^{(r)} g^{(j - r)} \|_2 \leq \| \langle y \rangle^{r + 1} Q_a^{(r)} \|_2
\| \langle y \rangle^{j - r} g^{(j - r)} \|_\infty \leq C_j W_j(g), \\
& \| \langle y \rangle^{j + 1} V_a^{(r)} g^{(j + 1 - r)} \|_2 \leq \| \langle y \rangle^r V_a^{(r)} \|_2
\| \langle y \rangle^{j + 1 - r} g^{(j + 1 - r)} \|_\infty \leq C_j W_j(g).
\end{aligned}
\]
For $0 \leq r \leq j$,
\[
\| \langle y \rangle^{j + 1} F^{(r)} \mathcal H g^{(j - r)} \|_2 \leq \| \langle y \rangle^r F^{(r)} \|_\infty \| \langle y \rangle^{j - r + 1} \mathcal H g^{(j - r)} \|_2 \leq C_j W_j(g).
\]
For $j = 0$, only the endpoint terms occur.

Finally, $\langle y \rangle /|V_a(y)| \leq C$ for $|y| \geq 1$. Division by the transport coefficient bounds the next weighted derivative.
\[
\| \langle y \rangle^{j + 2} g^{(j + 1)} \|_{L^2(|y| \geq 1)} \leq C_j d.
\]
On $[- 1, 1]$, the weight is bounded, so Step~3 applies.
\[
\| \langle y \rangle^{j + 2} g^{(j + 1)} \|_{L^2(- 1, 1)}
\leq 2^{(j + 2)/2} \| g^{(j + 1)} \|_2 \leq C_j d.
\]
Thus, for every $j \geq 0$,
\[
W_j(g) \leq C_j d \quad \Longrightarrow \quad
W_{j + 1}(g) \leq C_{j + 1} d.
\]
This proves \eqref{4.3}. The interval $|a - 1/2| < \delta_0$ remains the one fixed in Step~1.
\end{proof}

Under the rescaling \eqref{1.8}, \eqref{1.9}, the normalized equations read 
\[
\omega_\tau + (c_\ell y + a u) \omega' = c_A \omega + \omega \mathcal H \omega - J \mathcal H J + a B J' , 
\]
\[
J_\tau + (c_\ell y + a u ) J' = c_A J + \omega \mathcal H J - J \mathcal H \omega + a B \omega' . 
\]
Write $\omega = F_a + s$, $u = U_a + \mathcal P s$, and choose
\[
q_J = \frac{\mathcal H J(0) J'(0) }{F'(0)} \qquad
c_A = \bar c_A - \mathcal H s(0) + q_J, \, \, c_\ell = \bar c_\ell - a \mathcal H s(0) + a q_J.
\]
In this section the modulation errors are measured from the nearby profile.
\[
\delta_A := c_A - \bar c_A = - \mathcal H s(0) + q_J, \qquad
\delta_\ell := c_\ell - \bar c_\ell = - a \mathcal H s(0) + a q_J.
\]
These agree with the earlier definitions when $a = 1/2$. The origin trace equation is 
\[
\partial_\tau s_y (0, \tau ) = ( 1 - a) q_J s_y (0, \tau),
\]
and $s_y (0, 0) = 0$ implies $s_y (0, \tau ) = 0$. The linearized operators at $(F_a, 0)$ are
\[
\mathcal L_{s, a} s = - V_a s' + (\bar c_A + Q_a) s + F_a \mathcal H s - a F_a' \mathcal P s + \mathcal H s(0) (- F_a + a y F_a'),
\]
\[
\mathcal L_a h = - V_a h' + (\bar c_A - Q_a) h + F_a \mathcal H h + a F_a' \mathcal P h.
\]
The magnetic forcing is
\begin{equation} \label{4.4}
G_{J, a } = - J \mathcal H J + a B J' + q_J (F_a - a y F_a')
\end{equation}
where
\[
G_{J, a } (0) = G'_{J, a}(0) = 0. 
\]
Since $V/y$ is bounded above and below by positive constants, the domain from Proposition~\ref{2.1} can also be written as
\[
\mathcal D_\sigma = \{ h \in X : y h' \in X \}.
\]
\begin{Proposition}\label{4.5}
For sufficiently small $\left |a - \frac{1}{2} \right |$, the operator $\mathcal L_a$ on $X$ with domain $\mathcal D_\sigma$ has an isolated eigenvalue $\lambda_a$ of finite multiplicity such that
\[
\lim_{a \to \frac{1}{2}} \lambda_a = \lambda_*.
\]
The eigenvalue can be chosen so that
\[
\operatorname{Re} \lambda_a > \bar c_A + \bar c_\ell \left (\sigma + \frac{1}{2} \right),
\]
\[
\frac{27}{40} + \operatorname{Re} \lambda_a > 0, \qquad
\gamma_a = \frac{\alpha_a - \bar c_\ell + \operatorname{Re} \lambda_a}{\alpha_a} > \frac{1}{72}.
\]
\end{Proposition}

\begin{proof}
We realize $\mathcal L_a$ on the common domain $\mathcal D_\sigma$. By Lemma~\ref{4.1},
\[
c \| y h' \|_X \leq \| V_a h' \|_X \leq C \| y h' \|_X.
\]
As in Proposition~\ref{2.1},
\[
\|(\mathcal L_a + V_a \partial_y) h \|_X \leq C \| h \|_X.
\]
Hence
\[
c \bigl (\| h \|_X + \| y h' \|_X \bigr)
\leq \| h \|_X + \| \mathcal L_a h \|_X
\leq C \bigl (\| h \|_X + \| y h' \|_X \bigr).
\]
Completeness of $\mathcal D_\sigma$ in this norm implies that $\mathcal L_a$ is closed. Lemma~\ref{4.1} also bounds the operator difference.
\begin{equation} \label{4.6}
\| (\mathcal L_a - \mathcal L) h \|_X
\leq C \left |a - \frac12 \right |(\| h \|_X + \| y h' \|_X).
\end{equation}
For the potential term, use the bounds
\[
| \mathcal P h(y)| \leq C \langle y \rangle^{- 1/12} \| h \|_X,
\qquad
\left \| \langle y \rangle^{1/2} \left (a F_a' - \frac12 F' \right) \right \|_2
\leq C \left |a - \frac12 \right |.
\]

Choose a positively oriented circle $\Gamma$ enclosing only $\lambda_*$ in the spectrum of $\mathcal L$, with
\[
\Gamma \subset \{ z : \operatorname{Re} z > - 23/36 \}, \qquad
M_\Gamma := \sup_{z \in \Gamma}
\|(z - \mathcal L)^{- 1} \|_{X \to \mathcal D_\sigma} < \infty.
\]
By \eqref{4.6}, for $C|a - 1/2|M_\Gamma < 1/2$,
\[
\sup_{z \in \Gamma}
\|(\mathcal L_a - \mathcal L)(z - \mathcal L)^{- 1} \|_{X \to X} < \frac12.
\]
Thus
\[
z - \mathcal L_a
= \left [I - (\mathcal L_a - \mathcal L)(z - \mathcal L)^{- 1} \right](z - \mathcal L)
\]
is invertible uniformly on $\Gamma$. The Riesz projections
\[
P_a = \frac1{2 \pi i} \int_\Gamma (z - \mathcal L_a)^{- 1} \, d z
\]
therefore satisfy
\[
\| P_a - P_{1/2} \|_{X \to X} \longrightarrow 0 \qquad (a \to 1/2),
\]
and, for $a$ sufficiently close to $1/2$,
\[
1 \leq \operatorname{rank} P_a = \operatorname{rank} P_{1/2} < \infty.
\]
The enclosed spectrum consists of isolated eigenvalues of finite multiplicity. The same argument works on circles of radius $r$ about $\lambda_*$.
\[
\begin{gathered}
\forall r > 0 \  \hbox{sufficiently small}, \quad
\exists \delta_r > 0, \\
|a - 1/2| < \delta_r \quad \Longrightarrow \quad
\operatorname{spec}(\mathcal L_a) \cap \{ z : |z - \lambda_*| < r \} \neq \varnothing.
\end{gathered}
\]
We may therefore choose an eigenvalue with $\lambda_a \to \lambda_*$.

The rate bounds in Lemma~\ref{4.1} imply
\[
\bar c_A + \bar c_\ell \left (\sigma + \frac12 \right) \longrightarrow - \frac{23}{36},
\qquad
\gamma_a \longrightarrow \frac23 + \lambda_* > \frac1{36}.
\]
Since
\[
\lambda_* + \frac{23}{36} > 0, \qquad
\frac{27}{40} + \lambda_* > \frac{13}{360} > 0, \qquad
\frac1{36} > \frac1{72},
\]
the claimed strict inequalities hold after reducing the parameter interval.
\end{proof}

\subsection{Energy estimates} 

Let $d = \left |a - \frac{1}{2} \right |$. The weights $\phi$, $\psi$, and $\varrho$ are the same as in Section 3.

\begin{Lemma}\label{4.7}
There exist $\delta > 0$ and $C > 0$, independent of $a$ for $d \leq \delta$, such that
\[
\operatorname{Re} \langle h, \mathcal L_a h \rangle
\leq - \left (\frac{1}{4} - C d \right) \| h \|_2^2,
\]
\[
\operatorname{Re} \int (\mathcal L_a h) \bar h \varrho
\leq - \left (\frac{6}{25} - C d \right) \| h \|_{L^2(\varrho)}^2,
\]
for $h \in \mathcal D_\sigma$. For smooth odd $s$ with $s'(0) = 0$ and $\| s \|_{L^2(\phi)} + \| s^{(4)} \|_{L^2(\psi)} < \infty$,
\[
\operatorname{Re} \int (\mathcal L_{s, a} s) \bar s \phi
\leq - \left (\frac{1}{2} - C d \right) \| s \|_{L^2(\phi)}^2.
\]
\end{Lemma}

\begin{proof}
All constants below are uniform for $d \leq \delta \leq \frac14$. We compare the quadratic forms with their counterparts at $a = 1/2$.

\emph{Magnetic forms.} By Lemma~\ref{4.1},
\[
\| V_a' - V' \|_\infty + \left \|  \frac{V_a - V}{y}  \right \|_\infty + | \bar c_A + 1 | + \| Q_a - Q \|_\infty  + \| F_a - F \|_\infty \leq C d,
\]
and, for $\rho = 1, \varrho$, 
\[
\|   \sqrt{|y|} \rho^{\frac{1}{2}} \left ( a F_a' - \frac{F'}{2} \right)   \|_2 \leq C d. 
\]
The weight satisfies the bounds in Proposition~\ref{2.15}.
\[
c \langle y \rangle^{2 \sigma} \leq \varrho (y) \leq C \langle y \rangle^{2 \sigma}, \qquad
\left \| \frac{y \varrho'}{\varrho} \right \|_\infty < \infty.
\]
We have
\[
\begin{aligned}
(\mathcal L_a - \mathcal L ) h & = - \left ( V_a - V \right) h' + \left (  \bar c_A + 1 - (Q_a - Q ) \right) h  \\
& + \left ( F_a - F \right) \mathcal H h + \left ( a F_a' - \frac{F'}{2} \right) \mathcal P h,
\end{aligned}
\]
and
\[
\operatorname{Re} \int ( - \left ( V_a - V \right) h') \bar h \rho = \frac{1}{2} \int \left ( \left ( V_a' - V' \right) + \frac{ \left ( V_a - V \right) \rho'}{\rho} \right) |h|^2 \rho \leq C d \| h \|^2_{L^2 (\rho)}.
\]
Since $| \mathcal P h (y)| \leq \sqrt{|y|} \| h \|_2$, 
\[
\begin{aligned}
& \| \left ( \bar c_A + 1 - (Q_a - Q) \right) h \|_{L^2 (\rho)}
+ \| \left ( F_a - F \right) \mathcal H h \|_{L^2 (\rho)} \\
& \quad + \| \left ( a F_a' - \frac{F'}{2} \right) \mathcal P h \|_{L^2 (\rho)}
\leq C d \| h \|_{L^2 (\rho)}.
\end{aligned}
\]
These form bounds, together with \eqref{2.2} and \eqref{2.16}, prove the two magnetic estimates.

\emph{Scalar form.} For odd $s$ satisfying $s'(0) = 0 $, the scalar operator difference is
\[
\begin{aligned}
(\mathcal L_{s, a} - \mathcal L_s) s & = - \left ( V_a - V \right) s' + \left (  \bar c_A + 1 + (Q_a - Q )  \right) s + \left ( F_a - F \right) \left ( \mathcal H s - \mathcal H s(0) \right) \\
& - \left ( a F_a' - \frac{F'}{2} \right) \left ( \mathcal P s - y \mathcal H s(0) \right),
\end{aligned}
\]
\[
\begin{aligned}
& \left |  \operatorname{Re} \int (- \left ( V_a - V \right) s') \bar s \phi \right | \\
& \leq \frac{1}{2} \left \| \left ( V_a' - V' \right) + \left ( V_a - V \right) \frac{\phi '}{\phi } \right \|_\infty \| s \|_{L^2 (\phi)}^2  \leq C d \| s \|_{L^2 (\phi)}^2,
\end{aligned}
\]
and 
\[
| \mathcal H s(0) | + \| \mathcal P s \|_\infty + \| \langle y \rangle \mathcal H s \|_2 \leq C \| s \|_{L^2 (\phi)}
\]
where
\[
\frac{\mathcal H s - \mathcal H s(0)}{y} = \mathcal H \left ( \frac{s}{y} \right), \qquad
\mathcal P s(y) - y \mathcal H s(0) = \int_0^y \left ( \mathcal H s(x) - \mathcal H s(0) \right) d x,
\]
and 
\[
\left \|  \frac{\mathcal H s - \mathcal H s(0)}{y}   \right \|_2 + \left \|  \frac{\mathcal P s - y \mathcal H s(0)}{y^2}   \right \|_2  \leq C \| s \|_{L^2 (\phi)},
\]
For the transport coefficient,
\[
\frac{y \phi'} \phi = \frac{6 y^2}{1 + y^2} - 4 \in L^\infty.
\]
For $|y| \leq 1$,
\[
|(F_a - F)(y)| \leq C d|y|^3, \qquad
|(F_a' - F')(y)| \leq C d|y|^2.
\]
Consequently,
\[
\begin{aligned}
& \| \left ( F_a - F \right) \left ( \mathcal H s - \mathcal H s(0) \right) \|_{L^2 ( \phi; |y| \leq 1 )} \\
& \leq \left \|   y \phi^{\frac{1}{2}} \left ( F_a - F \right)    \right \|_{L^\infty (|y| \leq 1 )} \left \|  \frac{\mathcal H s - \mathcal H s(0)}{y}  \right \|_2  \leq C d \| s \|_{L^2 (\phi)},
\end{aligned}
\]
\[
\begin{aligned}
& \| \left ( F_a' - F' \right) \left ( \mathcal P s - y \mathcal H s(0) \right) \|_{L^2 ( \phi; |y| \leq 1 )} \\
& \leq \left \|   y^2 \phi^{\frac{1}{2}} \left ( F_a' - F' \right)    \right \|_{L^\infty (|y| \leq 1 )} \left \|  \frac{\mathcal P s - y \mathcal H s(0)}{y^2}  \right \|_2  \leq C d \| s \|_{L^2 (\phi)},
\end{aligned}
\]
\[
\begin{aligned}
& \| \left ( F_a - F \right) \left ( \mathcal H s - \mathcal H s(0) \right) \|_{L^2 ( \phi; |y| \geq 1 )} \\
& \leq C \left (   \| F_a - F \|_\infty  \| \langle y \rangle \mathcal H s \|_2 + | \mathcal H s(0)| \|  \langle y \rangle \left ( F_a - F \right)     \|_2  \right)  \leq C d \| s \|_{L^2 (\phi)},
\end{aligned}
\]
\[
\begin{aligned}
& \| \left ( F_a' - F' \right) \left ( \mathcal P s - y \mathcal H s(0) \right) \|_{L^2 ( \phi; |y| \geq 1 )} \\
& \leq C \left (   \| \mathcal P s \|_\infty  \| \langle y \rangle \left ( F_a' - F' \right) \|_2 + | \mathcal H s(0)| \|  \langle y \rangle^2 \left ( F_a' - F' \right)     \|_2  \right)  \leq C d \| s \|_{L^2 (\phi)},
\end{aligned}
\]
and
\[
d \| F ' \left ( \mathcal P s - y \mathcal H s(0) \right) \|_{L^2 (\phi)} \leq d \left \| y^2 \phi^{\frac{1}{2}} F' \right \|_\infty \left \| \frac{\mathcal P s - y \mathcal H s(0)}{y^2} \right \|_2 \leq C d \| s \|_{L^2 (\phi)}.
\]
Thus the scalar form differs from its base form by at most $C d \| s \|_{L^2(\phi)}^2$. Lemma~\ref{3.1}, applied to the real and imaginary parts, proves the scalar estimate.
The magnetic inequalities extend as quadratic form inequalities to $h \in X$ by density after integration by parts. The scalar pairing has the cutoff interpretation of Lemma~\ref{3.1}. The cutoff fluxes vanish since $V_a/y$ is bounded.
\end{proof}

Let
\[
E_4^2 = \| s^{(4)} \|_{L^2 (\psi)}^2 + \| J^{(4)} \|_{L^2 (\psi)}^2.
\]
With $\mu$ fixed as in Proposition~\ref{3.11}, define
\[
E^2 = \| s \|_{L^2 (\phi)}^2 + \| J \|_{L^2 (\varrho)}^2 + \mu E_4^2.
\]
\begin{Lemma}\label{4.8}
There exist $\delta > 0$, $E_* > 0$, and positive constants $C_0, C_2$, independent of $a$ for $d \leq \delta$, such that every smooth odd normalized solution with $s'(0) = 0$ and $E \leq E_*$ satisfies
\[
\frac{1}{2}(E^2)' \leq - \left (\frac{23}{100} - C_0 d \right) E^2 + C_2 E^3.
\]
Here $\mu > 0$ is fixed as in Proposition~\ref{3.11}.
\end{Lemma}

\begin{proof}
\emph{Linear form comparison.} Lemma~\ref{4.7} controls the order-zero terms. The fourth-order transport coefficients differ from their base values by
\[
\bar c_A + 1 \pm (Q_a - Q) - \frac72(V_a' - V')
+ \frac12(V_a - V) \frac{\psi'}{\psi},
\]
whose $L^\infty$ norm is at most $C d$. The plus sign corresponds to $s$, and the minus sign to $J$. Inserting \eqref{3.6} and \eqref{3.7},
\[
\bar c_A + Q_a - \frac72 V_a' + \frac12 V_a \frac{\psi'}{\psi}
\leq - \frac9{10} + C d, \quad \bar c_A - Q_a - \frac72 V_a' + \frac12 V_a \frac{\psi'}{\psi}
\leq - \frac{19}{78} + C d.
\]
The principal Hilbert forms satisfy
\[
\int \psi F f^{(4)} \mathcal H f^{(4)} = 0,
\]
and
\[
\left | \int \psi \left ( F_a - F \right) f^{(4)} \mathcal H f^{(4) }    \right | \leq C \| F_a - F \|_\infty  \| f^{(4)} \|_{L^2 (\psi)}^2 \leq C d E_4^2. 
\]
Indeed, for $f \in \{ s, J \}$, the cutoff argument in Lemma~\ref{1.21} verifies the required parity and moment conditions.
\[
f^{(4)}(- y) = - f^{(4)}(y), \qquad \int y f^{(4)}(y) \, d y = 0.
\]
Apply Lemma~\ref{1.20}(3) with $\eta = 2$.
\[
\| \mathcal H f^{(4)} \|_{L^2(\psi)}
\leq C \| f^{(4)} \|_{L^2(\psi)}.
\]
The remaining lower-order linear terms are combinations of
\[
\begin{gathered}
V_a^{(j)} f^{(5 - j)} \quad (2 \leq j \leq 4), \quad
Q_a^{(j)} f^{(4 - j)} \quad (1 \leq j \leq 4), \quad
F_a^{(j)} \mathcal H f^{(4 - j)} \quad (1 \leq j \leq 4), \\
F_a^{(j + 1)} \mathcal H f^{(3 - j)} \quad (0 \leq j \leq 3), \quad
F_a^{(5)} \mathcal P f, \quad
\mathcal H s(0) \partial_y^4(- F_a + a y F_a') \quad (f = s).
\end{gathered}
\]
Keep $\mu$ fixed. For the transport coefficients, Lemma~\ref{4.1} implies
\[
\|(V_a' - V') \|_\infty + \|(V_a - V)/y \|_\infty \leq C d,
\]
and the coefficient differences through order five satisfy \eqref{4.3}. In particular,
\[
\left \| \sqrt{|y|} \, \psi^{1/2} \partial_y^4 \left (a F_a' - \tfrac12 F' \right) \right \|_2 + \left \| \psi^{1/2} \partial_y^4 \left [- (F_a - F) + y \left (a F_a' - \tfrac12 F' \right) \right] \right \|_2 \leq C d.
\]
Integration by parts with the cutoff from Proposition~\ref{3.21} eliminates the fifth derivative.
\[
\left | \int \psi f^{(4)}[- (V_a - V) f^{(5)}] \right | = \frac12 \left | \int [\psi (V_a - V)]'|f^{(4)}|^2 \right | \leq C d \| f^{(4)} \|_{L^2(\psi)}^2.
\]
Let $\mathcal C_a$ be the finite set of coefficient differences
\[
\begin{gathered}
\partial_y^j(V_a - V) \quad (2 \leq j \leq 4), \qquad
\partial_y^j(Q_a - Q) \quad (1 \leq j \leq 4), \\
\partial_y^j(F_a - F) \quad (1 \leq j \leq 4), \qquad
\partial_y^j(a F_a' - \tfrac12 F') \quad (0 \leq j \leq 3).
\end{gathered}
\]
The weighted coefficient bounds and interpolation imply, for $c \in \mathcal C_a$ and $0 \leq k \leq 3$,
\[
\| \psi^{1/2} c \|_\infty \leq C d, \qquad \| c f^{(k)} \|_{L^2(\psi)}
+ \| c \mathcal H f^{(k)} \|_{L^2(\psi)}
\leq C d \bigl (\| f \|_2 + \| f^{(4)} \|_{L^2(\psi)} \bigr).
\]
The two endpoints satisfy
\[
\begin{aligned}
& \| \partial_y^4(a F_a' - \tfrac12 F') \mathcal P f \|_{L^2(\psi)} \leq \| \sqrt{|y|} \psi^{1/2} \partial_y^4(a F_a' - \tfrac12 F') \|_2 \| f \|_2
\leq C d \| f \|_2, \\
& | \mathcal H s(0)| \,
\| \partial_y^4[- (F_a - F) + y(a F_a' - \tfrac12 F')] \|_{L^2(\psi)} \leq C d \| s \|_{L^2(\phi)}.
\end{aligned}
\]
The principal Hilbert estimate completes the form bound.
\[
\begin{split}
& \left | \operatorname{Re} \int \psi s^{(4)} \partial_y^4 \bigl [(\mathcal L_{s, a} - \mathcal L_s) s \bigr] \right |
+ \left | \operatorname{Re} \int \psi J^{(4)} \partial_y^4 \bigl [(\mathcal L_a - \mathcal L) J \bigr] \right | \\
& \qquad \leq C d \left (E_4^2 + \| s \|_{L^2(\phi)}^2 + \| J \|_{L^2(\varrho)}^2 \right).
\end{split}
\]
Thus \eqref{3.10} retains its original constant $C_1$, with only an additional $C d$ in the lower-order coefficient.

\emph{Nonlinear terms.} The potentially derivative-losing magnetic transport terms cancel at the fifth derivative level.

\[
a \int \psi B ( s^{(4)} J^{(5)} + J^{(4)} s^{(5)} ) = - a \int (\psi B)' s^{(4)} J^{(4)}, 
\]
and
\[
\left |  a \int (\psi B)' s^{(4)} J^{(4)}   \right | \leq C \Bigg ( \| B' \|_\infty + \left \|   \left ( \frac{\psi ' }{\psi}  \right)  B \right \|_\infty \Bigg) E_4^2 \leq C E E_4^2. 
\]
For $f, h \in \{ s, J  \}$, apply Lemma~\ref{1.21} whenever $ 1 \leq j \leq 3$. 
\[
\| f^{(j)} \mathcal H h^{(4 - j)} \|_{L^2 (\psi )} \leq \|   \langle y \rangle^{\mathfrak b_j} f^{(j)}      \|_\infty \|   \langle y \rangle^{\mathfrak a_{4 - j}}    \mathcal H h^{(4 - j)}   \|_2 \leq C E^2,
\]
where $\mathfrak b_j + \mathfrak a_{4 - j} = \frac{265}{96} > 2$, and 
\[
\|  f \mathcal H h^{(4)}    \|_{L^2 (\psi )} + \|  f^{(4)} \mathcal H h    \|_{L^2 (\psi )} \leq C E^2. 
\]
The forcing term \eqref{4.4} has the origin cancellations
\[
G_{J, a } (0) = G'_{J, a } (0) = 0, 
\]
and satisfies the bounds
\[ 
\left \| \frac{G_{J, a}}{y^2}  \right \|_{L^2 { (|y| \leq 1   )}} \leq C \|  G_{J, a} '' \|_{L^2{(|y| \leq 1) }} \leq C E^2,
\]
\[
\|     |y| B J'      \|_{L^2{(|y| \geq 1) }} \leq \|  |y|^{\frac{1}{12}} B \|_{L^\infty{(|y| \geq 1) }}   \|   \langle    y  \rangle^{\mathfrak a_1} J'     \|_2   \leq C E^2,
\]
\[
\|     |y| J \mathcal H J      \|_{L^2{(|y| \geq 1) }}  + |q_J| \|    F_a - a y F'_a      \|_{L^2 (\phi; |y| \geq 1 )} \leq C E^2. 
\]
The inequalities $\frac1{12} + \mathfrak a_1 = \frac{49}{48} > 1$ and $\mathfrak a_0 + \mathfrak b_0 = \frac{43}{32} > 1$ therefore give $\| G_{J, a} \|_{L^2(\phi)} \leq C E^2$.

\emph{Combination with the fixed weight.} The three energy inequalities are
\[
\frac{1}{2} (\| s \|_{L^2 (\phi)}^2)' \leq - \left (  \frac{1}{2} - C d      \right) \| s \|_{L^2 (\phi)}^2 + C E^3,
\]
\[
\frac{1}{2} (\| J \|_{L^2 (\varrho)}^2)' \leq - \left (  \frac{6}{25} - C d      \right) \| J \|_{L^2 (\varrho)}^2 + C E^3,
\]
\[
\frac{1}{2} (E_4^2)' \leq - \left (  \frac{23}{100} - C d      \right) E_4^2 + (C_1 + C d) (\| s \|_{L^2 (\phi)}^2 + \| J \|_{L^2 (\varrho)}^2) + C E^3,
\]
Add the first two inequalities to $\mu$ times the third and use the choice of $\mu$ in Proposition~\ref{3.11}. This proves the asserted estimate.
Here $C_1$ is the constant in \eqref{3.10}. The constants $C_0, C_2$ may depend on the fixed $\mu$ and are uniform in $a$.
\end{proof}

For the fixed $\mu_3 > 0$ in Proposition~\ref{3.21}, let
\[
M_3^2 = \| s \|_{L^2(\phi)}^2 + \mu_3 \| s^{(3)} \|_{L^2(\psi)}^2.
\]

\begin{Proposition}\label{4.9}
There exists $\delta > 0$ such that for $\left |a - \frac{1}{2} \right | < \delta$ and every fixed real odd $f_a \in C_c^\infty (\R)$, there exists $\varepsilon_0(a, f_a) > 0$ with the following property. For $0 < \varepsilon < \varepsilon_0(a, f_a)$, the initial data
\[
(s, J)( \cdot , 0) = (0, \varepsilon f_a)
\]
give a smooth odd normalized solution for all $\tau \geq 0$, unique among odd solutions in $C_{\mathrm{loc}}([0, \infty); Z \times Z)$, with $s_y(0, \tau) = 0$. There exists $C_a > 0$, depending on $a$ and $f_a$ but independent of $\varepsilon$ and $\tau$, such that
\begin{equation} \label{4.10}
E (\tau) \leq C_a \varepsilon e^{- \frac{9 \tau }{40}}.
\end{equation}
\begin{equation} \label{4.11}
M_3 (\tau ) \leq C_a \varepsilon^2 e^{- \frac{9 \tau }{20}}. 
\end{equation}
\end{Proposition}

\begin{proof}
The local construction in Section~\ref{3.14} applies with the symmetric transport matrix
\[
a \begin{pmatrix} u  & - B \\ - B  & u \end{pmatrix}.
\]
It follows from Lemma~\ref{4.1} 
\[
(F_a, \varepsilon f_a) \in (Z \times Z) \cap (Y_m \times Y_m)
\quad \hbox{for every } m, \qquad F_a^{(3)} \in L^2(\psi).
\]
The physical-time estimate in Proposition~\ref{3.21} applies with the preceding weighted product bounds. On each compact local interval $I$, the scales satisfy
\[
0 < \inf_I A \leq \sup_I A < \infty, \qquad
0 < \inf_I L \leq \sup_I L < \infty.
\]
After rescaling,
\[
s^{(3)} \in C(I; L^2(\psi)), \qquad M_3 \in C(I), \qquad M_3(0) = 0.
\]

\emph{Total energy.} Choose $\delta$ and $0 < E_b \leq E_*$ so that
\[
C_0 \delta \leq \frac1{800}, \qquad C_2 E_b \leq \frac1{800}.
\]
For the fixed seed, choose
\[
0 < \varepsilon \leq
\frac{E_b}{2 \left [1 + \left (\| f_a \|_{L^2(\varrho)}^2
+ \mu \| f_a^{(4)} \|_{L^2(\psi)}^2 \right)^{1/2} \right]}.
\]
This choice also covers $f_a = 0$, and
\[
E(0) = \varepsilon \left (\| f_a \|_{L^2(\varrho)}^2
+ \mu \| f_a^{(4)} \|_{L^2(\psi)}^2 \right)^{1/2}
\leq \frac{E_b}{2}.
\]
Lemma~\ref{4.8} applies as long as $E \leq E_b$.
\[
\frac12(E^2)' \leq - \frac{91}{400} E^2,
\qquad
E(\tau) \leq E(0) e^{- 91 \tau /400} \leq \frac{E_b}{2}.
\]
Continuity closes the bootstrap on the interval of existence. Since $91/400 > 9/40$, this proves \eqref{4.10} there.

\emph{Improved scalar decay.} The third-order transport coefficient from Proposition~\ref{3.21} now satisfies
\[
\bar c_A + Q_a - \frac52 V_a' + \frac12 V_a \frac{\psi'}{\psi}
\leq - \frac{13}{18} + C d.
\]
As in that proposition, $s^{(3)}$ is even with zero integral. Apply Lemma~\ref{1.20}(2) with $\eta = 2$.
\[
\left | \int \psi (F_a - F) s^{(3)} \mathcal H s^{(3)} \right |
\leq C d \| s^{(3)} \|_{L^2(\psi)}^2, \qquad \| a B J^{(4)} \|_{L^2(\psi)} \leq C E^2.
\]
For the fixed $\mu_3$, the coefficient bounds through order four and the interpolation estimate in Proposition~\ref{3.21} control the scalar perturbation.
\[
\begin{split}
\Bigl |  & \operatorname{Re} \int [(\mathcal L_{s, a} - \mathcal L_s) s] \overline s \phi \\
& + \mu_3 \operatorname{Re} \int \psi \overline{s^{(3)}} \partial_y^3[(\mathcal L_{s, a} - \mathcal L_s) s] \Bigr |
\leq C_{\mu_3} d M_3^2.
\end{split}
\]
The base linear rate $49/100$ and the preceding bounds give
\[
\frac12(M_3^2)'
\leq - \left (\frac{49}{100} - C_{\mu_3} d - C_{\mu_3} E \right) M_3^2 + C E^2 M_3.
\]
Reduce $\delta, E_b$ subject to the preceding conditions and
\[
C_{\mu_3}(\delta + E_b) \leq \frac1{100}.
\]
Choose $\varepsilon$ by the bound above. Then
\[
M_3' \leq - \frac{12}{25} M_3 + C E^2,
\qquad
M_3(\tau) \leq C_a \varepsilon^2 \int_0^\tau
e^{- 12(\tau - r)/25 - 9 r/20} \, d r.
\]
Since $12/25 - 9/20 = 3/100 > 0$, this proves \eqref{4.11}.

\emph{Continuation.} The weighted energy calculation in \eqref{3.18} remains valid.
\[
\frac{d}{d \tau} \|(\omega, J) \|_{Y_4}
\leq C_a \|(\omega, J) \|_{Y_4}
\]
while $E \leq E_b$. Moreover,
\[
\|(\omega, J) \|_Z \leq C_a \left (
\|(\omega, J) \|_{Y_4} + E + \| \langle y \rangle^2 F_a^{(4)} \|_2 \right).
\]
For every finite $\tau_0$ in the existence interval,
\[
\sup_{0 \leq \tau \leq \tau_0} \|(\omega, J)(\tau) \|_{Y_4}
\leq e^{C_a \tau_0} \|(F_a, \varepsilon f_a) \|_{Y_4}, \qquad
\sup_{0 \leq \tau \leq \tau_0} \|(\omega, J)(\tau) \|_Z < \infty.
\]
By the continuation and regularity results in Section~\ref{3.14}, the normalized solution is global, smooth, and unique in the stated class. Estimates \eqref{4.10} and \eqref{4.11} hold for all $\tau \geq 0$.
\end{proof}

\subsection{Proof of finite-time blowup}

\begin{Lemma}\label{4.12}
For the eigenvalue $\lambda_a$ in Proposition~\ref{4.5}, there exists a nonzero $\ell_a \in X^*$ such that
\[
\ell_a(\mathcal L_a h) = \lambda_a \ell_a(h), \qquad h \in \mathcal D_\sigma.
\]
It has an odd representative $\zeta_a \in H^1(\R; \C)$ satisfying
\[
y \zeta_a' \in L^2(\R), \qquad (\mathcal H \zeta_a)' \in L^1(\R),
\]
and, for every $h \in X$,
\[
\ell_a(h) = \int \zeta_a h = - \int (\mathcal H \zeta_a)' \mathcal P h,
\qquad | \ell_a(h)| \leq C_a \| \mathcal P h \|_\infty.
\]
\end{Lemma}

\begin{proof}
Proposition~\ref{4.5} supplies a nonzero left eigenfunctional $\ell_a \in X^*$. Write its odd Riesz representative as
\[
\ell_a(h) = \int \zeta_a h,
\qquad \| \langle y \rangle^{- \sigma} \zeta_a \|_2 < \infty.
\]
Testing against compactly supported odd functions identifies the adjoint equation.
\[
R_a = V_a \zeta_a' - \left [\lambda_a - \bar c_A - \bar c_\ell + (1 - a) Q_a \right] \zeta_a = \mathcal H(F_a \zeta_a) - a \mathcal P(F_a' \zeta_a).
\]
The coefficient bounds in Lemma~\ref{4.1} imply
\[
\begin{gathered}
\| F_a \zeta_a \|_2 + \| F_a' \zeta_a \|_2 < \infty, \qquad \| y F_a' \zeta_a \|_1
\leq \| \langle y \rangle^{- \sigma} \zeta_a \|_2
\| y \langle y \rangle^\sigma F_a' \|_2 < \infty, \\
\| \mathcal P(F_a' \zeta_a) \|_{H^1}
\leq C \bigl (\| y F_a' \zeta_a \|_1 + \| F_a' \zeta_a \|_2 \bigr).
\end{gathered}
\]
In particular, $R_a \in L^2$.

\emph{Origin.} By oddness and Hardy's inequality,
\[
\mathcal H(F_a \zeta_a) = y \mathcal H \left (\frac{F_a}{y} \zeta_a \right),
\qquad
\left \| \frac{\mathcal P(F_a' \zeta_a)} y \right \|_{L^2(0, 1)}
\leq 2 \| \mathcal H(F_a' \zeta_a) \|_{L^2(0, 1)}.
\]
Thus $R_a/V_a \in L^2(0, 1)$. The origin and infinity exponents are
\[
\alpha_{\mathrm{adj}, a}
= \frac{\lambda_a + 2 Q_a(0) - 10/3}{5/3},
\qquad
\beta_{\mathrm{adj}, a}
= \frac{\lambda_a - \bar c_A - \bar c_\ell}{\bar c_\ell}.
\]
By Lemma~\ref{4.1} and Proposition~\ref{4.5}, after decreasing the parameter interval if necessary,
\[
\operatorname{Re} \alpha_{\mathrm{adj}, a} > \frac12,
\qquad
\operatorname{Re} \beta_{\mathrm{adj}, a} > \sigma - \frac12.
\]
For $0 < |y| < 1$, set
\[
b_{a, 0}(y) =
\frac{\lambda_a - \bar c_A - \bar c_\ell + (1 - a) Q_a(y)}{V_a(y)}
- \frac{\alpha_{\mathrm{adj}, a}} y.
\]
Then
\[
|b_{a, 0}(y)| \leq C_a|y|, \qquad
g_{a, 0} := \frac{R_a}{V_a} + b_{a, 0} \zeta_a \in L^2(- 1, 1),
\qquad
\zeta_a' - \frac{\alpha_{\mathrm{adj}, a}} y \zeta_a = g_{a, 0}.
\]
On $(0, 1)$, for some $c_{a, 0} \in \C$,
\[
\zeta_a(y) = c_{a, 0} y^{\alpha_{\mathrm{adj}, a}}
- y^{\alpha_{\mathrm{adj}, a}}
\int_y^1 t^{- \alpha_{\mathrm{adj}, a}} g_{a, 0}(t) \, d t.
\]
As in Proposition~\ref{2.17},
\[
\left \| y^{\alpha_{\mathrm{adj}, a} - 1}
\int_y^1 t^{- \alpha_{\mathrm{adj}, a}} g_{a, 0}(t) \, d t \right \|_{L^2(0, 1)}
\leq \frac{\| g_{a, 0} \|_{L^2(0, 1)}}
{\operatorname{Re} \alpha_{\mathrm{adj}, a} - 1/2}.
\]
Thus $\zeta_a/y, \zeta_a' \in L^2(0, 1)$. The odd extension satisfies
\[
\zeta_a \in H^1(- 1, 1), \qquad \zeta_a(0) = 0.
\]

\emph{Infinity.} Set
\[
r_a = \frac{y R_a}{V_a}
+ \left [\frac{y(\lambda_a - \bar c_A - \bar c_\ell + (1 - a) Q_a)}{V_a}
- \beta_{\mathrm{adj}, a} \right] \zeta_a.
\]
For $y \geq 1$, the bounds $U_a \in L^\infty$ and $|Q_a(y)| \leq C_a \langle y \rangle^{- 1}$ control the coefficient.
\[
\left | \frac{y(\lambda_a - \bar c_A - \bar c_\ell + (1 - a) Q_a)}{V_a}
- \beta_{\mathrm{adj}, a} \right | \leq \frac{C_a} y.
\]
Since $\sigma < 1$,
\[
\| r_a \|_{L^2(1, \infty)}
\leq C_a \left (\| R_a \|_2 +
\| \langle y \rangle^{- \sigma} \zeta_a \|_2 \right) < \infty.
\]
The adjoint equation becomes
\[
\zeta_a' - \frac{\beta_{\mathrm{adj}, a}} y \zeta_a = \frac{r_a} y.
\]
The homogeneous solution is excluded by
\[
\int_1^\infty y^{2 \operatorname{Re} \beta_{\mathrm{adj}, a} - 2 \sigma} \, d y = \infty,
\qquad \langle y \rangle^{- \sigma} \zeta_a \in L^2.
\]
Hence
\[
\zeta_a(y) = - y^{\beta_{\mathrm{adj}, a}}
\int_y^\infty t^{- \beta_{\mathrm{adj}, a} - 1} r_a(t) \, d t.
\]
By the tail Hardy estimate,
\[
\begin{aligned}
\| \zeta_a \|_{L^2(1, \infty)}
& \leq \frac{\| r_a \|_{L^2(1, \infty)}}{\operatorname{Re} \beta_{\mathrm{adj}, a} + 1/2}, \\
\| y \zeta_a' \|_{L^2(1, \infty)}
& \leq | \beta_{\mathrm{adj}, a}| \| \zeta_a \|_{L^2(1, \infty)}
+ \| r_a \|_{L^2(1, \infty)}.
\end{aligned}
\]
These bounds also hold on $(- \infty, - 1)$ by oddness. Therefore
\[
\zeta_a \in H^1(\R), \quad y \zeta_a' \in L^2, \qquad
\zeta_a'(- y) = \zeta_a'(y), \qquad
\int \zeta_a' = \zeta_a(+ \infty) - \zeta_a(- \infty) = 0.
\]
The derivative of the Hilbert transform is integrable by Lemma~\ref{1.20}(2) with $1/2 < \eta < 1$.
\[
\| \mathcal H \zeta_a' \|_1
\leq \| \langle y \rangle^{- \eta} \|_2
\| \langle y \rangle^\eta \mathcal H \zeta_a' \|_2
\leq C_\eta \| \langle y \rangle \zeta_a' \|_2 < \infty.
\]
Since $\mathcal H \zeta_a \in H^1$ and $\mathcal P h \in L^\infty$,
\[
\ell_a(h) = \int \mathcal H \zeta_a(\mathcal P h)', \qquad
\lim_{y \to \pm \infty} \mathcal H \zeta_a(y) \mathcal P h(y) = 0.
\]
Then we obtain
\begin{equation} \label{4.13}
\ell_a(h) = - \int (\mathcal H \zeta_a)' \mathcal P h,
\qquad
| \ell_a(h)| \leq C_a \| \mathcal P h \|_\infty.
\end{equation}
\end{proof}

\begin{proof}[Proof of Theorem~\ref{1.6}]
\emph{Prepared data and projected growth.} By Lemma~\ref{4.12}, choose
\[
f_a \in C_c^\infty (\R; \R), \qquad f_a(- y) = - f_a(y), \qquad \ell_a(f_a) \neq 0.
\]
Normalize $\ell_a$ so that $\ell_a(f_a) = 1$. Proposition~\ref{4.9} constructs the global normalized solution with
\[
(s, J)( \cdot , 0) = (0, \varepsilon f_a)
\]
and bounds \eqref{4.10}, \eqref{4.11}. Write
\[
J_\tau = \mathcal L_a J + \mathcal N_a,
\]
where
\[
\mathcal N_a = - (\delta_\ell y + a \mathcal P s) J'
+ \delta_A J + s \mathcal H J - J \mathcal H s + a B s'.
\]
As in Lemma~\ref{3.19}, define
\[
\ell_a(\mathcal N_a)
= \lim_{R \to \infty} \int \chi (y/R) \zeta_a \mathcal N_a.
\]
The same integration by parts is justified by Lemma~\ref{4.12} and bounds the adjoint pairing.
\[
| \ell_a(\mathcal N_a)|
\leq C_a \left (\| s \|_{L^2(\phi)} + \| s \|_{H^1} + \| J \|_{H^2}^2 \right) \| J \|_X
\leq C_a \varepsilon^3 e^{- 27 \tau /40}.
\]
For the even cutoff $\chi_R(y) = \chi (y/R)$, with $\chi = 1$ on $[- 1, 1]$ and supported in $[- 2, 2]$, the transport cutoff error is bounded, on each $0 \leq \tau \leq \tau_0 < \infty$, by
\[
\left | \int \chi_R'(V_a + \delta_\ell y + a \mathcal P s) \zeta_a J \right |
\leq C_{a, \tau_0} \| \zeta_a \|_{L^2(|y| \geq R)} \| J \|_2
\longrightarrow 0.
\]
The cutoff errors vanish as $R \to \infty$, so the projection satisfies
\[
\frac{d}{d \tau} \ell_a(J) = \lambda_a \ell_a(J) + \ell_a(\mathcal N_a),
\qquad \ell_a(J( \cdot , 0)) = \varepsilon.
\]
Since $27/40 + \operatorname{Re} \lambda_a > 0$, the variation-of-constants integral is bounded uniformly in time.
\[
\left |e^{- \lambda_a \tau} \ell_a(J(\tau)) - \varepsilon \right | \leq C_a \varepsilon^3 \int_0^\tau e^{- (27/40 + \operatorname{Re} \lambda_a) r} \, d r \leq \frac{C_a \varepsilon^3}{27/40 + \operatorname{Re} \lambda_a}.
\]
Choose $\varepsilon$ so that
\[
\frac{C_a \varepsilon^2}{27/40 + \operatorname{Re} \lambda_a} < \frac12.
\]
Then \eqref{4.13} implies
\[
| \ell_a(J(\tau))| \geq \frac12 \varepsilon e^{\operatorname{Re} \lambda_a \tau},
\qquad
\| \mathcal P J(\tau) \|_\infty \geq c_a \varepsilon e^{\operatorname{Re} \lambda_a \tau}.
\]

\emph{Modulation and physical blowup.} By \eqref{4.11},
\[
| \delta_A| + | \delta_\ell | \leq C_a \varepsilon^2 e^{- 9 \tau /20}.
\]
Hence the limiting factors
\[
A_\infty = \exp \int_0^\infty \delta_A(r) \, d r,
\qquad
L_\infty = \exp \left (- \int_0^\infty \delta_\ell (r) \, d r \right)
\]
satisfy
\[
e^{- C_a \varepsilon^2} \leq A_\infty, L_\infty \leq e^{C_a \varepsilon^2}.
\]
In particular,
\[
\sup_{0 < \varepsilon \leq \varepsilon_0(a)}
\max \{ A_\infty, A_\infty^{- 1}, L_\infty, L_\infty^{- 1} \}
\leq e^{C_a \varepsilon_0(a)^2}.
\]
With $\alpha_a = - \bar c_A > 0$, the scales determined by \eqref{1.9} satisfy
\[
\begin{aligned}
A(\tau)  & = A_\infty e^{- \alpha_a \tau}
\exp \left (- \int_\tau^\infty \delta_A(r) \, d r \right)
\sim A_\infty e^{- \alpha_a \tau}, \\
L(\tau)  & = L_\infty e^{- \bar c_\ell \tau}
\exp \left (\int_\tau^\infty \delta_\ell (r) \, d r \right)
\sim L_\infty e^{- \bar c_\ell \tau}.
\end{aligned}
\]
Consequently,
\[
T = \int_0^\infty A(r) \, d r < \infty,
\qquad
T - t(\tau) \sim \frac{A_\infty}{\alpha_a} e^{- \alpha_a \tau}.
\]
Restore tildes for normalized fields. Lemma~\ref{4.1} ensures $U_a \neq 0$, and
\[
\| \widetilde u - U_a \|_\infty
= \| \mathcal P s \|_\infty \leq C_a \varepsilon^2 e^{- 9 \tau /20}.
\]
Using
\[
(u, B)(x, t(\tau))
= \frac L A(\widetilde u, \widetilde B)(x/L, \tau),
\]
we obtain
\[
\| u(t(\tau)) \|_\infty
\sim \frac{L_\infty}{A_\infty} \| U_a \|_\infty
e^{(\alpha_a - \bar c_\ell) \tau}, \qquad \| B(t(\tau)) \|_\infty
\geq c_a \varepsilon
e^{(\alpha_a - \bar c_\ell + \operatorname{Re} \lambda_a) \tau}.
\]
The preceding time asymptotic proves \eqref{1.7}, with the exponents stated in the theorem. The uniform bounds on $A_\infty, L_\infty$ make these constants independent of $\varepsilon$.

\emph{Finite energy before blowup.} At every finite $\tau$, the normalized profiles belong to $Y_4$, so $y \widetilde \omega, y \widetilde J \in L^2$. The odd-profile Hardy estimate from the proof of Theorem~\ref{1.4} then bounds the physical energy after rescaling.
\[
\| u(t(\tau)) \|_2^2 + \| B(t(\tau)) \|_2^2
\leq \frac{4 L^3}{A^2}
\left (\| y \widetilde \omega \|_2^2 + \| y \widetilde J \|_2^2 \right) < \infty.
\]
\end{proof}

\end{document}